\documentclass[10pt,a4paper]{amsart}
\usepackage{amsfonts,layout}
\usepackage{bm,epsfig}
\usepackage{graphicx,tikz,subfig,color}

\newcount\minute
\newcount\hour
\def\currenttime{%
	\minute\time
	\hour\minute
	\divide\hour60
	\the\hour:\multiply\hour60\advance\minute-\hour\the\minute}
\def\draftnote{{\it \today \quad  \currenttime \hfill  tex-file :   \jobname}}

\catcode`@=11
\@addtoreset{equation}{section}

\catcode`@=12
\newtheorem{Theorem}{Theorem}[section]

\newtheorem{Proposition}{Proposition}[section]
\newtheorem{Lemma}{Lemma}[section]

\newtheorem{Corollary}{Corollary}[section]

\newtheorem{Remark}{Remark}[section]

\newtheorem{Hyp.}{Hyp.}[section]

\newcommand{\RR}{\mathbb{R}}

\newcommand{\eps}{\varepsilon}

\usetikzlibrary{decorations.markings}
\tikzset{->-/.style={decoration={markings,mark=at position .6 with {\arrow[scale=1.5]{stealth}}},postaction={decorate}}}

\begin{document}

\title[]{
Analysis of a two-layer energy balance model: long time behaviour and greenhouse effect}

\author{P. Cannarsa} 
\address{Dipartimento di Matematica, Universit\`a di Roma "Tor Vergata",
Via della Ricerca Scientifica, 00133 Roma, Italy}
\email{cannarsa@axp.mat.uniroma2.it}

\author{V. Lucarini} 
\address{Department of Mathematics and Statistics, University of Reading,
Reading RG66AX, United Kingdom,
Centre for the Mathematics of Planet Earth, University of Reading,
Reading RG66AX, United Kingdom}
\email{v.lucarini@reading.ac.uk}

\author{P. Martinez} 
\address{Institut de Math\'ematiques de Toulouse; UMR 5219, Universit\'e de Toulouse; CNRS \\ 
UPS IMT F-31062 Toulouse Cedex 9, France} \email{patrick.martinez@math.univ-toulouse.fr}

\author{C. Urbani} 
\address{Dipartimento di Scienze Tecnologiche e dell'Innovazione, Universitas Mercatorum,
Piazza Mattei 10, 00186 Roma, Italy}
\email{cristina.urbani@unimercatorum.it}

\author{J. Vancostenoble} 
\address{Institut de Math\'ematiques de Toulouse; UMR 5219, Universit\'e de Toulouse; CNRS \\ 
UPS IMT F-31062 Toulouse Cedex 9, France} \email{judith.vancostenoble@math.univ-toulouse.fr}

\thanks{P. Cannarsa and C. Urbani were partly supported by Istituto Nazionale di Alta Matematica (GNAMPA Research Projects) and by the MIUR Excellence Department Project awarded to the Department of Mathematics, University of Rome Tor Vergata, CUP E83C18000100006. C. Urbani also acknowledges Accademia Italiana dei Lincei - Postdoc scholarship "Beniamino Segre". V. Lucarini acknowledges the support provided by the Horizon 2020 project TiPES (Grant No. 820970), by the Marie Curie ITN CriticalEarth (Grant No. 956170), and by the EPSRC project EP/T018178/1. P. Martinez and J. Vancostenoble were supported by the Agence Nationale de la Recherche, Project TRECOS, under grant ANR-20-CE40-0009. This work was partly supported in the framework of LIA COPDESC}

\begin{abstract}
We study a two-layer energy balance model, that allows for vertical exchanges between a surface layer and the atmosphere. 
The evolution equations of the surface temperature and the atmospheric temperature are coupled by the emission of infrared radiation by one level, that emission being partly captured by the other layer, and the effect of all non-radiative vertical exchanges of energy.
Therefore, an essential parameter is the absorptivity of the atmosphere, denoted $\varepsilon _a$. The value of $\varepsilon _a$ depends critically on greenhouse gases: increasing concentrations of $\text{CO}_2$ and $\text{CH}_4$ lead to a more opaque atmosphere with higher values of $\epsilon_a$. First we prove that global existence of solutions of the system holds if and only if $\varepsilon _a \in (0,2)$, and blow up in finite time occurs if $\varepsilon _a >2$. (Note that the physical range of values for $\varepsilon _a$ is $(0,1]$). Next, we explain the long time dynamics for $\varepsilon _a \in (0,2)$, and we prove that all solutions converge to some equilibrium point. Finally, motivated by the physical context, we study the dependence of the equilibrium points with respect to the involved parameters, and we prove in particular that the surface temperature increases monotonically with respect to $\varepsilon _a$. 
This is the key mathematical manifestation of the greenhouse effect.
\end{abstract}

\maketitle



\section{Introduction}

\subsection{Energy Balance Models} \hfill

The climate is a multiphase  system featuring variability over many temporal and spatial scales. Its evolution can be written in terms of extremely complicated conservation laws for energy, momentum, and chemical species for three-dimensional (3D) fields \cite{Peixoto1992, Lucarini2014}. Given such a level of complexity, it is far from trivial to relate data, theories, and numerical models \cite{Held}. Indeed, the theoretical and numerical investigation of the climate system relies on the use of models that differ wildly in terms of scope, details, and overall complexity, ranging from extremely low dimensional models to Earth system models, which are some of the heaviest users of high performance computing facilities \cite{GL2020}. 

A simple yet extremely valuable approach to the study of the climate system comes from the use of  Energy Balance Models (EBMs), which had originally been introduced in the sixties independently by  Budyko \cite{Budyko} and Sellers \cite{sellers}. Such models describe in a very simplified yet effective way the evolution of the zonally averaged temperature on the Earth's surface, thus reducing the problem to a single 1D field. The planet receives radiation from the Sun (mostly in the form of visible and ultraviolet radiation); part of this radiation is scattered back to space through an elastic process where no energy is exchanged, and part is absorbed, mostly at surface. Then, radiation is emitted by the planet, mostly in the form of infrared radiation. The incoming solar radiation is unequally distributed over the surface of the planet, hence the balance between absorbed and emitted radiation will depend on latitude. A variety of physical processes, mostly associated with the large-scale motion of the geophysical fluids (the atmosphere and the ocean) are responsible for transporting heat from warm to cold regions, thus acting effectively as agents of diffusion. An EBM evolution equation reads as
\begin{equation}\label{1DEBM}
\gamma \Bigl[ \frac{\partial T}{\partial t} -  k \frac{\partial }{\partial x} \Bigl( (1-x^2) \frac{\partial T}{\partial x}\Bigr) \Bigr] = \mathcal R_s -\mathcal{R}_e,
\end{equation}
where $T(t,x)$ is the surface temperature, measured in Kelvin degrees, at colatitude $\theta=\sin^{-1}x$, $x\in(-1,1)$ is the space variable, $t>0$ is the time variable, $\mathcal{R}_s$ and $\mathcal{R}_e$ are the average amount of solar energy flowing into and out a unit area of the Earth surface per unit time. The constant $\gamma$ represents the effective heat capacity (which is the energy needed to raise the temperature by one kelvin), while the quantity $k\gamma=D$ is the effective thermal conductivity, which controls the efficacy of the latitudinal diffusion of energy. As hinted above, this is a very simplified way to represent the effect of the action of the geophysical fluids in the climate system.

The fundamental laws of thermodynamics impose that the amount of energy radiated from Earth to space depends on the temperature. As a first approximation, we can assume that the Earth emits as a black body with a surface temperature $T$. Therefore, we assume that function $\mathcal{R}_e$ follows the Stefan-Boltzmann law
\begin{equation*}
\mathcal{R}_e(T)=\sigma_BT^4
\end{equation*}
where $\sigma_B=5.67\cdot 10^{-8}Wm^{-2}K^{-4}$ is the Stefan-Boltzmann constant.

The energy absorbed by the Earth is a fraction of the incoming solar radiation $Q$
\begin{equation*}
\mathcal{R}_s=Q(t,x) \beta (T),
\end{equation*}
where $\beta$ is the effective coalbedo. The effective coalbedo 
depends on many local factors as cloud cover, composition of the Earth's atmosphere, presence of ice on the Earth's surface, etc., and, by and large, has to do with the color of the planet as seen from space: darker hues absorb more solar radiation than light ones. It is possible to provide a reasonable parametrization of the coalbedo as a function of the temperature via piecewice linear function of the form
\begin{equation*}
\beta(T)=\begin{cases}
\beta_-&T\leq T_-\\
\beta_-+\frac{T-T_-}{T_+-T_-}(\beta_+-\beta_-)&T\in[T_-,T_+]\\
\beta_+&T\geq T_+.
\end{cases}
\end{equation*}
Indeed, the polar regions, where the temperature is lower, can be covered by snow and ice and have a higher cloud cover, leading to a smaller coalbedo with respect to equatorial region, which are free of snow and ice and covered with land and open water.  Typical reference values for the parameters of the equation above are
\begin{equation*}
\beta_-=0.3,\quad \beta_+=0.7,\quad T_-=250\,K,\quad T_+=280\,K,
\end{equation*}
see, for instance, \cite[Chapter 2]{KE}.

The solar radiation $Q$ can be taken of the form
\begin{equation}\label{Q}
Q(t,x)=r(t)q(x)
\end{equation}
where $r(t)$ is a positive, possibly periodic, function allowing for seasonal cycle and $q(x)$ is the latitudinal-dependent insolation function, which depends on the geometry of the Sun-Earth system \cite{Peixoto1992}.
\subsection{Multistability and Critical Transitions}\hfill

Many authors studied the well-posedness, uniqueness of solutions, asymptotic behaviour, existence of periodic solutions, free boundary problem and numerical approximations of these models. We recall the results of North and co-workers \cite{North75,North81,north}, Ghil \cite{Ghil76}, Held and Suarez \cite{Held-Suarez}, Diaz and co-authors \cite{Diaz-Hetzer,DHT,Diaz-Tello}, Hetzer \cite{Hetzer96globalexistence, Hetzer1, hetzer2011global}, and many others. Chen and Ghil \cite{Chen1996} studied in detail a more sophisticated version of the problem above, comprising of an atmosphere described by an EBM coupled to a ocean described through (approximate) fluid dynamical equations, finding low-frequency variability associated with the occurrence of a Hopf bifurcation.

We remark that the EBM described above, despite its simplicity, has been instrumental for discovering the multistability of the Earth's climate. Indeed, as anticipated by Budyko \cite{Budyko} and Sellers \cite{sellers} and analysed in detail by Ghil \cite{Ghil76}, the model allows for the presence of two competing asymptotic states for the same values of the parameters. Such states correspond to the current warm climate and the so-called snowball state, characterised by global glaciation and surface temperatures of the order of 220 K. Paleoclimatic evidences collected in the '90s \cite{Hoffman2000} have shown that, indeed, our planet has spent in the distant past many million years in snowball conditions, the departure from which has allowed the evolution of multicellular life \cite{Pierrehumber2011,GL2020}. Between the two competing climates, one can find a saddle solution which lives on an invariant set that belongs to the boundary between the two competing basins of attraction; see discussion in Lucarini and Bodai and references therein \cite{Lucarini-Bodai2017, Lucarini-Bodai2019, Lucarini-Bodai2020}. 

We remark that one of the key manifestations of the multistability of the climate  is the existence - in models and observations - of critical transitions, which lead the system to qualitative (and \textit{de facto} irreversible) changeovers from one regime of operation to a qualitatively different one \cite{Feudel}. A paradigmatic example of such critical transition is the saddle-node bifurcation whereby the warm state and the snowball states becomes the only viable attractor  \cite{Ghil76}. 

In Earth system sciences such critical behaviour is associated with the so-called tipping points (TPs) \cite{Lenton}. Indeed, the history of the Earth's climate features periods of relatively smooth response to perturbations alternating with rapid changes due to TPs \cite{GL2020,PaleoJump,Boers2022}. We are now at risk of experiencing within our lifetime the collapse of the Amazon Forest \cite{Boulton2022} (forest to savannah transition) or of the Atlantic meridional overturning circulation \cite{Boers2021} (transition from vigorous to very weak circulation). The nearing of a tipping points is flagged by the increased sensitivity of a system to perturbations and by the increase in the correlation time of generic signals \cite{Lenton2012}; see a more complete theory in \cite{Chekroun2020,Santos2022} and extension to a time-dependent framework in \cite{Boettner2022}.

\subsection{From one to two Layers}\hfill

The EBM described in the equations above can be improved by increasing the vertical resolution. Indeed, considering various vertical layers it is possible to represent, at least approximately, the very important vertical exchanges processes occurring between surface and the atmosphere, and, possibly, between different atmospheric levels (e.g. troposphere and stratosphere) \cite{Hartmann}. Hence, instead of considering just one vertically homogeneous layer, a more accurate description of the climate system can be obtained via the following two-layer energy balance model (2LEBM), that allows for vertical exchanges between a surface layer and the atmosphere:

\begin{equation}
\label{2layer-pbm}
\begin{cases}
\gamma _a \Bigl[ \frac{\partial T_a}{\partial t} -  k_a \frac{\partial }{\partial x} \Bigl( (1-x^2) \frac{\partial T_a}{\partial x}\Bigr) \Bigr] \\
\quad\quad= -\lambda (T_a - T_s) + \varepsilon _a \sigma _B \vert T_s \vert ^3 T_s  - 2 \varepsilon _a \sigma _B \vert T_a \vert ^3 T_a + \mathcal R_a , \\
\gamma _s \Bigl[ \frac{\partial T_s}{\partial t} -  k_0 \frac{\partial }{\partial x} \Bigl( (1-x^2) \frac{\partial T_s}{\partial x}\Bigr) \Bigr] \\
\quad\quad= -\lambda (T_s - T_a) - \sigma _B \vert T_s \vert ^3 T_s  + \varepsilon _a \sigma _B \vert T_a \vert ^3 T_a + \mathcal R_s , \\
(1-x^2) \frac{\partial T_a}{\partial x} _{\vert x = \pm 1} = 0 = (1-x^2) \frac{\partial T_s}{\partial x} _{\vert x = \pm 1} ,\\
T_a (0,x)=T_a ^{(0)} (x) ,\\
T_s (0,x)=T_s ^{(0)} (x) .

\end{cases}
\end{equation}
$T_a$ represents the temperature of the atmospheric layer while $T_s$ stands for the surface temperature. The energy coupling between the two layers occurs through two different terms. One involves the emission of infrared radiation by one level, that emission being captured by the other layer. The other one is linear with the difference between the temperature of the two layers and describes succinctly the effect of all non-radiative vertical exchanges of energy due to the action of the geophysical fluids (see \cite[Chapter 10]{North-Kim}). Notice that when the atmospheric temperature is lower than the surface one, this term tends to warm up the atmosphere and cool down the surface. 
Note that the atmosphere is assumed to have, in general, non-unitary absorptivity $\epsilon_a$, because it is treated as a grey rather than black body. One needs to keep in mind that out of fundamental physical principles $0\leq\epsilon_a\leq1$. The value of $\epsilon_a$ depends critically on greenhouse gases: increasing concentrations of $\text{CO}_2$ and $\text{CH}_4$ lead to a more opaque atmosphere with higher values of $\epsilon_a$; see an instructive discussion in \cite{Hartmann}. Indeed, $\eps_a$ measures the greenhouse effect: an estimate of $\eps_a$ for a basic Energy Balance Model for present-day conditions gives $\eps_a\approx0.62$ \cite[Chapter 2]{KE}).

Similarly to the one-layer model, $\mathcal R_s$ is the solar radiation absorbed at the surface. It is a fraction of the incoming solar flux $Q$
\begin{equation}
\mathcal R_s (t) = Q(t,x) \, \beta _s (T_s),
\end{equation}
where $\beta _s$ is the coalbedo function. In general, $\beta _s$ is modelled as a nondecreasing positive and bounded function (as, for instance, the piecewise linear function showed for the one-layer case). We also introduce the term $\mathcal R_a$, which represents the solar radiation absorbed by the atmospheric layer
\begin{equation}
\mathcal R_a (t) = Q(t,x) \, \beta _a (T_a),
\end{equation}
and is much smaller than $\mathcal R_s$, because the atmosphere is almost transparent in the visible range. Note that most of such absorption occurs in the stratosphere, whereas the troposphere, the atmospheric layer that is closest to surface and that contains most of the mass, plays a lesser role, unless pollutants like black carbon are present. Indeed, as well known, the atmosphere is a system warmed \textit{from below}, because the external forcing coming from the absorption of the solar radiation acts prevalently at surface \cite{Peixoto1992}. We choose for the incoming flux $Q(t,x)$ the representation \eqref{Q}.
Finally, the generalized Neumann boundary condition arises naturally when one performs the change of variable
$\theta = \sin ^{-1} x$ between the colatitude $\theta$ and the new space variable $x$.


Note that related two slabs or two boxes models have been studied before in \cite{Gregory2000, Held2010, Stommell961}. However, in these works the authors considered coupled linear evolution equations with no diffusion, or diffusion only on one layer.

We also remark that atmospheres can be \textit{very} opaque to infrared radiation (much more than the Earth's) as a result of their composition and/or sheer mass thereof. The most obvious example is Venus, where the surface pressure is 90 times larger than the Earth's and the atmosphere is overwhelmingly composed of $\text{CO}_2$. The planet Venus is conjectured, in fact, to have undergone a runaway greenhouse transition in a now distant past \cite{Ingersoll,Kastings,DelGenio}. 

The runaway greenhouse effect emerges when the surface warming leads to excessive evaporation of surface water, with the resulting water vapour contributing to further increasing the opacity of the atmosphere, up to full evaporation of the available water, eventual loss of water vapour to space, and transition to a fundamentally different climate (a divergent behaviour, at all practical levels). In a much weaker form, the water vapour feedback contributes to a great amplification of the greenhouse effect on the Earth with respect to what would be realised in absence of water \cite{Peixoto1992}.


\subsection{Outline of the main findings} \hfill

In this paper we study a simplified version of the two-layer system given in \eqref{2layer-pbm} where we neglect the effect of latitudinal variation of the fields, so that the system of partial differential equations can be reduced to the following autonomous ODE problem:
\begin{equation}
\label{2layer-pbm-EDO}
\begin{cases}
\gamma _a T_a' = -\lambda (T_a - T_s) + \varepsilon _a \sigma _B \vert T_s \vert ^3 T_s - 2 \varepsilon _a \sigma _B \vert T_a \vert ^3 T_a + \mathcal R_a (T_a) , \\
\gamma _s T_s' = -\lambda (T_s - T_a) - \sigma _B \vert T_s \vert ^3 T_s  + \varepsilon _a \sigma _B \vert T_a \vert ^3 T_a + \mathcal R_s (T_s), \\
T_a (0)=T_a ^{(0)} ,\\
T_s (0)=T_s ^{(0)} .
\end{cases}
\end{equation}
It appears that the parameter $\varepsilon _a$ plays a major role in the qualitative behaviour of the solution of \eqref{2layer-pbm-EDO}. 
\begin{itemize} 
\item If $\varepsilon _a \in (0,2)$, we prove that the solution is global in time, remains positive and bounded, and converges to some equilibrium point (with positive components),
see Proposition \ref{prop-wellposed} in section \ref{sec-global}.
 
\item on the other hand, if $\varepsilon _a > 2$, we prove that there is blow up in finite time, at least for some solutions, see Proposition \ref{prop-blowup} in section \ref{sec-blowup}.
\end{itemize}

Next, when $\varepsilon _ a\in (0,2)$, in order to analyze the influence of the different parameters, we focus on the problem where $\mathcal R_a =0$ and $\mathcal R _s$ is piecewise linear: we assume that
there exist $q>0$, $T_{s,+} > T_{s,-} >0$ and $\beta _{s,+} > \beta _{s,-} >0$ such that $\mathcal R_s (T_s) = q \beta _s (T_s)$ where
\begin{equation}
\label{beta}
 \beta _s (T) =
\begin{cases}
 \beta _{s,-} & T \leq T_{s,-} , \\
 \beta _{s,-} + (\beta _{s,+} - \beta _{s,-}) \frac{T-T_{s,-}}{T_{s,+} - T_{s,-}} &T \in [T_{s,-}, T_{s,+}], \\
\beta _{s,+} &T \geq T_{s,+} .
\end{cases} 
\end{equation}
Then we prove the following results:
\begin{itemize}

\item For all $\lambda \geq 0$, all solutions converge to some equilibrium point, see Proposition \ref{prop-comp-lambda-eq-cv}. All the equilibrium points $(T_a,T_s)$ remain in a compact subset independent of $\lambda$, and satisfy $T_a< T_s$. Moreover, there is at most one warm equilibrium ($T_s > T_{s,+}$, where $T_{s,+}$ appears in \eqref{beta}),
and at most one cold equilibrium
($T_s < T_{s,-}$, where $T_{s,-}$ appears in \eqref{beta}), and these equilibrium points are asymptotically exponentially stable, see Proposition \ref{prop-comp-lambda-eq-gene}. Furthermore, there exists at most a finite number of equilibria (Proposition \ref{prop-comp-lambda-eq-gene}) and we are able to describe the asymptotic behaviour of the solutions of our problem.

Moreover, we have the following monotonicity properties: 
assume that $(T_a ^{(\varepsilon_a, \lambda)},T_s ^{(\varepsilon_a, \lambda)})$ is a warm or a cold equilibrium (hence assume that $T_s ^{(\varepsilon_a, \lambda)} \notin [T_{s,-}, T_{s,+}]$), then,
as a function of $\varepsilon _a$ and $\lambda$, we prove that the surface temperature $T_s ^{(\varepsilon_a, \lambda)}$ 
\begin{itemize}
\item increases as $\varepsilon_a$ grows (see Proposition \ref{prop-comp-lambda=0}),
\item and decreases as $\lambda$ grows (see Proposition \ref{prop-comp-lambda>0});
\end{itemize}
and the atmosphere temperature $T_a ^{(\varepsilon_a, \lambda)}$
\begin{itemize}
\item increases as $\varepsilon _a$ grows, at least if $\varepsilon _a \in [1,2)$ (see Proposition \ref{prop-comp-lambda=0}),
\item increases as $\lambda$ grows if $\varepsilon _a \in [0,1)$ 
and decreases as $\lambda$ grows if $\varepsilon _a \in (1,2)$ (see Proposition \ref{prop-comp-lambda>0}).
\end{itemize}
As a consequence, an increment of $\eps_a$ causes a rise of $T_s$ (see Corollary \ref{cor-greenhouse}).

\item When $\lambda =0$ (or $\lambda>0$ and $\eps_a\in(0,1]$) we can describe more precisely the stability of the system: there are exactly one, two or three equilibrium points, and we are able to detail the asymptotic behaviour of all solutions. Indeed, in this case, the number of equilibrium points is perfectly determined by the values of $\varepsilon _a$, $\sigma_B$, $q$, and the parameters  $T_{s,-}$, $T_{s,+}$, $\beta _{s,-}$, $\beta _{s,+}$ appearing in the structure of function $\beta_s$, see section \ref{sec-comp-asympt-lambda=0}. Moreover, if there are three equilibrium points (one warm, one cold and one intermediate), then we have proved that the intermediate equilibrium is unstable, i.e. we clarify the multistable nature of the climate system.

\end{itemize}


\medskip

This paper is structured as follows:

\begin{itemize}
\item In Section \ref{sec2}, we state the main results:
\begin{itemize}
\item global existence when $\varepsilon _a \in (0,2)$ (see Proposition \ref{prop-wellposed}),
\item blow up in finite time if $\varepsilon _a >2$ (see Proposition \ref{prop-blowup}),
\item asymptotic behaviour if $\varepsilon _a \in (0,2)$ (see Propositions \ref{prop-comp-lambda-eq-cv} and \ref{prop-comp-lambda-eq-gene}),
\item monotonicity of the equilibrium points with respect to parameters $\varepsilon _a$ and $\lambda$ (see Propositions \ref{prop-comp-lambda=0} and \ref{prop-comp-lambda>0}) and application to the asymptotic behaviour (Corollary \ref{cor-greenhouse}).
\end{itemize}

\item In Section \ref{sec3}, we mention some open problems.

\item In Section \ref{sec-global-proof}, we study the well-posedness  of the problem (existence, uniqueness, positivity of solutions).

\item In Section \ref{sec-comp-asympt-lambda=0}, we analyze the behaviour of all solutions when $\lambda=0$ and $\mathcal R_a=0$.

\item In Section \ref{sec7}, we consider the case of $\lambda>0$ and $\mathcal R_a=0$.

\item In Section \ref{sec8}, we derive several results showing the sensitivity of the equilibria with respect to parameters $\lambda$ and $\varepsilon_a$.

\item In Appendix \ref{sec-conv-Phi1} we derive  upper bounds for the number of equilibria which, in the physical case $\eps_a\in(0,1]$, is equal to three.

\item In Appendix \ref{sec-blowup-proof}, we show that solutions may blow up in finite time.

\end{itemize}


\section{Statement of the main results}
\label{sec2}


\subsection{Global existence, positivity and boundedness for $\varepsilon _a \in (0,2)$} 
\label{sec-global}\hfill

We make the following assumptions:
\begin{itemize}
\item let the coefficients $\lambda$, $q$ and $\sigma _B$ be such that
\begin{equation}
\label{global-coeffs}
\lambda \geq 0, \quad q>0, \quad \sigma _B >0 .
\end{equation}

\item  let $\beta_a, \beta _s: \mathbb R \to \mathbb R$ be globally Lipschitz continuous and such that $\beta _a \geq 0$ and $\beta _s >0$, and define
\begin{equation}
\label{global-coalb}
\mathcal R_a = q \beta _a (T_a) , \quad \mathcal R_s = q \beta _s (T_s).
\end{equation}

\item let the initial conditions satisfy
\begin{equation}
\label{global-ci}
T_a ^{(0)} \geq 0 , \quad  T_s ^{(0)} \geq 0.
\end{equation}

\item let $\varepsilon _a$ be in the following range
\begin{equation}
\label{global-eps}
\varepsilon _a \in (0,2).
\end{equation}
 \end{itemize}
\begin{Proposition}
\label{prop-wellposed}

Under the assumptions \eqref{global-coeffs}-\eqref{global-eps}, problem \eqref{2layer-pbm-EDO} admits a unique solution, which is defined and bounded for any $t\in[0,+\infty)$. Moreover 
\begin{equation}
\label{global-posi}
\forall\, t \in (0,+\infty), \quad T_a (t) > 0 \quad \text{ and } \quad T_s (t) >0.
\end{equation}
\end{Proposition}
The proof of Proposition \ref{prop-wellposed} is given in Section \ref{sec-global-proof}.


\subsection{Asymptotic behaviour of the solutions when $\varepsilon _a \in (0,2)$ and $\lambda \geq 0$}\hfill

When $\varepsilon _a \in (0,2)$, the solution of \eqref{2layer-pbm-EDO} is global in time. 
\begin{Proposition}
 \label{prop-comp-lambda-eq-cv}

Consider $\varepsilon _a \in (0,2)$ , $\lambda \geq 0$,
$T_s ^{(0)} \geq 0$ and $T_ s^{(0)} \geq 0$. Then, the solution of \eqref{2layer-pbm-EDO} converges to an equilibrium point.
 \end{Proposition}
 
Proposition \ref{prop-comp-lambda-eq-cv} follows from a general result concerning cooperative systems (see Smith \cite{Smith}). We complete this general convergence result in two directions:

\begin{itemize}
\item first, explaining how the convergence occurs when $\lambda =0$ and $\mathcal{R}_a = 0$ (see Section \ref{sec-comp-asympt-lambda=0},
where we give a complete description of the basins of attraction of the different equilibrium points for $\lambda =0$, and Section \ref{sec7} for$\lambda >0$),

\item finally, proving some properties of the equilibrium points in the general case $\lambda \geq 0$:
\end{itemize}

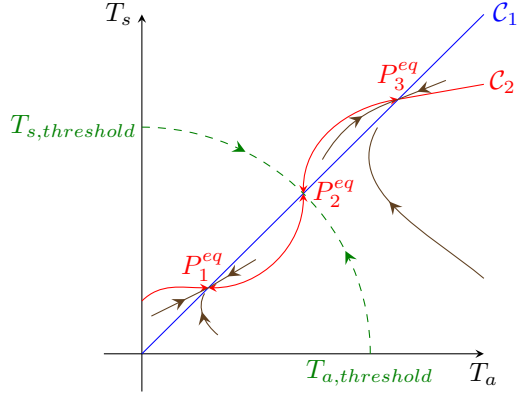
\begin{figure}[h!]
\centering

\begin{tikzpicture}[scale=0.5]
\draw[-stealth] (-1,0) -- (9,0)node[below]{$T_a$};
\draw[-stealth] (0,-1) -- (0,9)node[left]{$T_s$};
\draw[blue] (0,0) -- (9,9)node[right]{$\mathcal C_1$};
\draw[red,-stealth] (0,1.4) to[in=180] (1.75,1.75);
\node[red] at (1.6,2.3){$P_1 ^{eq}$};
\draw[red,stealth-stealth] (1.75,1.75) to[out=0,in=-90] (4.25,4.25)node[right]{$ P_2 ^{eq}$};
\draw[red,stealth-stealth] (4.25,4.25) to[out=90,in=190] (6.75,6.75)node[above]{$P_3 ^{eq}$};
\draw[red] (6.75,6.75) -- (9,{.176*(9-6.75)+6.75});
\draw[->-,brown!50!black] (.25,1) -- (1.75,1.75);
\draw[->-,brown!50!black] (3,2.5) -- (1.75,1.75);
\draw[->-,green!50!black,dashed] (0,{sqrt(2)*(4.25)}) arc (90:45:{sqrt(2)*(4.25)});
\draw[->-,green!50!black,dashed] ({sqrt(2)*(4.25)},0) arc (0:45:{sqrt(2)*(4.25)});
\draw[->-,brown!50!black] (4.75,5.15) to[out=60,in=200] (6.75,6.75);
\draw[->-,brown!50!black] (8,7.25) -- (6.75,6.75);
\draw[red] (9.5,7.2)node{$\mathcal C_2$};
\draw[->-,brown!50!black] (2,.5) to[out=140,in=240] (1.65,1.6);
\draw[->-,brown!50!black] (9,2) to[out=140,in=240] (6.2,6);
\draw[green!50!black] (6,-0.5)node{$T_{a,threshold}$};
\draw[green!50!black] (-1.8,6)node{$T_{s,threshold}$};
\end{tikzpicture}
\caption{In the phase space we sketch by black and green arrows the convergence of initial conditions according to their belonging to the basins of attraction of the three equilibria $(2^{-1/4}T^*_{s,1},T^*_{s,1})$, $(2^{-1/4}T^*_{s,2}, T^*_{s,2})$ and $(2^{-1/4}T^*_{s,3}, T^*_{s,3})$, solutions of (6.1).
}
\label{fig:basins}
\end{figure}

\begin{Proposition}
 \label{prop-comp-lambda-eq-gene}
Assume that $\varepsilon _a \in (0,2)$, $\lambda \geq 0$, $q>0$, $\beta _s >0$ is given by \eqref{beta}, and $\mathcal R_a = 0$. Then 
\begin{itemize}
\item problem \eqref{2layer-pbm-EDO} has at least one equilibrium point $(T_a,T_s)$ and at most a finite number of equilibria. In particular, there are at most three equilibrium points for $\lambda=0$ and any $\eps_a\in(0,2)$, and for any $\lambda>0$ and $\eps_a\in(0,1]$,
\item all equilibrium points $(T_a,T_s)$ of problem \eqref{2layer-pbm-EDO}
belong to a compact subset of $(0,+\infty)^2$ which is independent of $\lambda$,
\item problem \eqref{2layer-pbm-EDO} has at most one cold equilibrium $(T_a,T_s)$ (that is, $T_s \leq T_{s,-}$), and at most one warm equilibrium $(T_a,T_s)$ (namely, $T_s \geq T_{s,+}$),
\item a cold equilibrium point $(T_a,T_s)$ is asymptotically exponentially stable,
\item a warm equilibrium point $(T_a,T_s)$ is asymptotically exponentially stable.
\end{itemize}
\end{Proposition}
The asymptotic behaviour of the solution of our system, for the case $\lambda=0$, is summarized in the phase space shown in Figure \ref{fig:basins} below (see the proof in Section \ref{sec-mono-3eq} for $\lambda =0$ and in Section \ref{sec-asympt-lambda>0} for $\lambda >0$).
 

\subsection{The influence of the parameters $\varepsilon _a$ and $\lambda$ on the equilibrium points} \hfill

Since all solutions converge to equilibrium points, it is interesting to study the behaviour of such points with respect to the parameters $\varepsilon _a$ and $\lambda$. We consider the case where $\mathcal R_a =0$, and $\beta_s$ is given by \eqref{beta}. Our results are the following:

\begin{Proposition}
\label{prop-comp-lambda=0}
Fix $\lambda \geq 0$. 
Assume that $(T_a ^{eq,\varepsilon _a ^*}, T_s ^{eq,\varepsilon _a ^*})$
is a warm [respectively cold] equilibrium point of problem \eqref{2layer-pbm-EDO} with $\varepsilon _a = \varepsilon _a ^*$. 
More precisely, assume that $T_s ^{eq,\varepsilon _a ^*} \notin [T_{s,-}, T_{s,+}]$.

Then, there exists a unique warm [respectively cold] equilibrium point $(T_a ^{eq,\varepsilon _a }, T_s ^{eq,\varepsilon _a})$ of problem \eqref{2layer-pbm-EDO} for $\varepsilon _a$ close to $\varepsilon _a ^*$. This equilibrium is also asymptotically exponentially stable. Moreover, the functions $\varepsilon _a \mapsto T_s ^{eq,\varepsilon _a }$ and $\varepsilon _a \mapsto T_a ^{eq,\varepsilon _a }$ are locally analytic, and the following monotonicity properties hold:
\begin{itemize}
\item locally, function $\varepsilon _a \mapsto T_s ^{eq,\varepsilon _a }$ is increasing. Hence, the surface temperature equilibrium increases as $\varepsilon _a$ increases;
\item locally, function $\varepsilon _a \mapsto T_a ^{eq,\varepsilon _a }$ is increasing if $\varepsilon _a ^* \in (1,2)$. Therefore, for such range of $\varepsilon^*_a$, the atmosphere temperature equilibrium increases as $\varepsilon _a$ increases.
\end{itemize}
\end{Proposition}

Note that we were not able to determine the variations of $T_a ^{eq,\varepsilon _a }$ when $\varepsilon _a ^* \in (0,1)$. However, we have established not only the sign of the derivative of $T_s ^{eq,\varepsilon _a }$ with respect to $\varepsilon _a$, but also its value, which is interesting to predict the evolution of $T_s ^{eq,\varepsilon _a }$ with respect to $\varepsilon _a$. From the previous proposition we deduce the following
\begin{Corollary}
    \label{cor-greenhouse}
    Fix $\lambda \geq 0$. Let $(T_a ^{eq,\varepsilon _a ^*}, T_s ^{eq,\varepsilon _a ^*})$ be a warm equilibrium of problem \eqref{2layer-pbm-EDO} associated to the parameter $\varepsilon _a ^*$. Assume that $T_s ^{eq,\varepsilon _a ^*} > T_{s,+}$. Then
    \begin{itemize}
        \item for all $\varepsilon _a \in [\varepsilon _a ^*, 2)$, there exists a unique warm equilibrium of problem \eqref{2layer-pbm-EDO} with parameter $\varepsilon _a$. Furthermore, function $\varepsilon _a \mapsto T_s ^{eq,\varepsilon _a }$ is analytic and increasing on $[\varepsilon _a ^*, 2)$ and function $\varepsilon _a \mapsto T_a ^{eq,\varepsilon _a }$ is analytic on $[\varepsilon _a ^*, 2)$;
        \item  given $\varepsilon _a ^+ > \varepsilon _a ^*$, the solution of problem \eqref{2layer-pbm-EDO} with parameter $\varepsilon _a ^+$ and initial conditions $(T_a ^{eq,\varepsilon _a ^*}, T_s ^{eq,\varepsilon _a ^*})$ converge to the warm equilibrium $(T_a ^{eq,\varepsilon _a ^+}, T_s ^{eq,\varepsilon _a ^+})$.
    \end{itemize}
\end{Corollary}
We can interpret the second item of Corollary \ref{cor-greenhouse} as follows: if the absorptivity parameter increases, jumping from $\varepsilon _a ^*$ to $\varepsilon _a ^+$, then
the former warm equilibrium $(T_a ^{eq,\varepsilon _a ^*}, T_s ^{eq,\varepsilon _a ^*})$, which is no more an equilibrium, converges to the new warm equilibrium, $(T_a ^{eq,\varepsilon _a ^+}, T_s ^{eq,\varepsilon _a ^+})$, associated to $\eps_a^+$. Furthermore, since $T_s ^{eq,\varepsilon _a ^+} > T_s ^{eq,\varepsilon _a ^*}$, the surface temperature rises as the pollution parameter increases.


Let us now analyse the dependence of the equilibria on the coupling parameter $\lambda$.
\begin{Proposition}
\label{prop-comp-lambda>0}
Fix $\varepsilon _a \in (0,2)$.
Assume that $(T_a ^{eq,\lambda ^*}, T_s ^{eq,\lambda ^*})$
is a warm [respectively cold] equilibrium point of problem \eqref{2layer-pbm-EDO} with $\lambda = \lambda ^* \geq 0$.
More precisely, assume that $T_s ^{eq,\lambda ^*} \notin [T_{s,-}, T_{s,+}]$.

Then, there exists a unique equilibrium point $(T_a ^{eq,\lambda }, T_s ^{eq,\lambda })$ of problem \eqref{2layer-pbm-EDO} for $\lambda$ close to $\lambda ^*$. Such equilibrium is also asymptotically exponentially stable. Furthermore, the functions $\lambda \mapsto T_s ^{eq,\lambda }$ and $\lambda \mapsto T_a ^{eq,\lambda }$ are locally analytic, and the following monotonicity properties hold:
\begin{itemize}
\item locally, function $\lambda \mapsto T_s ^{eq,\lambda }$ is decreasing. Thus, the surface temperature equilibrium decreases as $\lambda$ increases;
\item locally, function $\lambda \mapsto T_a ^{eq,\lambda }$ is
\begin{itemize}
\item  increasing if $\varepsilon _a \in (0,1)$,
\item  decreasing if $\varepsilon _a \in (1,2)$.
\end{itemize}
Hence, the atmosphere temperature equilibrium behaves monotonically with respect to $\lambda$, and the monotonicity depends on $\varepsilon _a$.
\end{itemize}
\end{Proposition}
Observe that we have determined not only the sign of the derivative of $T_s ^{eq,\lambda }$ and $T_a ^{eq,\lambda }$ with respect to $\lambda$, but also their values. This information can be useful to predict the evolution of $T_s ^{eq,\lambda }$ and $T_a ^{eq,\lambda }$ with respect to $\lambda$.


\subsection{Blow up in finite time when $\varepsilon _a >2$}
\label{sec-blowup}\hfill

We complete the results of Proposition \ref{prop-comp-lambda-eq-cv} on the global existence and boundedness of solutions of \eqref{2layer-pbm-EDO} for $\eps_\alpha\in(0,2)$ by studying the case of $\eps_\alpha>2$, which, as discussed above, has a mathematical motivation but a less solid physical interpretation.  

We prove the following result about the ODE system \eqref{2layer-pbm-EDO}. When $\varepsilon _a >2$, we give 
\begin{itemize} 
\item a precise result on a simplified model (assuming that $\lambda =0$ and $\mathcal R_a=0$): all solutions blow up in finite time.
\item a general result directly on \eqref{2layer-pbm-EDO} (assuming only that $\lambda \geq 0$): there exist some solutions that blow up in finite time.
\end{itemize}

Let us consider the following assumptions:
\begin{itemize}
\item let $\lambda$, $q$ and $\sigma _B$ satisfy
\begin{equation}
\label{global-coeffs-bl}
\lambda = 0, \quad q>0, \quad \sigma _B >0 ,
\end{equation}

\item let $\beta _a = 0$ and $\beta _s >0$ be defined by \eqref{beta}, so that
\begin{equation}
\label{global-coalb-bl}
\mathcal R_a = q \beta _a (T_a)=0 , \quad \mathcal R_s = q \beta _s (T_s),
\end{equation}

\item let $\varepsilon _a$ be such that
\begin{equation}
\label{global-eps-bl}
\varepsilon _a > 2.
\end{equation}
 \end{itemize}
\begin{Proposition}\hfill
\label{prop-blowup}

a) Under assumptions \eqref{global-ci},\eqref{global-coeffs-bl}-\eqref{global-eps-bl}, problem \eqref{2layer-pbm-EDO} admits a unique maximal solution that blows up in finite time.

b) Under the following weaker assumptions 
\begin{itemize}
\item let $\lambda >0$, $q>0$, $\sigma _B >0$,
\item let $\beta _a \geq 0$ and $\beta _s >0$ be globally Lipschitz continuous,
\item let $\varepsilon _a >2$,
\end{itemize}
and \eqref{global-ci}, there exist solutions of problem \eqref{2layer-pbm-EDO} that blow up in finite time.

\end{Proposition}
The proof of Proposition \ref{prop-blowup} is postponed in Appendix \ref{sec-blowup-proof}.


\section{Extensions and open questions}
\label{sec3}

\subsection{Some open questions directly related to our results}\hfill

\begin{figure}[h!]
    \centering
    \begin{tikzpicture}
        \node at (0,0){\includegraphics[width=.3\textwidth]{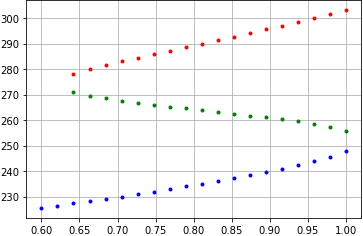}};
        \node at (-2.8,.2)[right]{$\scriptsize{T_s^{eq}}$};
        \node at (0,-1.5)[below]{$\scriptsize{\eps_a}$};
        \node at (6,0){\includegraphics[width=.3\textwidth]{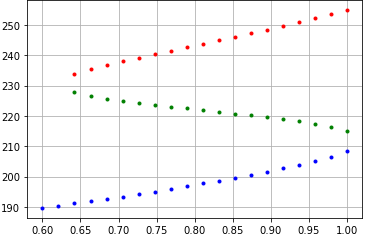}};
        \node at (3.2,.2)[right]{$\scriptsize{T_a^{eq}}$};
        \node at (6,-1.5)[below]{$\scriptsize{\eps_a}$};
    \end{tikzpicture}
    \caption{Sensitivity of equilibrium $T_s^{eq}$ (left)  and  $T_a^{eq}$ (right)  with respect to $\varepsilon_a$ in the case $\lambda=0$}
    \label{diag-lambda0}
\end{figure}
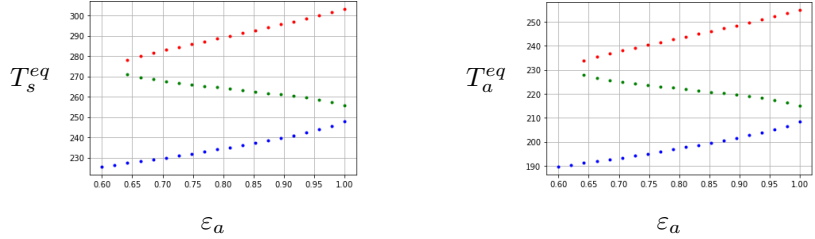

There are some properties that we were not able to prove which are of interest:
\begin{itemize}
\item the monotonicity of $\varepsilon_a \mapsto T_a^{eq, \varepsilon_a}$, where $T_a^{eq, \varepsilon_a}$ is a warm equilibrium, in the case $\varepsilon_a^\star \in [0,1)$ (see Proposition \ref{prop-comp-lambda=0}). From numerical computations (see right columns of  Figures \ref{diag-lambda0}, \ref{diag-lambda1}, \ref{diag-lambda2}), we expect that $\varepsilon_a \mapsto T_a^{eq, \varepsilon_a}$  is also increasing in these cases. 

\item we obtained in Proposition \ref{prop-comp-lambda=0} {\sl local} results  concerning the behaviour of warm (or cold) equilibria, and {\sl global} ones in Corollary \ref{cor-greenhouse} concerning the behaviour of the solution $(T_a, T_s)$ when $\varepsilon _a$ jumps to a {\it higher} value (that is, if suddenly there is an increase of concentration of $\text{CO}_2$). However, it would be interesting to go further and, for instance, to understand the behaviour of $(T_a,T_s)$ starting from a warm equilibrium point in the case of a {\sl lowering} of the value of $\varepsilon_a$.  Our numerical tests suggest that there is some hysteresis phenomenon, with the existence of a tipping point for which the solution of our system could jump from a warm equilibrium to a cold one. An analytic proof of such phenomenon would be of great concern. 
\end{itemize}

\begin{figure}[h!]
 \centering
    \begin{tikzpicture}
        \node at (0,0){\includegraphics[width=.3\textwidth]{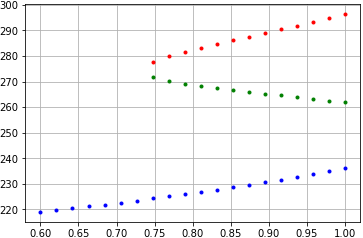}};
        \node at (-2.8,.2)[right]{$\scriptsize{T_s^{eq}}$};
        \node at (0,-1.5)[below]{$\scriptsize{\eps_a}$};
        \node at (6,0){\includegraphics[width=.3\textwidth]{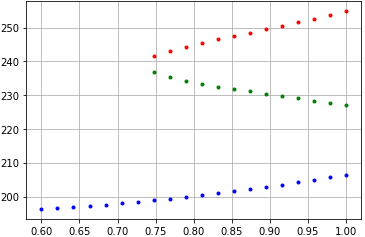}};
        \node at (3.2,.2)[right]{$\scriptsize{T_a^{eq}}$};
        \node at (6,-1.5)[below]{$\scriptsize{\eps_a}$};
    \end{tikzpicture}
\caption{Sensitivity of equilibrium $T_s^{eq}$ (left) and  $T_a^{eq}$ (right) with respect to $\varepsilon_a$ in the case $\lambda=1$}\label{diag-lambda1}
\end{figure}
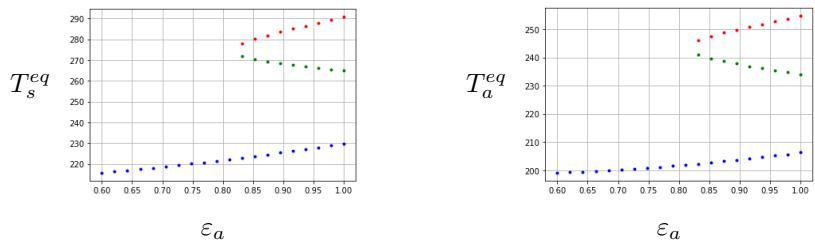
\begin{figure}[h!]
 \centering
    \begin{tikzpicture}
        \node at (0,0){\includegraphics[width=.3\textwidth]{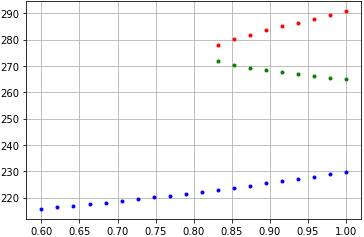}};
        \node at (-2.8,.2)[right]{$\scriptsize{T_s^{eq}}$};
        \node at (0,-1.5)[below]{$\scriptsize{\eps_a}$};
        \node at (6,0){\includegraphics[width=.3\textwidth]{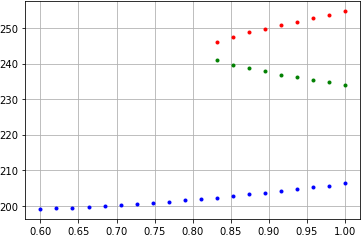}};
        \node at (3.2,.2)[right]{$\scriptsize{T_a^{eq}}$};
        \node at (6,-1.5)[below]{$\scriptsize{\eps_a}$};
    \end{tikzpicture}
\caption{Sensitivity of equilibrium $T_s^{eq}$ (left) and  $T_a^{eq}$ (right) with respect to $\varepsilon_a$ in the case $\lambda=2$}\label{diag-lambda2}
\end{figure}

\subsection{A periodic extension} \hfill

It would be natural to investigate problem \eqref{2layer-pbm-EDO} in presence of a seasonal effect:
$$ \mathcal R_a = r(t) q \beta _a (T_a), \quad \mathcal R_s = r(t) q \beta _s (T_s),$$
where the function $r$ is positive and periodic in time.
We expect the equilibrium points to be replaced by cycles, but a careful analysis of this model should be carried out. 

\subsection{Inverse problem question} \hfill

It would also be interesting to investigate the problem of determining the values of $\varepsilon _a$ and $\lambda$, as well as the insolation parameter $q$ (appearing in the expression of $\mathcal R_s$), from (the fewest possible) measurements of the solution. We refer to \cite{C-R-bis} for the determination of two coefficients in a reaction-diffusion equation (invasion model), and to \cite{Roques20140349, Sellers-manifold, memorySB} concerning Budyko-Sellers parabolic equation (possibly involving memory effects).


\section{Well-posedness}
\label{sec-global-proof}

In the following we prove Propositions \ref{prop-wellposed} and \ref{prop-comp-lambda-eq-cv}.

\subsection{Local existence} \hfill
\label{sec-local-existence}

Let us introduce the function $F: \mathbb R^2 \to \mathbb R^2$ defined by
\begin{equation}
\label{def-Fnonlin}
F \begin{pmatrix} T_a \\ T_s \end{pmatrix}= \begin{pmatrix} \frac{1}{\gamma _a} \Bigl[ -\lambda (T_a - T_s) + \varepsilon _a \sigma _B \vert T_s \vert ^3 T_s  - 2 \varepsilon _a \sigma _B \vert T_a \vert ^3 T_a + \mathcal R_a \Bigr] \\ \frac{1}{\gamma _s} \Bigl[ -\lambda (T_s - T_a) - \sigma _B \vert T_s \vert ^3 T_s  + \varepsilon _a \sigma _B \vert T_a \vert ^3 T_a + \mathcal R_s \Bigr] \end{pmatrix} .
\end{equation}
So, problem \eqref{2layer-pbm-EDO} can be recast into the form
$$ \begin{cases}
\begin{pmatrix} T_a \\ T_s \end{pmatrix} ' = F \begin{pmatrix} T_a \\ T_s \end{pmatrix}  , \\
\begin{pmatrix} T_a \\ T_s \end{pmatrix} (0) = \begin{pmatrix} T_a ^{(0)} \\ T_s ^{(0)} \end{pmatrix} .
\end{cases} $$
Then, the Cauchy-Lipschitz theorem implies that there exists a unique maximal solution $(T_a,T_s)$ defined on the maximal time interval $ (\tau _{a,s} ^-, \tau _{a,s} ^+)$, where $-\infty \leq \tau _{a,s} ^- < 0 < \tau _{a,s} ^+ \leq + \infty$.

\subsection{Positivity of the solution}\hfill

Now, since we are investigating positive initial conditions, let us prove the following result.

\begin{Lemma}
\label{lem-postivity}
Assume that $T_a ^{(0)} \geq 0$ and $T_s ^{(0)} \geq 0$. Then, as long as the solution $(T_a,T_s)$ exists, we have $ T_a (t) >0$ and $T_s(t)>0$ for positive times.
\end{Lemma}

\noindent {\it Proof of Lemma \ref{lem-postivity}.} Let us first consider initial conditions $(T^{(0)}_a,T^{(0)}_s)$ lying on the boundary of the set
$$ \mathcal Q := \Bigl\{ \begin{pmatrix} T_a \\ T_s \end{pmatrix} , T_a > 0, T_s > 0 \Bigr\} ,$$
and let us study the behaviour of $F$.

In the case $T^{(0)}_a>0$ and $T^{(0)}_s=0$, we have
$$ F \begin{pmatrix} T_a >0 \\ T_s =0 \end{pmatrix}
= \begin{pmatrix} \frac{1}{\gamma _a} \Bigl[ -\lambda T_a  - 2 \varepsilon _a \sigma _B T_a ^4 + \mathcal R_a \Bigr] \\ \frac{1}{\gamma _s} \Bigl[ \lambda T_a + \varepsilon _a \sigma _B T_a ^4 + \mathcal R_s \Bigr] \end{pmatrix} .$$
Observe that the second component is positive, and so $F$ points inward the set $\mathcal Q$. Therefore, the associated solution stays in $\mathcal Q$ for at least a certain period of time (small enough).

If $T^{(0)}_a=0$ and $T^{(0)}_s>0$, we get
$$ F \begin{pmatrix} T_a =0 \\ T_s >0 \end{pmatrix}
= \begin{pmatrix} \frac{1}{\gamma _a} \Bigl[ \lambda T_s + \varepsilon _a \sigma _B T_s ^4  + \mathcal R_a \Bigr] \\ \frac{1}{\gamma _s} \Bigl[ -\lambda T_s - \sigma _B  T_s ^4  + \mathcal R_s \Bigr] \end{pmatrix} .$$
In this case the first component is positive and this implies that $F$ points inward the set $\mathcal Q$. Once again, we find that the associated solution lives in $\mathcal Q$ for a certain amount of time.

Finally, if $T^{(0)}_a=0$ and $T^{(0)}_s=0$, one has that
$$ F \begin{pmatrix} T_a =0 \\ T_s =0 \end{pmatrix}
= \left( \begin{array}{c} \frac{1}{\gamma _a} \mathcal R_a  \\ \frac{1}{\gamma _s}  \mathcal R_s  \end{array} \right) .$$
We recall that $\mathcal R_s (0) >0$. If also $\mathcal R_a (0) >0$, then once again we deduce that $F$ points inward the set $\mathcal Q$. If instead $\mathcal R_a (0)= 0$, then $T_a '(0)=0= T_a (0)$, that implies $T_a (t)= o(t)$. On the other hand, $T_s '(0) = \frac{1}{\gamma _s}  \mathcal R_s (0)$, hence close to $t=0$ we have that $T_s(t) = \frac{\mathcal R_s (0)}{\gamma _s} t + o(t)$. We note that
\begin{equation*}
\begin{split}
T_a ' (t)& = 
\frac{1}{\gamma _a} \Bigl[ -\lambda (T_a - T_s) + \varepsilon _a \sigma _B \vert T_s \vert ^3 T_s  - 2 \varepsilon _a \sigma _B \vert T_a \vert ^3 T_a \Bigr] + \frac{1}{\gamma _a} \mathcal R_a
\\ 
&\geq \frac{1}{\gamma _a} \Bigl[ -\lambda (T_a - T_s) + \varepsilon _a \sigma _B \vert T_s \vert ^3 T_s  - 2 \varepsilon _a \sigma _B \vert T_a \vert ^3 T_a \Bigr] ,
\end{split}
\end{equation*}
where we have used that $ \mathcal R_a$ is nonnegative. Le us now treat separately the cases $\lambda >0$ and $\lambda =0$.

If $\lambda >0$, we have
\begin{equation*}
\begin{split}
\frac{1}{\gamma _a} \Bigl[ -\lambda (T_a - T_s) &+ \varepsilon _a \sigma _B \vert T_s \vert ^3 T_s  - 2 \varepsilon _a \sigma _B \vert T_a \vert ^3 T_a \Bigr]
\\
&\,\,\,\,= \frac{1}{\gamma _a} \Bigl[ \lambda T_s \Bigr]
+ \frac{1}{\gamma _a}  \Bigl[ -\lambda T_a + \varepsilon _a \sigma _B \vert T_s \vert ^3 T_s  - 2 \varepsilon _a \sigma _B \vert T_a \vert ^3 T_a \Bigr]
\\
&\,\,\,\,= \frac{1}{\gamma _a} \Bigl[ \lambda \frac{\mathcal R_s (0)}{\gamma _s} t + o(t) \Bigr] + \frac{1}{\gamma _a} \Bigl[  o(t) \Bigr] = \frac{1}{\gamma _a} \Bigl[ \lambda \frac{\mathcal R_s (0)}{\gamma _s} t + o(t) \Bigr] ,
\end{split}
\end{equation*}
which implies that
$$  T_a '(t) \geq \frac{1}{\gamma _a} \Bigl[ \lambda  \frac{\mathcal R_s (0) }{\gamma _s }  t + o(t) \Bigr] \quad \text{ if } \lambda > 0 .$$
By integration, we get that $T_a (t) >0$ for $t>0$ small. Since also $T_s>0$ for $\lambda$ small enough, we conclude that the solution stays in $\mathcal Q$ for a certain amount of time.

The same property is true  if $\lambda =0$: indeed, in this case, we have that for $t>0$ small enough
\begin{equation*}
\begin{split}
\frac{1}{\gamma _a} \Bigl[ -\lambda (T_a - T_s) + \varepsilon _a \sigma _B \vert T_s \vert ^3 T_s & - 2 \varepsilon _a \sigma _B \vert T_a \vert ^3 T_a \Bigr]\\
&= \frac{1}{\gamma _a} \Bigl[ \varepsilon _a \sigma _B \vert T_s \vert ^3 T_s  - 2 \varepsilon _a \sigma _B \vert T_a \vert ^3 T_a \Bigr] 
\\
&= \frac{1}{\gamma _a} \Bigl[ \varepsilon _a \sigma _B \Bigl( \frac{\mathcal R_s (0)}{\gamma _s} t + o(t) \Bigr) ^4   - 2 \varepsilon _a \sigma _B o(t) ^4 \Bigr] \\
&= \frac{1}{\gamma _a} \Bigl[ \varepsilon _a \sigma _B  \frac{\mathcal R_s (0) ^4}{\gamma _s ^4} t^4 + o(t^4)  \Bigr]  .
\end{split}
\end{equation*}
Therefore 
$$  T_a ' (t) \geq \frac{1}{\gamma _a} \Bigl[ \varepsilon _a \sigma _B  \frac{\mathcal R_s (0) ^4}{\gamma _s ^4} t^4 + o(t^4)  \Bigr] ,$$
and by integrating in time, we obtain that $T_a (t) >0$ for $t>0$ small. Recalling that $T_s (t) >0$ for $t>0$ small, then the solution remains in $\mathcal Q$ for an interval of time small enough. 

Let us now consider $\left( \begin{array}{c} T_a ^{(0)} \\ T_s ^{(0)} \end{array} \right) \in \overline{\mathcal Q}$.
We claim that the associated solution stays in $\mathcal Q$ for all positive times (as long as it exists).
We have already proved that the solution belongs to $\mathcal Q$ for small positive times. 
Now we proceed by contradiction: assume that the solution leaves $\mathcal Q$ and let $t_0$ be the first exit time. Then, either $(T_a (t_0) >0, T_s(t_0)=0)$, or $(T_a (t_0) =0, T_s(t_0)>0)$, or $(T_a (t_0) =0, T_s(t_0)=0)$. However, in the first two cases, $F$ would point inward, a contradiction with the minimality of $t_0$. In the latter case we would have that $T_s '(t_0) = \frac{1}{\gamma _s} \mathcal R_s (0) >0$. However, since $T_s (t_0)=0$, we would get that
$$ T_s (t) = \frac{\mathcal R_s (0)}{\gamma _s} (t-t_0) + o(t-t_0) ,$$
which implies that $T_s$ would be negative before $t_0$. This fact contradicts the minimality of $t_0$. Therefore, as long as the solution exists, it remains inside $\mathcal Q$. \qed


\subsection{Bounds on the solution for $\varepsilon _a \in (0,2)$}\hfill
\label{sec-boundedness}

We will prove the following

\begin{Lemma}
\label{lem-bounds}
Assume that $\varepsilon _a \in (0,2)$. Let $\mu \in (\varepsilon _a ^{1/4}, 2^{1/4})$ and $M_a$ be large enough such that
$$ (T_a ^{(0)} , T_s ^{(0)}) \in [0,M_a] \times [0, \mu M_a] .$$
Then, the solution $(T_a,T_s)$ of \eqref{2layer-pbm-EDO} does not leave the rectangle $[0,M_a] \times [0, \mu M_a]$ for positive times.
\end{Lemma}

From the above Lemma it is easy to deduce the following

\begin{Corollary}\label{coroll-bound-max-sol}
Assume that $\varepsilon _a \in (0,2)$.
Then, if $T_a ^{(0)}$ and $T_s ^{(0)}$ are nonnegative, 
the maximal solution of \eqref{2layer-pbm-EDO} exists on $[0,+\infty)$.
\end{Corollary}

\noindent {\it Proof of Lemma \ref{lem-bounds}.}  
Let the initial condition be on the rectangle $[0,M_a]\times [0,M_s]$. On the right hand side of the rectangle, we have
\begin{equation*}
F \begin{pmatrix} M_a \\ T_s \end{pmatrix}
=\begin{pmatrix}\frac{1}{\gamma _a} \Bigl[ -\lambda (M_a - T_s) + \varepsilon _a \sigma _B \vert T_s \vert ^3 T_s  - 2 \varepsilon _a \sigma _B M_a ^4 + \mathcal R_a (M_a) \Bigr] \\ \cdots \end{pmatrix},
\end{equation*}
while on the top side, we have
\begin{equation*}
F \begin{pmatrix} T_a \\ M_s \end{pmatrix} 
= \begin{pmatrix} \cdots \\ \frac{1}{\gamma _s} \Bigl[ -\lambda (M_s - T_a) - \sigma _B M_s ^4  + \varepsilon _a \sigma _B \vert T_a \vert ^3 T_a + \mathcal R_s (M_s) \Bigr] \end{pmatrix} .
\end{equation*}
We are interested in the sign of first component of $F \left( \begin{array}{c} M_a \\ T_s \end{array} \right)$ and 
in the sign of the second component of $ F \left( \begin{array}{c} T_a \\ M_s \end{array} \right)$ to check if $F$ points inward the rectangle on these two sides.
We note that
\begin{multline*}
 -\lambda (M_a - T_s) + \varepsilon _a \sigma _B \vert T_s \vert ^3 T_s  - 2 \varepsilon _a \sigma _B M_a ^4 + \mathcal R_a (M_a)
\\
\leq -\lambda (M_a - M_s) + \varepsilon _a \sigma _B M_s ^4  - 2 \varepsilon _a \sigma _B M_a ^4 + \mathcal R_a (M_a)
\\
= \varepsilon _a \sigma _B \Bigl[  M_s ^4 -2 M_a ^4 \Bigr]  -\lambda (M_a - M_s) + \mathcal R_a (M_a),
\end{multline*}
and
\begin{multline*}
 -\lambda (M_s - T_a) - \sigma _B M_s ^4  + \varepsilon _a \sigma _B \vert T_a \vert ^3 T_a + \mathcal R_s (M_s)
\\
\leq -\lambda (M_s - M_a) - \sigma _B M_s ^4  + \varepsilon _a \sigma _B M_a ^4 + \mathcal R_s (M_s) 
\\
= \sigma _B \Bigl[ \varepsilon _a  M_a ^4 - M_s ^4 \Bigr] -\lambda (M_s - M_a) + \mathcal R_s (M_s).
\end{multline*}
To simplify the analysis, we set
$$ M_s = \mu M_a ,$$
with $\mu$ to be chosen later. Then 
$$ \Bigl[  M_s ^4 -2 M_a ^4 \Bigr] = (\mu ^4 - 2 ) M_a ^4 
\text{ and }  \Bigl[ \varepsilon _a  M_a ^4 - M_s ^4 \Bigr] = (  \varepsilon _a - \mu ^4) M_a ^4 .$$
Therefore, if $\varepsilon _a \in (0,2)$, by choosing $\mu \in (\varepsilon _a ^{1/4}, 2^{1/4})$, we have
$$ \mu ^4 - 2 < 0 \quad \text{ and } \quad  \varepsilon _a - \mu ^4 < 0 .$$
Then, since $\mathcal R_a$ and $\mathcal R_s$ are globally Lipschitz, $F$ points inward the rectangle $[0,M_a]\times [0,\mu M_a]$ on all the sides if $M_a$ is chosen sufficiently large. We 
finally obtain that this rectangle is invariant: if $\mu \in (\varepsilon _a ^{1/4}, 2^{1/4})$ and $M_a$ sufficiently large, then
\begin{equation*}
    \begin{pmatrix} T_a ^{(0)} \\ T_s ^{(0)} \end{pmatrix} \in [0,M_a]\times [0,\mu M_a]\quad \implies \quad \begin{pmatrix} T_a (t)  \\ T_s (t)  \end{pmatrix} \in (0,M_a)\times (0,\mu M_a)  
\end{equation*}
as soon as $t>0$ and as long as the solution exists. Thus, we further deduce that
\begin{equation*}
    \varepsilon _a \in (0,2) , T_a ^{(0)} \geq 0, T_s ^{(0)} \geq 0 \quad
    \implies \quad \text{ the solution $\begin{pmatrix} T_a \\ T_s \end{pmatrix}$
is bounded} ,
\end{equation*}
which yields global existence (for positive times). This concludes the proof of Proposition \ref{prop-wellposed}.
\qed


\subsection{Proof of Proposition \ref{prop-comp-lambda-eq-cv}: general asymptotic behaviour} \hfill

We first recall some definitions from \cite{Smith}.
A set $D \subset \RR^n$ is \emph{$p$-convex} (see \cite[page 33]{Smith}) if for all  $x, y \in D$ such that $x_i \leq y_i$ for all  $i=1,\dots,n$
$$[x,y] \subset D .$$
Given a $p$-convex set $D$, a $C^1$-system of differential equations on $D$
\begin{equation*}
x' _i (t) =F_i(x_1 (t),\dots,x_n (t))=F_i(x (t))\qquad i=1,\dots,n
\end{equation*}
is called \emph{cooperative} if
\begin{equation*}
\frac{\partial F_i}{\partial x_j}(x) \geq 0\qquad i\neq j , \quad x \in D ,
\end{equation*}
and \emph{competitive} if the reverse inequalities hold (see \cite[page 34]{Smith}).

Furthermore, we recall the following result (see \cite[Theorem 2.2, page 3]{Smith}).
\begin{Theorem}\label{Teo-Smith}
Consider a cooperative or competitive system on a $p$-convex set $D \subset \RR^2$.

If $t\mapsto x(t)$ is a solution defined on $[0,+\infty)$, then there exists $T\geq0$ such that $x(t)$ is monotone for $t\geq T$.

Moreover, if the solution $x(t)$ is bounded, then it converges to some equilibrium point.
\end{Theorem}
Let us consider the $2$-convex set $\overline{\mathcal{Q}}=[0,+\infty) \times [0,+\infty)$. Let $\varepsilon _a \in (0,2)$ and $\lambda\geq 0$ and observe that
 $$ T_s \mapsto \frac{1}{\gamma_a}\left[\lambda T_s+\varepsilon _a \sigma _B \vert T_s \vert ^3 T_s -\lambda T_a - 2 \varepsilon _a \sigma _B \vert T_a \vert ^3 T_a + \mathcal R_a (T_a)\right]$$
and
$$ T_a \mapsto \frac{1}{\gamma_s}\left[\lambda T_a- \sigma _B \vert T_s \vert ^3 T_s  -\lambda T_s+ \varepsilon _a \sigma _B \vert T_a \vert ^3 T_a + \mathcal R_s (T_s)\right]$$
are nondecreasing on $[0,+\infty)$ (indeed it is easy to check that the derivatives of these maps are nonnegative). Thus, system \eqref{2layer-pbm-EDO} is cooperative if the vector field $F$ is $C^1$. However, $F=(f_1,f_2)$ is globally Lipschitz continuous because of functions $\beta _a$ and $\beta _s$. Thus, since we cannot apply Theorem \ref{Teo-Smith} to our system, we need to verify if the monotonicity property of the solution still holds under our weaker assumptions.

Let us consider the following subset of $\overline{\mathcal Q}$
$$
P_+ := \{ (T_a, T_s) \in \overline{\mathcal Q}\,:\, f_1 (T_a,T_s) \geq 0 \text{ and } f_2 (T_a,T_s) \geq 0 \} ,
$$
$$
P_- := \{ (T_a, T_s)\in \overline{\mathcal Q}\,:\, f_1 (T_a,T_s) \leq 0 \text{ and } f_2 (T_a,T_s) \leq 0 \} . $$
We are going to prove the following

\begin{Lemma}
\label{lem-mono-lip}
If $(T_a ^{(0)}, T_s ^{(0)}) \in P_+$ then each component $t\mapsto T_a(t)$, $t\mapsto F_s(t)$ of the solution of \eqref{2layer-pbm-EDO} is nondecreasing, that is, the solution of \eqref{2layer-pbm-EDO} nondecreasing.

Similarly, if $(T_a ^{(0)}, T_s ^{(0)}) \in P_-$ then the solution of \eqref{2layer-pbm-EDO} is nonincreasing.
\end{Lemma}
Thanks to the above result we can prove the following Theorem.
\begin{Theorem}
\label{lem-cv-lip}
Given $(T_a ^{(0)}, T_s ^{(0)}) \in \overline{\mathcal Q}$, the solution $(T_a(t),T_s(t))$ of \eqref{2layer-pbm-EDO} converges to some equilibrium point. Moreover, there exists $T\geq0$ such that $t\mapsto T_a(t)$ and $t\mapsto T_s(t)$ are monotone for $t\geq T$.
\end{Theorem}

\noindent {\it Proof of Lemma \ref{lem-mono-lip}.} 
We prove the result using a regularization argument. Given $n\geq 1$, there exist $\beta _{a,n}, \beta _{s,n}: \RR \to \RR$ of class $C^1$, globally Lipschitz continuous, such that $\beta _{a,n} \geq 0$, $\beta _{s,n} >0$ and
$$ \sup _{\RR} \vert \beta _{a,n} - \beta _a \vert \leq \frac{1}{n}, \quad 
\text{ and } \quad \sup _{\RR} \vert \beta _{s,n} - \beta _s \vert \leq \frac{1}{n} .$$
Given an initial condition $(T_a ^{(0)}, T_s ^{(0)}) \in \overline{\mathcal Q}$, consider the associated solution
$(T_{a,n}, T_{s,n})$ of the regularized problem
$$ \begin{cases}
\left( \begin{array}{c} T_a \\ T_s \end{array} \right) ' = F_n \left( \begin{array}{c} T_a \\ T_s \end{array} \right)  , \\\\
\left( \begin{array}{c} T_a \\ T_s \end{array} \right) (0) = \left( \begin{array}{c} T_a ^{(0)} \\ T_s ^{(0)} \end{array} \right) ,
\end{cases} $$
where $F_n$ is defined by replacing $\beta _a$ and $\beta _s$ by $\beta _{a,n}$ and $\beta _{s,n}$, respectively, in \eqref{def-Fnonlin}. The regularized vector field $F_n$ is $C^1$ and cooperative (indeed,
its first component, denoted by $f_{n,1}$, is nondecreasing with respect to $T_s$ and the second component, denoted by $f_{n, 2}$, is nondecreasing with respect to $T_a$.)

We now consider $(T_a ^{(0)}, T_s ^{(0)}) \in P_+$ and we assume that
$$f_1 (T_a ^{(0)},T_s ^{(0)}) > 0  ,$$
$$f_2 (T_a ^{(0)},T_s ^{(0)}) > 0 . $$
Since it holds that
$$f_{n,1} (T_a ^{(0)},T_s ^{(0)}) \to f_1 (T_a ^{(0)},T_s ^{(0)}) ,$$
$$f_{n,2} (T_a ^{(0)},T_s ^{(0)}) \to f_2 (T_a ^{(0)},T_s ^{(0)}) , $$
as $n\to \infty$, then for $n$ large enough
$$f_{n,1} (T_a ^{(0)},T_s ^{(0)}) > 0 ,$$
$$f_{n,2} (T_a ^{(0)},T_s ^{(0)}) > 0. $$
Recalling that the regularized system is cooperative, we get that
\begin{equation*}
\begin{split}
&\frac{d}{dt} \Bigl( f_{n,1} (T_{a,n}(t), T_{s,n}(t)) \  f_{n,2} (T_{a,n}(t), T_{s,n}(t)) \Bigr) =
\\
&\,= \Bigl( \frac{\partial f_{n,1}}{\partial T_a} (T_{a,n}(t), T_{s,n}(t)) + \frac{\partial f_{n,2}}{\partial T_s} (T_{a,n}(t), T_{s,n}(t)) \Bigr)\\
&\qquad\qquad\qquad\qquad \cdot\Bigl( f_{n,1} (T_{a,n}(t), T_{s,n}(t)) \  f_{n,2} (T_{a,n}(t), T_{s,n}(t)) \Bigr) 
\\
&\quad+ \frac{\partial f_{n,1}}{\partial T_s} (T_{a,n}(t), T_{s,n}(t)) \, f_{n,2}^2 (T_{a,n}(t), T_{s,n}(t))\\
&\quad
+ \frac{\partial f_{n,2}}{\partial T_a} (T_{a,n}(t), T_{s,n}(t)) \, f_{n,1}^2 (T_{a,n}(t), T_{s,n}(t))
\\
&\,\geq ( \nabla\cdot F_n \begin{pmatrix} T_{a,n} (t) \\  T_{s,n} (t) \end{pmatrix})  \Bigl( f_{n,1} (T_{a,n}(t), T_{s,n}(t))\, f_{n,2} (T_{a,n}(t), T_{s,n}(t)) \Bigr) .
\end{split}
\end{equation*}
By integrating the above first order differential inequality, we deduce that the quantity 
$$ f_{n,1} (T_{a,n}(t), T_{s,n}(t)) \  f_{n,2} (T_{a,n}(t), T_{s,n}(t)) $$
remains positive. Therefore, for every $t\geq0$
$$ f_{n,1} (T_{a,n}(t), T_{s,n}(t)) > 0 ,$$
$$f_{n,2} (T_{a,n}(t), T_{s,n}(t)) > 0 . $$
This implies that the solution $(T_{a,n}, T_{s,n})$ of the regularized problem is increasing in time and
since for every $t\geq0$ it holds that
$$ T_{a,n} (t) \to _{n \to \infty} T_a (t),$$
$$T_{s,n} (t) \to _{n \to \infty} T_s (t), $$
we deduce that the solution $(T_a, T_s)$ of the original problem is nondecreasing in time.

We now prove that the same property is true if 
$$ f_1 (T_a ^{(0)},T_s ^{(0)}) = 0  ,$$
$$f_2 (T_a ^{(0)},T_s ^{(0)}) > 0. $$
Consider an initial condition $(T_a ^{(0)},T_s ^{(0)})$ such that the above relations hold.
Then, since $f_1$ is strictly increasing with respect to the second variable, we have that for all $\eta >0$
$$  f_1 (T_a ^{(0)},T_s ^{(0)} + \eta) >0 .$$
Moreover,
$$ f_2 (T_a ^{(0)},T_s ^{(0)} +\eta ) \to f_2 (T_a ^{(0)},T_s ^{(0)}) \quad \text{ as } \eta \to 0 ^+ ,$$
and therefore we have that for $\eta >0$ small enough
$$ f_1 (T_a ^{(0)},T_s ^{(0)}+\eta) > 0  ,$$
$$f_2 (T_a ^{(0)},T_s ^{(0)}+\eta) > 0. $$
Then, the solution $(T_a ^{(\eta)}, T_s ^{(\eta)})$ corresponding to the initial condition $(T^{(0)}_a,T^{(0)}_s+\eta)$ is nondecreasing in time thanks to the previous step. Furthermore, since for every $t\geq0$ it holds that
$$T_a ^{(\eta)} (t) \to _{\eta \to 0^+} T_a (t),$$
$$T_s ^{(\eta)} (t) \to _{\eta \to 0^+} T_s (t),$$
we conclude that the solution $(T_a, T_s)$ of the original problem is nondecreasing in time, that is, what we wanted to prove.
 
If we now consider the case
$$ f_1 (T_a ^{(0)},T_s ^{(0)}) > 0  ,$$
$$ f_2 (T_a ^{(0)},T_s ^{(0)}) = 0 ,$$
we can proceed similarly to the previous case by introducing the approximated initial condition $(T_a ^{(0)} +\eta,T_s ^{(0)})$ and deduce that the solution of the original problem is again nondecreasing. This concludes the case $(T_a ^{(0)} ,T_s ^{(0)}) \in P_+$.

One can adapt the strategy proposed for $(T^{(0)}_a,T^{(0)}_s)\in P_+$ to the case of an initial condition in $P_-$ and conclude the proof of Lemma \ref{lem-mono-lip}.\qed
 
\medskip

\noindent {\it Proof of Theorem \ref{lem-cv-lip}.} If $(T_a ^{(0)},T_s ^{(0)}) \in P_+$, then the solution is nondecreasing, and bounded thanks to Lemma \ref{lem-bounds}. Therefore, it converges to some limit, that has to be an equilibrium point. Similarly, if $(T_a ^{(0)},T_s ^{(0)}) \in P_-$ then the solution converges to some limit, that again has to be an equilibrium point.
This solves the case of an initial condition that verifies the property
$$ f_1 (T_a ^{(0)},T_s ^{(0)}) \, 
 f_2 (T_a ^{(0)},T_s ^{(0)}) \geq 0. $$
 We now assume that 
$$ f_1 (T_a ^{(0)},T_s ^{(0)}) \, 
 f_2 (T_a ^{(0)},T_s ^{(0)}) < 0. $$
Then, either it holds that
 $$ \forall t \geq 0, \quad f_1 (T_a (t),T_s (t)) \,  f_2 (T_a (t),T_s (t)) < 0 ,$$
or there exists $\tau\geq0$ such that
 $$ f_1 (T_a (\tau),T_s (\tau)) \,  f_2 (T_a (\tau),T_s (\tau)) = 0 .$$
In the second case, the solution has entered the set $P_+ \cup P_-$ and from the previous analysis we deduce that it converges monotonically to some equilibrium point. In the first case, each function $t \mapsto f_1 (T_a (t),T_s (t))$ and $t \mapsto f_2 (T_a (t),T_s (t))$
has a precise sign for every $t\geq0$. Therefore $T_a$ and $T_s$ are monotone. Because of the boundedness the solution $(T_a,T_s)$ converges with some monotonicity to a limit which is an equilibrium point. This concludes the proof of Theorem \ref{lem-cv-lip}. \qed


\section{Analysis of the asymptotic behaviour for $\lambda =0$ and $\mathcal R_a =0$}
\label{sec-comp-asympt-lambda=0}

We assume here that 
$\lambda = 0$, $q>0$, $\sigma _B >0$, $\beta _a =0$ (considering that $\mathcal R_a$ is negligible with respect to $\mathcal R_s$), $\beta$ is given by \eqref{beta},
and $\varepsilon _a \in (0,2)$ (since the solutions are unbounded when $\varepsilon _a \geq 2$).

We already know that any solution converges to some equilibrium point. The goal here is to be more precise about the monotonicity of the solution.

\subsection{There are one, two or three equilibrium points} \hfill
\label{sub-sec-eq123}


Equilibrium points of \eqref{2layer-pbm-EDO} when $\lambda = 0$ and $\mathcal R_a =0$ are $(T_a,T_s)$ solutions of 
\begin{equation}
\label{equil} \begin{cases}
\varepsilon _a \sigma _B \vert T_s \vert ^3 T_s  - 2 \varepsilon _a \sigma _B \vert T_a \vert ^3 T_a = 0 , \\
- \sigma _B \vert T_s \vert ^3 T_s  + \varepsilon _a \sigma _B \vert T_a \vert ^3 T_a + \mathcal R_s (T_s) = 0 .
\end{cases}
\end{equation}
Since we are interested in the behaviour of positive initial conditions, we look for nonnegative equilibrium points. From the first equation of \eqref{equil}, we have
$$ T_s ^4 = 2 T_a ^4,$$
which gives
\begin{equation}
\label{equil-relTaT0}
T_s = 2^{1/4} T_a .
\end{equation}
Substituting this last identity in the second equation of \eqref{equil}, we obtain
$$\sigma _B (-1  + \frac{\varepsilon _a}{2}) T_s ^4  + \mathcal R_s (T_s)= 0,$$
that is equivalent to
\begin{equation}
\label{equil-eqT0}
\sigma _B (1  - \frac{\varepsilon _a}{2}) T_s ^4 = q \beta _s (T_s) .
\end{equation}
Observe that this last equation is of the same type than the one satisfied by equilibrium points of a one-layer model. By considering $\beta _s: \mathbb R \to \mathbb R$ as in \eqref{beta}, that is, continuous, positive, nondecreasing, constant on the intervals $[0,T_{s,-}]$ and $[T_{s,+},+\infty)$ and linear on $[T_{s,+},T_{s,+}]$, we deduce that there can be exactly one, two, or three equilibrium points, depending on the values of the parameters $\sigma _B$, $\varepsilon _a$, $q$, and those that appear in $\beta _s$. No other situations are possible. Indeed it is easy to see that function
\begin{equation}
\label{def-g}
g: [0, +\infty) \to \mathbb R , \quad g(T_s )= \sigma _B (1  - \frac{\varepsilon _a}{2}) T_s ^4 - q \beta _s (T_s) 
\end{equation} 
is continuous, negative at $T_s =0$ and goes to $+\infty$ as $T_s\to +\infty$. Thus, the mean value theorem ensures that there exists at least a point $T_s\in(0,+\infty)$ such that \eqref{equil-eqT0} is satisfied. Moreover, $g$ is strictly increasing on $[0,T_{s,-}]$, hence $g$ can be equal to $0$ on $[0,T_{s,-}]$ at most once. One can use the same argument to prove that also in $[T_{s,+}, +\infty)$ $g$ can have at most one zero. Finally, $g$ is strictly convex on $[T_{s,-},T_{s,+}]$, an so it can assume the value $0$ at most twice on this set. Such argument implies that $g$ can be equal to $0$ at most four times on $\mathbb R_+$. However, if it takes twice the value $0$ on $[T_{s,-},T_{s,+}]$, this would mean that the curve of $T_s \mapsto \sigma _B (1  - \frac{\varepsilon _a}{2}) T_s ^4$ intersects twice the segment $\{ (T_s,q\beta _s (T_s)), T_s \in [T_{s,-},T_{s,+}] \}$. In this case it cannot intersect anymore the half-line $\{(T_s, q\beta _{s,+}), T_s \in [T_{s,+}, +\infty)\}$ because, by convexity, it will remain above the half line on which the previous segment lies. Thus, this observation reduces the number of possible zeros of $g$ on $\mathbb R_+$ to three. We describe such values in the following pictures:


\begin{itemize}
\item the intersection between the graph of $T\mapsto \sigma _B (1  - \frac{\varepsilon _a}{2}) T ^4$
and the graph of $T\mapsto q \beta _s (T)$ is exactly one point (see Figure \ref{fig1:region})
\begin{figure}[h!]
\centering
\subfloat[]
{\begin{tikzpicture}[scale=.5]
\draw[-stealth] (-1,0) -- (7,0);
\draw[-stealth] (0,-1) -- (0,5);
\draw (0,1)node[left]{$q\beta_{s,-}$} -- (2,1);
\draw (2,1) -- (4,3);
\draw (4,3) -- (6,3);
\draw[dashed] (0,3)node[left]{$q\beta_{s,+}$} -- (4,3);
\draw[domain=0:5.25, smooth, variable=\x,green!50!black] plot ({\x}, {.15*\x*\x});
\node[green!50!black] at (4.7,4.5) {$\scriptstyle{T\mapsto \sigma_B\left(1-\frac{\eps_a}{2}\right)T^4}$};
\draw[dashed] (4.49,3) -- (4.49,0)node[below]{$\quad\scriptstyle{T_{s,1}^*}$};
\node at (6.8,3) {$q\beta _s$};
\end{tikzpicture}
}\,\,
\subfloat[]{
\begin{tikzpicture}[scale=.5]
\draw[-stealth] (-1,0) -- (7,0);
\draw[-stealth] (0,-1) -- (0,5);
\draw (0,1)node[left]{$q\beta_{s,-}$} -- (2,1);
\draw (2,1) -- (4,3);
\draw (4,3) -- (6,3);
\draw[dashed] (0,3)node[left]{$q\beta_{s,+}$} -- (4,3);
\draw[domain=0:4.5, smooth, variable=\x,green!50!black] plot ({\x}, {.21*\x*\x});
\node[green!50!black] at (4.3,4.5) {$\scriptstyle{T\mapsto \sigma_B\left(1-\frac{\eps_a}{2}\right)T^4}$};
\draw[dashed] (3.34,2.3) -- (3.34,0)node[below]{$\scriptstyle{T_{s,1}^*}$};
\node at (6.8,3) {$q\beta _s$};
\end{tikzpicture}
}\,\,
\subfloat[]{
\begin{tikzpicture}[scale=.5]
\draw[-stealth] (-1,0) -- (7,0);
\draw[-stealth] (0,-1) -- (0,5);
\draw (0,1)node[left]{$q\beta_{s,-}$} -- (2,1);
\draw (2,1) -- (4,3);
\draw (4,3) -- (6,3);
\draw[dashed] (0,3)node[left]{$q\beta_{s,+}$} -- (4,3);
\draw[domain=0:1.9, smooth, variable=\x,green!50!black] plot ({\x}, {1.1*\x*\x});
\node[green!50!black] at (3.2,4.5) {$\scriptstyle{T\mapsto \sigma_B\left(1-\frac{\eps_a}{2}\right)T^4}$};
\draw[dashed] (.95,1) -- (.95,0)node[below]{$\scriptstyle{T_{s,1}^*}$};
\node at (6.8,3) {$q\beta _s$};
\end{tikzpicture}
\,\,
}
\caption{In this figure we show the possible cases of one intersection between the curves $T\mapsto \sigma_B(1-\frac{\eps_a}{2})T^4$ and $T\mapsto q\beta_s(T)$. In (a) there is a unique ``warm" equilibrium, in (c) a unique ``cold" equilibrium and in (b) an equilibrium at an intermediate surface temperature between $T_{s,-}$ and $T_{s,+}$.}\label{fig1:region}
\end{figure}
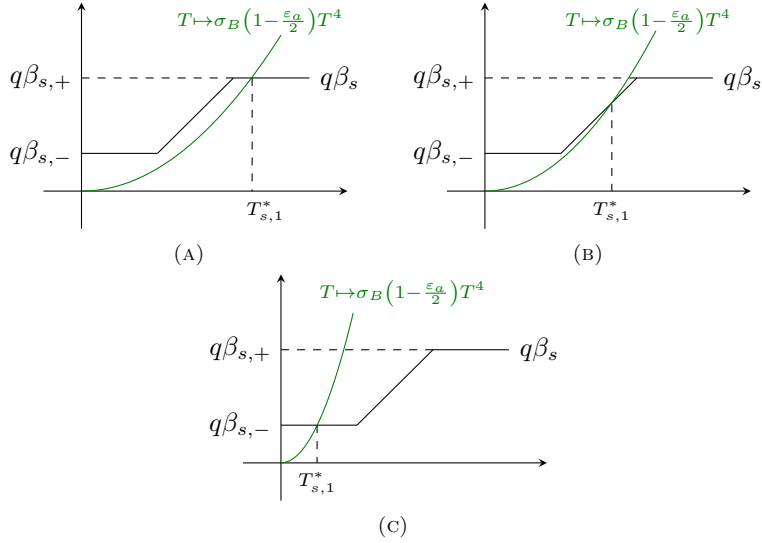

\noindent or the \lq\lq degenerate" cases where the intersection is either $T_{s,-}$ or $T_{s,+}$,
		
\item the graph of $T\mapsto \sigma _B (1  - \frac{\varepsilon _a}{2}) T_s ^4$ intersects the graph of $T\mapsto q \beta _s (T)$ three times (see Figure \ref{fig2:region})
\begin{figure}[h!]
\centering
\subfloat[]{
\begin{tikzpicture}[scale=.5]
\draw[-stealth] (-1,0) -- (6,0);
\draw[-stealth] (0,-1) -- (0,5);
\draw (0,1)node[left]{$q\beta_{s,-}$} -- (2.5,1);
\draw (2.5,1) -- (3,3);
\draw (3,3) -- (5,3);
\draw[dashed] (0,3)node[left]{$q\beta_{s,+}$} -- (3,3);
\draw[domain=0:4.5, smooth, variable=\x,green!50!black] plot ({\x}, {.2*\x*\x});
\node[green!50!black] at (4.3,4.5) {$\scriptstyle{T\mapsto \sigma_B\left(1-\frac{\eps_a}{2}\right)T^4}$};
\draw[dashed] (2.25,1) -- (2.25,0)node[below]{$\scriptstyle{T_{s,1}^*}\quad$};
\draw[dashed] (2.58,1.3) -- (2.58,0)node[below]{$\quad\scriptstyle{T_{s,2}^*}$};
\draw[dashed] (3.86,3) -- (3.86,0)node[below]{$\scriptstyle{T_{s,3}^*}$};
\node at (5.8,3) {$q\beta _s$};
\end{tikzpicture}
}\quad
\subfloat[]{
\begin{tikzpicture}[scale=.5]
\draw[-stealth] (-1,0) -- (6,0);
\draw[-stealth] (0,-1) -- (0,5);
\draw (4,3) -- (5,3);
\draw (1.25, .25) -- (4,3);
\draw (0,.25)node[left]{$q\beta_{s,-}$} -- (1.25,.25);
\draw[dashed] (0,3)node[left]{$q\beta_{s,+}$} -- (4,3);
\draw[domain=0:4.5, smooth, variable=\x,green!50!black] plot ({\x}, {.2*\x*\x});
\node[green!50!black] at (4,4.5) {$\scriptstyle{T\mapsto \sigma_B\left(1-\frac{\eps_a}{2}\right)T^4}$};
\draw[dashed] (1.1,.25) -- (1.1,0)node[below]{$\scriptstyle{T_{s,1}^*}\quad$};
\draw[dashed] (1.4,.4) -- (1.4,0)node[below]{$\quad \scriptstyle{T_{s,2}^*}$};
\draw[dashed] (3.6,2.6) -- (3.6,0)node[below]{$\quad \scriptstyle{T_{s,3}^*}$};
\node at (5.8,3) {$q\beta _s$};
\end{tikzpicture}
}	
\caption{In (a) and (b) we show possible cases of three intersection between the curves $T\mapsto \sigma_B(1-\frac{\eps_a}{2})T^4$ and $T\mapsto q\beta_s(T)$. In particular, in (a) there are one ``cold", one intermediate and one ``warm" equilibrium, while in (b) there are one ``cold" and two intermediate equilibrium.}
\label{fig2:region}
\end{figure}
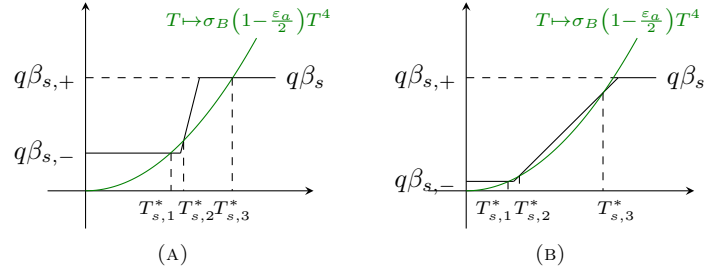

\item the \lq\lq degenerate" cases where the intersections between the graph of $T\mapsto \sigma _B (1  - \frac{\varepsilon _a}{2}) T_s ^4$
and the graph of $T\mapsto q \beta _s (T)$ are exactly two (see Figure \ref{fig3:region})

\begin{figure}[h!]
\centering
\subfloat[]{
\begin{tikzpicture}[scale=.4]
\draw[-stealth] (-1,0) -- (6,0);
\draw[-stealth] (0,-1) -- (0,5);
\draw (0,.8)node[left]{$q\beta_{s,-}$} -- (1.95,.8);
\draw (1.95,0.8) -- (4,3);
\draw (4,3) -- (5.5,3);
\draw[dashed] (0,3)node[left]{$q\beta_{s,+}$} -- (4,3);
\draw[domain=0:4.5, smooth, variable=\x,green!50!black] plot ({\x}, {.21*\x*\x});
\node[green!50!black] at (3.9,4.8) {$\scriptstyle{T\mapsto \sigma_B\left(1-\frac{\eps_a}{2}\right)T^4}$};
\draw[dashed] (3.18,2.1) -- (3.2,0)node[below]{$\scriptstyle{T_{s,2}^*}$};
\draw[dashed] (1.95,.8) -- (1.95,0)node[below]{$\scriptstyle{T_{s,1}^*}$};
\node at (6.8,3) {$q\beta _s$};
\end{tikzpicture}
}\quad
\subfloat[]{
\begin{tikzpicture}[scale=.4]
\draw[-stealth] (-1,0) -- (6,0);
\draw[-stealth] (0,-1) -- (0,5);
\draw (3.6,2.6) -- (5,2.6);
\draw (1.4, .4) -- (3.6,2.6);
\draw (0,.4) node[left]{$q\beta_{s,-}$} -- (1.4,.4) ;
\draw[dashed] (0,2.6)node[left]{$q\beta_{s,+}$} -- (3.6,2.6);
\draw[domain=0:4.5, smooth, variable=\x,green!50!black] plot ({\x}, {.2*\x*\x});
\node[green!50!black] at (3.9,4.8) {$\scriptstyle{T\mapsto \sigma_B\left(1-\frac{\eps_a}{2}\right)T^4}$};
\draw[dashed] (1.4,.4) -- (1.4,0)node[below]{$\scriptstyle{T_{s,1}^*}$};
\draw[dashed] (3.6,2.6) -- (3.6,0)node[below]{$\scriptstyle{T_{s,2}^*}$};
\node at (6.2,2.6) {$q\beta _s$};
\end{tikzpicture}
}
\quad
\subfloat[]{
\begin{tikzpicture}[scale=.4]
\draw[-stealth] (-1,0) -- (6,0);
\draw[-stealth] (0,-1) -- (0,5);
\draw (4,2.6) -- (5,2.6);
\draw (1.6, .4) -- (4,2.6);
\draw (0,.4) node[left]{$q\beta_{s,-}$} -- (1.6,.4) ;
\draw[dashed] (0,2.6)node[left]{$q\beta_{s,+}$} -- (3.6,2.6);
\draw[domain=0:4.5, smooth, variable=\x,green!50!black] plot ({\x}, {.2*\x*\x});
\node[green!50!black] at (3.9,4.8) {$\scriptstyle{T\mapsto \sigma_B\left(1-\frac{\eps_a}{2}\right)T^4}$};
\draw[dashed] (1.4,.4) -- (1.4,0)node[below]{$\scriptstyle{T_{s,1}^*}$};
\draw[dashed] (2.4,1.1) -- (2.4,0)node[below]{\,\,$\scriptstyle{T_{s,2}^*}$};
\node at (6.2,2.6) {$q\beta _s$};
\end{tikzpicture}
}
\caption{In these pictures we represents three possible cases of two intersections between the curves $T\mapsto \sigma_B(1-\frac{\eps_a}{2})T^4$ and $T\mapsto q\beta_s(T)$. In (a) and in (c) there are one ``cold" and one intermediate equilibrium and in (b) one ``cold" and one ``warm" equilibrium.}\label{fig3:region}
\end{figure}
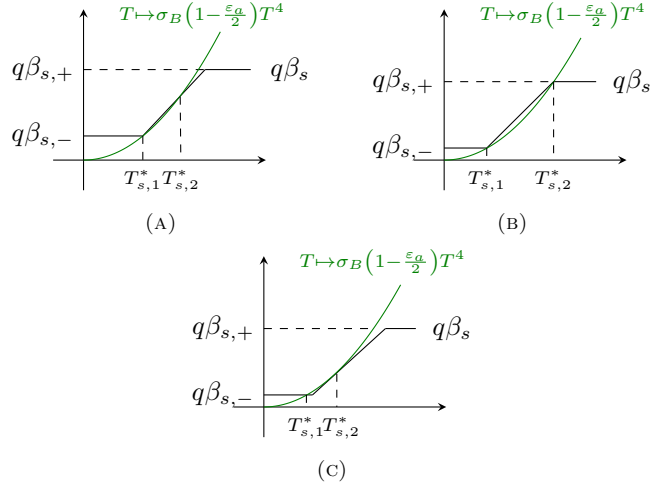

\end{itemize}

In the case where we have only one equilibrium point, all solutions converge to such equilibrium. When there are more than one equilibrium points, every solution converges to one of them. In the following we precise the nature of such equilibria.

We will also show that the equilibrium points have some monotonicity property with respect to the parameter $\varepsilon _a$. More precisely, we will prove that
\begin{itemize}
\item in situations as those described in Figure \ref{fig1:region}, (A) the equilibrium point has the same monotony of $\varepsilon _a$, that is, if $\varepsilon _a$ increases the equilibrium point increases;

\item in situations as those of Figure \ref{fig2:region}, (A): the equilibrium points enjoy the following property: if $\varepsilon _a$ increases, the \lq\lq cold" and the \lq\lq warm" equilibrium points increase, while the intermediate decrease.
\end{itemize}
We will extend later on these properties to the case $\lambda >0$.

For the sake of shortness, in what follows we analyse the most common cases of one a three equilibria. With similar techniques one can study the case of two intersection between the curves $q\beta_s(T)$ and $\sigma_B(1-\frac{\varepsilon_a}{2})T^4$.


\subsection{Monotonicity and convergence when there is only one equilibrium}\hfill
\label{sec-mono-1eq}

Here we study the case where 
\eqref{equil-eqT0} has one and only one solution that we denote by $T_{s,1} ^*$. This case corresponds to a unique equilibrium point of \eqref{2layer-pbm-EDO}, which is $(2^{-1/4} T_{s,1} ^*, T_{s,1} ^*)$. Note that we have
$$ \begin{cases}
\sigma _B (1  - \frac{\varepsilon _a}{2}) T_s ^4 < q\beta _s (T_s) & \text{ if } T_s < T_{s,1} ^* , \\
\sigma _B (1  - \frac{\varepsilon _a}{2}) T_s ^4 = q\beta _s (T_s) & \text{ if } T_s = T_{s,1} ^* , \\
\sigma _B (1  - \frac{\varepsilon _a}{2}) T_s ^4 > q\beta _s (T_s) & \text{ if } T_s > T_{s,1} ^* .
\end{cases}
$$
In the phase plane, let us consider once again the line 
$$\mathcal C_1: T_s = 2^{1/4} T_a$$ 
and the set
$$ \mathcal C_2:= \{(T_a, T_s) \in \mathbb R^2, - \sigma _B \vert T_s \vert ^3 T_s  + \varepsilon _a \sigma _B \vert T_a \vert ^3 T_a + q \beta _s (T_s) =0 \}.$$
Let us analyse the set $\mathcal C_2$. We claim that, given $T_s$, there exists one and only one value $T_a$, denoted $T_a ^{(2)} (T_s)$, such that
$$ \varepsilon _a \sigma _B \vert T_a \vert ^3 T_a  = \sigma _B \vert T_s \vert ^3 T_s - q\beta _s (T_s) .$$
Indeed, the function 
\begin{equation}
\label{def-psi}
\psi: \mathbb R \to \mathbb R, \quad \psi (T_a)= - \sigma _B \vert T_s \vert ^3 T_s  + \varepsilon _a \sigma _B \vert T_a \vert ^3 T_a + q \beta _s (T_s) 
\end{equation}
is strictly increasing and has infinite limits as $T_a \to \pm \infty$. Moreover, if $T_s \geq 0$, we have
\begin{equation*}
    \psi (2^{-1/4} T_s) 
= - \sigma _B T_s ^4  + \frac{\varepsilon _a}{2} \sigma _B T_s ^4 + q\beta _s (T_s)= q\beta _s(T_s) - \sigma _B (1 - \frac{\varepsilon _a}{2}) T_s ^4 .
\end{equation*}
Therefore 
$$ \begin{cases}
\psi (2^{-1/4} T_s) >0 &\text{ if } T_s < T_{s,1} ^* , \\
\psi (2^{-1/4} T_s) =0 &\text{ if } T_s = T_{s,1} ^*, \\
\psi (2^{-1/4} T_s) < 0 &\text{ if } T_s > T_{s,1} ^* .
\end{cases} $$
Since $0= \psi (T_a ^{(2)} (T_s))$, the monotonicity of $\psi$
gives
$$ \begin{cases}
T_a ^{(2)} (T_s) < 2^{-1/4} T_s &\text{ if } T_s < T_{s,1} ^* , \\
T_a ^{(2)} (T_s) = 2^{-1/4} T_s &\text{ if } T_s = T_{s,1} ^*, \\
T_a ^{(2)} (T_s) > 2^{-1/4} T_s &\text{ if } T_s > T_{s,1} ^* .
\end{cases} $$
Moreover, $T_a ^{(2)} (T_s)=0$ if $\sigma _B \vert T_s \vert ^3 T_s = q\beta _s (T_s)$. The latter equality has at least once solution and at most three, as we have seen in Section \ref{sub-sec-eq123}. More precisely, in Section \ref{sub-sec-eq123}, we studied the intersections between the graph of $T_s \mapsto q\beta _s (T_s)$ and the graph of $ T_s \mapsto \sigma _B (1 - \frac{\varepsilon _a}{2}) T_s ^4$, while here we are interested in the intersections between the graph of $T_s \mapsto q\beta _s (T_s)$ with the graph of $ T_s \mapsto \sigma _B T_s ^4$. However, the argument of convexity still applies and thus we can have one, two or three intersections (depending on the values of $\sigma _B$, $q$, and the parameters appearing in $\beta _s$). Therefore, we shall analyse these different cases. Finally, for $T_s$ large, we have
$$ T_a ^{(2)} (T_s) \sim \varepsilon _a ^{-1/4} T_s \quad \text{ as } T_s \to \infty .$$
Let us introduce the following subdomains of $\mathcal Q$:
\begin{equation*}
\mathcal Q_1 := \{ (T_a, T_s) \in \mathcal Q,\quad
T_s > 2^{1/4} T_a\quad
\text{and}\quad
- \sigma _B T_s ^4  + \varepsilon _a \sigma _B T_a ^4 + q\beta _s (T_s) > 0 
\} ,
\end{equation*}
\begin{equation*}
\mathcal Q_2 := \{ (T_a, T_s) \in \mathcal Q,\quad
T_s > 2^{1/4} T_a\quad
\text{and}\quad
- \sigma _B T_s ^4  + \varepsilon _a \sigma _B T_a ^4 + q\beta _s (T_s) < 0  \} ,
\end{equation*}
\begin{equation*}
\mathcal Q_3 := \{ (T_a, T_s) \in \mathcal Q,\quad
T_s < 2^{1/4} T_a\quad
\text{and}\quad
- \sigma _B T_s ^4  + \varepsilon _a \sigma _B T_a ^4 + q\beta _s (T_s) < 0 \} ,
\end{equation*}
\begin{equation*}
\mathcal Q_4 := \{ (T_a, T_s) \in \mathcal Q,\quad
T_s < 2^{1/4} T_a\quad
\text{and}\quad
- \sigma _B T_s ^4  + \varepsilon _a \sigma _B T_a ^4 + q\beta _s (T_s) > 0 \} .
\end{equation*}
Observe that $\mathcal Q_1$, $\mathcal Q_2$, $\mathcal Q_3$ and $\mathcal Q_4$ are open subdomains of $\mathcal Q$, and they are separated by $\mathcal C_1$ and $\mathcal C_2$.

Figure \ref{fig4:region} below represents in phase space the case of one solution, $T_{s,1}$, of the equation $\sigma _B \vert T_s \vert ^3 T_s = q\beta _s (T_s)$. Whereas, in Figure \ref{fig4:region-ter} we sketch the sets $\mathcal{Q}_1$, $\mathcal{Q}_2$, $\mathcal{Q}_3$ and $\mathcal{Q}_4$ when $\sigma _B \vert T_s \vert ^3 T_s = q\beta _s (T_s)$
admits three solutions, denoted by $T_{s,1}$, $T_{s,2}$ and $T_{s,3}$.

\begin{figure}[h!]
\centering
\begin{tikzpicture}[scale=.6]
\draw[-stealth] (-1,0) -- (7,0)node[below]{$T_a$};
\draw[-stealth] (0,-1) -- (0,7)node[left]{$T_s$};
\draw[blue] (0,0) -- (7,7)node[right]{$\mathcal C_1$};
\draw[red] (0,1.25) to[out=60,in=200] (3,3) to (7,{(1/sqrt(3))*7});
\draw[red,-stealth] (.625,2) -- (1.125,2);
\draw[red,-stealth] (1.525,2.5) -- (2.025,2.5);
\draw[blue,-stealth] (1,1) -- (1,1.5);
\draw[blue,-stealth] (1.5,1.5) -- (1.5,2);
\node at (.5,1.15) {$Q_1$};
\draw[red,-stealth] (4.5,3.389) -- (4,3.389);
\draw[red,-stealth] (5.5,3.65) -- (5,3.65);
\draw[blue,-stealth] (5,5) -- (5,4.5);
\draw[blue,-stealth] (6,6) -- (6,5.5);
\node at (5.5,4.25) {$Q_3$};
\node at (2,4.5) {$Q_2$};
\draw[dashed,green!50!black] (0,3)node[left]{$T_{s,1} ^*$} -- (7,3);
\draw (-.1,1.3) node[left]{$T_{s,1}$}-- (.1,1.3);
\node[red] at (7.8,4.2) {$\mathcal C_2$};
\node at (4,1.5) {$Q_4$};
\end{tikzpicture}
\caption{In the phase space we represent the case of one equilibrium point $(2^{-1/4}T^*_{s,1}T^*_{s,1})$, which is solution of (6.1), and therefore intersection of the curves $\mathcal{C}_1\,:\,T_s=2^{1/4}T_a$ and the one defined by the set $\mathcal{C}_2=\{(T_a,T_s)\in\RR^2\,:\,-\sigma_B|T_s|^3T_s+\eps_a\sigma_B|T_a|^3T_a+q\beta_s(T_s)=0\}$. In particular, here we consider the case of one solution of the equation $\sigma_BT^4_s=q\beta_s(T_s)$. We subdivide the set $\mathcal{Q}=\{(T_a,T_s)\,:\, T_a>0,\,T_s>0\}$ in the subsets $\mathcal{Q}_1$, $\mathcal{Q}_2$, $\mathcal{Q}_3$ and $\mathcal{Q}_4$. We describe with arrows the behaviour of the vector field on the boundaries of $\mathcal{Q}_1$, $\mathcal{Q}_2$, $\mathcal{Q}_3$ and $\mathcal{Q}_4$.
}\label{fig4:region}
\end{figure}
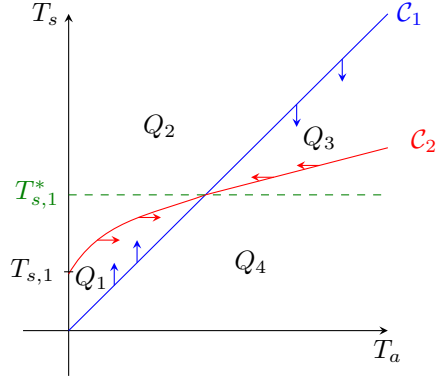

\begin{figure}[h!]
\centering

\begin{tikzpicture}[scale=.6]
\draw[-stealth] (-1,0) -- (7,0)node[below]{$T_a$};
\draw[-stealth] (0,-1) -- (0,7)node[left]{$T_s$};
\draw[blue] (0,0) -- (6,7)node[right]{$\mathcal C_1$};
\draw[red](0,3) -- (6.5,5);
\node[red] at (7.2,5.2) {$\mathcal C_2$};
\draw[red,-stealth] (1.65,3.5) -- (2.15,3.5);
\draw[red,-stealth] (.65,3.2) -- (1.15,3.2);
\draw[red,-stealth] (5,4.54) -- (4.5,4.54);
\draw[red,-stealth] (6,4.85) -- (5.5,4.85);
\draw[red,-stealth] (.38,2.3) -- (.88,2.3);
\draw[blue,-stealth] (2,2.32) -- (2,2.82);
\draw[blue,-stealth] (1.5,1.75) -- (1.5,2.25);
\draw[blue,-stealth] (4.5,5.25) -- (4.5,4.75);
\draw[blue,-stealth] (5,5.82) -- (5,5.32);
\draw[red] (0,1) arc (-90: 90: .7);
\node at (.3,1.7) {$Q_2$};
\node at (1.2,2.7) {$Q_1$};
\node at (2.5,5.5) {$Q_2$};
\node at (6,5.5) {$Q_3$};
\node at (4.5,1.7) {$Q_4$};
\draw (-.1,1) node[left]{$T_{s,1}$}-- (.1,1);
\draw (-.1,2.4) node[left]{$T_{s,2}$}-- (.1,2.4);
\draw (-.1,3) node[left]{$T_{s,3}$}-- (.1,3);
\draw[dashed,green!50!black] (0,4.1)node[left]{$T_{s,1} ^*$} -- (7,4.1);
\end{tikzpicture}
\caption{In the phase space we represent the case of one equilibrium point $(2^{-1/4}T^*_{s,1}T^*_{s,1})$, which is solution of (6.1), and therefore intersection of the curves $\mathcal{C}_1\,:\,T_s=2^{1/4}T_a$ and the one defined by the set $\mathcal{C}_2=\{(T_a,T_s)\in\RR^2\,:\,-\sigma_B|T_s|^3T_s+\eps_a\sigma_B|T_a|^3T_a+q\beta_s(T_s)=0\}$. In particular, here we consider the case of three solutions of the equation $\sigma_BT^4_s=q\beta_s(T_s)$. We subdivide the set $\mathcal{Q}=\{(T_a,T_s)\,:\, T_a>0,\,T_s>0\}$ in the subsets $\mathcal{Q}_1$, $\mathcal{Q}_2$, $\mathcal{Q}_3$ and $\mathcal{Q}_4$. 
We describe with arrows the behaviour of the vector field on the boundaries of $\mathcal{Q}_1$, $\mathcal{Q}_2$, $\mathcal{Q}_3$ and $\mathcal{Q}_4$.
}
\label{fig4:region-ter}
\end{figure}
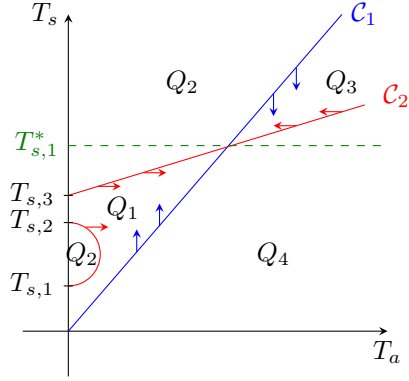

We are now going to study the behaviour of the solutions of \eqref{2layer-pbm-EDO} when the initial condition varies on $\mathcal{Q}$. We consider the case of one solution of $\sigma _B \vert T_s \vert ^3 T_s = q\beta _s (T_s)$. The remain case can be treated analogously.

\begin{figure}[h!]
\centering

\begin{tikzpicture}[scale=0.6]
\draw[-stealth] (-1,0) -- (7,0)node[below]{$T_a$};
\draw[-stealth] (0,-1) -- (0,7)node[left]{$T_s$};
\draw[blue] (0,0) -- (7,7);
\draw[red] (0,1.25) to[out=60,in=200] (3,3) to (7,{(1/sqrt(3))*7});
\draw[->-,dashed,brown!50!black] (1,5) -- (3,3) ;
\draw[->-,dashed,brown!50!black] (5,1) -- (3,3);
\draw[->-,brown!50!black] (.5,1.25) -- (3,3);
\draw[->-,brown!50!black] (5.5,{(7/10)*5+1.25}) -- (3,3);
\draw[->-,brown!50!black] (2,.5) to[out=140,in=240] (1.65,1.85);
\draw[->-,brown!50!black] (.5,2.5) to[out=-30,in=210] (1.75,2.25);
\draw[->-,brown!50!black] (4,6.5) to[out=-45,in=45] (4.5,4.25);
\draw[->-,brown!50!black] (6,1.5) to[out=110,in=30] (3.5,3.25);
\draw[dashed,green!50!black] (0,3)node[left]{$T_{s,1} ^*$} -- (7,3);
\node[red] at (10,4.5) {$\eps_a\sigma_BT_a^4 = \sigma_BT_s^4-q\beta _s (T_s)$};
\node[blue] at (8.8,7.2) {$T_s = 2^{1/4}T_a$};
\end{tikzpicture}
\caption{In the phase space we sketch by black arrows the convergence of initial conditions lying in each subset $\mathcal{Q}_1$, $\mathcal{Q}_2$, $\mathcal{Q}_3$ and $\mathcal{Q}_4$, determined by the the curves $\mathcal{C}_1\,:\,T_s=2^{1/4}T_a$ and the one defined by the set $\mathcal{C}_2=\{(T_a,T_s)\in\RR^2\,:\,-\sigma_B|T_s|^3T_s+\eps_a\sigma_B|T_a|^3T_a+q\beta_s(T_s)=0\}$, to the unique equilibrium point $(2^{-1/4}T^*_{s,1},T^*_{s,1})$, solution of (6.1).  In particular, here we consider the case of one solution of the equation $\sigma_BT^4_s=q\beta_s(T_s)$.
}
\label{fig4:region-comp}
\end{figure}
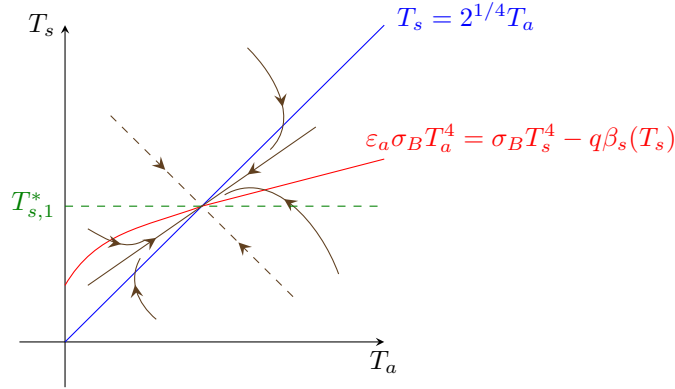

We have that:
\begin{itemize}
\item $(T_a ^{(0)}, T_s ^{(0)}) \in \mathcal Q_1$: the solution $(T_a,T_s)$ cannot leave $\mathcal Q_1$ since the vector field points inward the boundary. Therefore the solution cannot leave $\mathcal Q_1$. Furthermore, since $T_a ' >0$ and $T_s ' >0$ in $\mathcal{Q}_1$, then $T_a$ and $T_s$ are increasing. However, $T_a$ and $T_s$ are also bounded, which means that $(T_a,T_s)$ converges to the unique equilibrium point.

\item $(T_a ^{(0)}, T_s ^{(0)}) \in \mathcal Q_3$: also in this case the vector field points inward $\mathcal{Q}_3$. Thus, $(T_a,T_s)$ remains in $\mathcal{Q}_3$. Moreover, $T'_a<0$ and $T'_s<0$ in $\mathcal{Q}_3$ and so the solution converges to the equilibrium point;

\item $(T_a ^{(0)}, T_s ^{(0)}) \in \mathcal Q_2$: in this region $T_a$ is increasing and $T_s$ is decreasing. Therefore, we can have three possible scenarios:
\begin{itemize}
\item the solution never leaves $\mathcal Q_2$: it converges to the unique equilibrium point;
\item there exists a minimal value $\tau >0$ such that \\
$(T_a (\tau), T_s (\tau)) \in \overline{\mathcal Q_2} \cap \overline{\mathcal Q_1}$. In this case the vector field $F$ drives the solution inside $\mathcal Q_1$, and the solution converges monotonically (increasingly) to the equilibrium point;

\item there exists a minimal value $\tau >0$ such that such that \\
$(T_a (\tau), T_s (\tau)) \in \overline{\mathcal Q_2} \cap \overline{\mathcal Q_3}$. Thus the vector field $F$ drives the solution inside $\mathcal Q_3$ from which the solution converges monotonically (decreasingly) to the equilibrium point.

\end{itemize}

\item $(T_a ^{(0)}, T_s ^{(0)}) \in \mathcal Q_4$: the situation is similar to the previous case.
\end{itemize}
Therefore, the solution converges to the unique equilibrium point with some monotonicity: if the initial condition lies in $\mathcal Q_2$ or $\mathcal Q_4$, then the solution first enters in $\mathcal Q_1$ or $\mathcal Q_3$ and then it converges monotonically to the equilibrium (see Figure \ref{fig4:region-comp}).


\subsection{Monotonicity when there are three equilibrium points} \hfill
\label{sec-mono-3eq}

This is the most interesting case from a physical point of view. We assume that equation \eqref{equil-eqT0} has three solutions $T_{s,1} ^* \in (0,T_{s,-})$, $T_{s,2} ^* \in (T_{s,-}, T_{s,+})$, and $T_{s,3} ^* > T_{s,+}$, see Figure \ref{fig2:region}, case (A). Hence, problem \eqref{2layer-pbm-EDO} has three equilibrium points: $P^{eq} _1= (T_{a,1} ^*, T_{s,1} ^*)$, $P^{eq} _2= (T_{a,2} ^*, T_{s,2} ^*)$ and $P^{eq} _3=(T_{a,3} ^*, T_{s,3} ^*)$, with $T_{a,i} ^* = 2^{-1/4} T_{s,i} ^*$ ($i=1,2,3$).

Since the system is cooperative and the solutions are bounded, every solution will converge to an equilibrium point (see \cite{Smith}). We shall study in details the nature of such equilibria.

\subsubsection{Local stability of the equilibrium points}\hfill
\label{sec-local-stab}

First we look at the stability of the equilibrium points: 
since $\beta ' _s (T_{s,1} ^*)=0$, we have
$$  DF \left( \begin{array}{c} T_{a,1} ^* \\ T_{s,1} ^* \end{array} \right)
= \left( \begin{array}{cc} -\frac{8 \varepsilon _a \sigma _B}{\gamma _a} (T_{a,1} ^*)  ^3 & \frac{4 \varepsilon _a \sigma _B}{\gamma _a} (T_{s,1}  ^*)^3 \\
\frac{4 \varepsilon _a \sigma _B}{\gamma _s} (T_{a,1} ^*) ^3 & -\frac{4  \sigma _B}{\gamma _s} (T_{s,1} ^*) ^3 \end{array} \right) .$$
The trace of this matrix is negative, while the determinant is equal to 
\begin{multline*}\frac{32}{\gamma _a \gamma _s} \varepsilon _a \sigma _B ^2 (T_{a,1}^*) ^3 \, (T_{s,1} ^*) ^3 - \frac{16}{\gamma _a \gamma _s}\varepsilon _a ^2 \sigma _B ^2 (T_{s,1} ^*) ^3 \, (T_{a,1} ^*) ^3\\
= \frac{16}{\gamma _a \gamma _s}  (2 - \varepsilon _a) \varepsilon _a \sigma _B ^2 (T_{a,1} ^*) ^3 \, (T_{s,1} ^*) ^3 
\end{multline*}
and it is positive. Then, we conclude that the two eigenvalues have a negative real part which implies that the equilibrium point $P^{eq} _1$ is asymptotically exponentially stable.

The same can be proved for the equilibrium point $P^{eq} _3$. Concerning the equilibrium point $P^{eq} _2$, we compute the Jacobian matrix at $(T_{a,2} ^*,T_{s,2} ^*)$:
\begin{equation*}
DF \left( \begin{array}{c} T_{a,2} ^* \\ T_{s,2} ^* \end{array} \right)
= \left( \begin{array}{cc} -\frac{8 \varepsilon _a \sigma _B}{\gamma _a} (T_{a,2} ^*) ^3 & \frac{4 \varepsilon _a \sigma _B}{\gamma _a} (T_{s,2} ^*) ^3 \\
\frac{4 \varepsilon _a \sigma _B}{\gamma _s} (T_{a,2} ^*) ^3 & \frac{1}{\gamma _s} [-4 \sigma _B  (T_{s,2} ^*) ^3 + q\beta _s '(T_{s,2} ^*)] \end{array} \right) .
\end{equation*}
The determinant of $DF(T_{a,2} ^*,T_{s,2} ^*)$ is given by
\begin{equation*} 
\begin{split}
\gamma _a \gamma _s\det DF \begin{pmatrix} T_{a,2} ^*  \\ T_{s,2} ^* \end{pmatrix}&=16  (2 - \varepsilon _a) \varepsilon _a \sigma _B ^2 (T_{a,2}^*)  ^3 \, (T_{s,2} ^*) ^3
- 8 \varepsilon _a \sigma _B (T_{a,2} ^*)^3 \, q\beta _s '(T_{s,2}^*)
\\
&= 8 \varepsilon _a \sigma _B (T_{a,2} ^*) ^3 \, \Bigl( \frac{d}{d T_s} \Bigl[ (1-\frac{\varepsilon _a}{2})\sigma _B T_s ^4 - q\beta _s(T_s) \Bigr] \Bigr) _{\vert T_s = T_{s,2}^*}  .
\end{split}
\end{equation*}
The function $T_s \mapsto (1-\frac{\varepsilon _a}{2}) \sigma _B T_s ^4 - q\beta (T_s)$ is equal to $0$ at $T_{s,2}^*$,
and its derivative at $T_{s,2}^*$ is negative when there are $3$ equilibrium points (otherwise, by convexity, there would not be a third intersection relative to $T_{s,3}$). Therefore $\det DF (P^{eq} _2) <0$, 
and this implies that the equilibrium point $P^{eq} _2$ is unstable.

\subsubsection{The decomposition in subdomains} \hfill

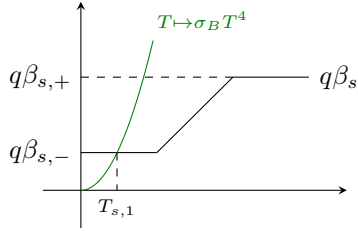
\begin{figure}[h!]
\centering
\begin{tikzpicture}[scale=.5]
\draw[-stealth] (-1,0) -- (7,0);
\draw[-stealth] (0,-1) -- (0,5);
\draw (0,1)node[left]{$q\beta_{s,-}$} -- (2,1);
\draw (2,1) -- (4,3);
\draw (4,3) -- (6,3);
\draw[dashed] (0,3)node[left]{$q\beta_{s,+}$} -- (4,3);
\draw[domain=0:1.9, smooth, variable=\x,green!50!black] plot ({\x}, {1.1*\x*\x});
\node[green!50!black] at (3.2,4.5) {$\scriptstyle{T\mapsto \sigma_BT^4}$};
\draw[dashed] (.95,1) -- (.95,0)node[below]{$\scriptstyle{T_{s,1}}$};
\node at (6.8,3) {$q\beta _s$};
\end{tikzpicture}
\caption{In the picture we show the case of a unique intersection between the curves $T\mapsto \sigma_BT^4$ and $T\mapsto q\beta_s(T)$. In particular, the solution of the equation $\sigma_BT^4_s=q\beta_s(T_s)$, denoted by $T_{s,1}$, is attained for $T<T_{s,-}$.
}\label{fig9}
\end{figure}

Let us consider again the line 
$$\mathcal C_1: T_s = 2^{1/4} T_a$$ and the curve
$$ \mathcal C_2:= \{(T_a, T_s), - \sigma _B \vert T_s \vert ^3 T_s  + \varepsilon _a \sigma _B \vert T_a \vert ^3 T_a + q\beta _s (T_s) =0 \}.$$
We recall that equation
$$ \sigma _B \vert T_s \vert ^3 T_s = q\beta _s (T_s) ,$$
can have at most 3 solutions, denoted by $T_{s,i}$. Let us assume that there is only one solution, $T_{s,1}$, as in Figure \ref{fig9}.

Let $\psi$ be the function defined in \eqref{def-psi}, and $T_a ^{(2)} (T_s)$ be the unique value for which
$\psi (T_a ^{(2)} (T_s)) =0$, that is, $(T_a,T_s) \in \mathcal C_2$. Note that
\begin{equation*}
\begin{split}
    T_s < T_{s,1} &\iff T_a ^{(2)} (T_s) < 0, \\
T_s = T_{s,1} &\iff T_a ^{(2)} (T_s) = 0, \\
T_s > T_{s,1} &\iff T_a ^{(2)} (T_s) > 0 .
\end{split}
\end{equation*}

Moreover, as we observed in section \ref{sec-mono-1eq}, it holds that
$$ \psi (2^{-1/4} T_s)
= q \beta _s (T_s) - \sigma _B (1 - \frac{\varepsilon _a}{2}) T_s ^4 .$$ 
Therefore
 \begin{itemize}
\item if $T_s \in (0,T_{s,1}^*)$ then $\psi (2^{-1/4} T_s) >0$ (see Figure \ref{fig2:region} (A)). Hence $T_a ^{(2)} (T_s) < 2^{-1/4} T_s$ and we deduce that $\mathcal C_2$ is on the left of the line $\mathcal C_1$;
\item if $T_s \in (T_{s,1}^*, T_{s,2}^*)$ we obtain that $T_a ^{(2)} (T_s) > 2^{-1/4} T_s$, and thus $\mathcal C_2$ is on the right of the line $\mathcal C_1$;
\item if $T_s \in (T_{s,2}^*, T_{s,3}^*)$ then $T_a ^{(2)} (T_s) < 2^{-1/4} T_s$. Therefore $\mathcal C_2$ is on the left of the line $\mathcal C_1$;
\item if $T_s >T_{0,3}^*$ we get that $T_a ^{(2)} (T_s) > 2^{-1/4} T_s$ which means that $\mathcal C_2$ is on the right of the line $\mathcal C_1$.
\end{itemize}
Let us now consider the following subdomains of $\mathcal Q$:
\begin{equation*}
    \mathcal Q_1 := \{ (T_a, T_s) \in \mathcal Q\,:\,
T_s \in (0, T_{s,1}^*), \,
T_s > 2^{1/4} T_a ,\,
- \sigma _B T_s ^4  + \varepsilon _a \sigma _B T_a ^4 + q\beta _s(T_s) > 0  \} ,
\end{equation*}

\begin{equation*}
\mathcal Q_1 ' := \{ (T_a, T_s) \in \mathcal Q\,:\,
T_s \in (T_{s,1}^*, T_{s,2}^*),\,
T_s < 2^{1/4} T_a ,
- \sigma _B T_s ^4  + \varepsilon _a \sigma _B T_a ^4 + q\beta _s(T_s) < 0  \} ,
\end{equation*}

\begin{equation*}
\mathcal Q_2 := \{ (T_a, T_s) \in \mathcal Q\,:\,
T_s > 2^{1/4} T_a,
- \sigma _B T_s ^4  + \varepsilon _a \sigma _B T_a ^4 + q\beta  _s(T_s) < 0 \} ,
\end{equation*}

\begin{equation*}
\mathcal Q_3 := \{ (T_a, T_s) \in \mathcal Q\,:\,
T_s >T_{s,3}^*,\,
T_s < 2^{1/4} T_a, 
- \sigma _B T_s ^4  + \varepsilon _a \sigma _B T_a ^4 + q\beta _s (T_s) < 0 \} ,
\end{equation*}

\begin{equation*}
\mathcal Q_3 ' := \{ (T_a, T_s) \in \mathcal Q\,:\,
T_s \in (T_{s,2}^*, T_{s,3}^* ),\,
T_s > 2^{1/4} T_a,
- \sigma _B T_s ^4  + \varepsilon _a \sigma _B T_a ^4 + q\beta _s (T_s) > 0 \} ,
\end{equation*}

\begin{equation*}
\mathcal Q_4 := \{ (T_a, T_s) \in \mathcal Q\,:\,
T_s < 2^{1/4} T_a,
- \sigma _B T_s ^4  + \varepsilon _a \sigma _B T_a ^4 + q\beta _s (T_s) > 0 \} 
\end{equation*}
which are sketched in Figure \ref{fig8:region} below.

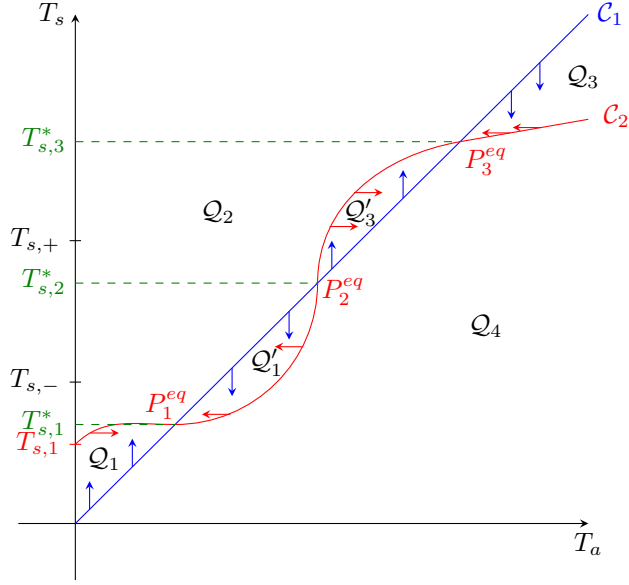
\begin{figure}[h!]
\centering

\begin{tikzpicture}[scale=0.75]
\draw[-stealth] (-1,0) -- (9,0)node[below]{$T_a$};
\draw[-stealth] (0,-1) -- (0,9)node[left]{$T_s$};
\draw[blue] (0,0) -- (9,9)node[right]{$\mathcal C_1$};
\draw[red] (0,1.4) to[in=180] (1.75,1.75);
\draw[red] (1.75,1.75) to[out=0,in=-90] (4.25,4.25);
\draw[red] (4.25,4.25) to[out=90,in=190] (6.75,6.75);
\draw[red] (6.75,6.75) -- (9,{.176*(9-6.75)+6.75});
\draw[red,-stealth] (.265,1.6) -- (.765,1.6);
\draw[blue,-stealth] (1,1) -- (1,1.5);
\draw[blue,-stealth] (.25,.25) -- (.25,.75);
\node at (.5,1.15) {$\mathcal Q_1$};
\draw[red,-stealth] (4,3.125) -- (3.5,3.125);
\draw[blue,-stealth] (3.75,3.75) -- (3.75,3.25);
\draw[blue,-stealth] (2.75,2.75) -- (2.75,2.25);
\draw[red,-stealth] (2.72,1.93) -- (2.22,1.93);
\draw[blue,-stealth] (4.5,4.5) -- (4.5,5);
\draw[blue,-stealth] (5.75,5.75) -- (5.75,6.25);
\draw[blue,-stealth] (7.65,7.65) -- (7.65,7.15);
\draw[blue,-stealth] (8.15,8.15) -- (8.15,7.65);
\draw[red,-stealth] (4.46,5.25) -- (4.96,5.25);
\draw[red,-stealth] (4.89,5.85) -- (5.39,5.85);
\draw[red,-stealth] (8.173,7) -- (7.673,7);
\draw[red,-stealth] (7.63,6.905) -- (7.13,6.905);
\draw[green!50!black,dashed] (0,1.75)node[left]{$T_{s,1}^*$} -- (1.75,1.75);
\draw[green!50!black,dashed] (0,4.25)node[left]{$T_{s,2}^*$} -- (4.25,4.25);
\draw[green!50!black,dashed] (0,6.75)node[left]{$T_{s,3}^*$} -- (6.75,6.75);
\node at (5,5.5) {$\mathcal Q_3 '$};
\node at (8.9,7.9) {$\mathcal Q_3$};
\node at (2.5,5.5) {$\mathcal Q_2$};
\node at (3.35,2.85) {$\mathcal Q_1 '$};
\draw (-.1,2.5) node[left]{$T_{s,-}$}-- (.1,2.5);
\draw (-.1,5) node[left]{$T_{s,+}$}-- (.1,5);
\node at (7.2,3.5) {$\mathcal Q_4$};
\draw[red] (9.5,7.2)node{$\mathcal C_2$};
\draw[red] (-.1,1.4) node[left]{$T_{s,1}$}-- (.1,1.4);
\node[red] at (1.6,2.1){$P_1 ^{eq}$};
\node[red] at (4.7,4.1){$P_2 ^{eq}$};
\node[red] at (7.2,6.4){$P_3 ^{eq}$};
\end{tikzpicture}
\caption{In the phase space we consider the case of three intersections between the curves $\mathcal{C}_1\,:\,T_s=2^{1/4}T_a$ and the one defined by the set $\mathcal{C}_2=\{(T_a,T_s)\in\RR^2\,:\,-\sigma_B|T_s|^3T_s+\eps_a\sigma_B|T_a|^3T_a+q\beta_s(T_s)=0\}$. We denote by $(2^{-1/4}T^*_{s,1}, T^*_{s,1})$, $(2^{-1/4}T^*_{s,2}, T^*_{s,2})$ and $(2^{-1/4}T^*_{s,3}, T^*_{s,3})$ the three equilibria and by $\mathcal{Q}_1$, $\mathcal{Q}'_1$, $\mathcal{Q}_2$, $\mathcal{Q}_3$, $\mathcal{Q}'_3$, $\mathcal{Q}_4$ the subsets of $\mathcal{Q}=\{(T_a,T_s)\,:\, T_a>0,\,T_s>0\}$ bordered by $\mathcal{C}_1$ and $\mathcal{C}_2$. We describe with arrows the behaviour of the vector field on the boundaries of the aforementioned subsets. Notice that we are considering the case in which the equation $\sigma_BT^4_s=q\beta_s(T_s)$ has a unique solution.
}
\label{fig8:region}
\end{figure}

As before, we look at the vector field $F$ on the boundary of each subdomain. On the line $\mathcal C_1$, we have
$$ 
F \begin{pmatrix} T_a = 2^{-1/4} T_s\\ T_s \end{pmatrix} 
= \begin{pmatrix} 0 \\ \frac{1}{\gamma _s} \Bigl[ q\beta _s (T_s) - \sigma _B (1- \frac{\varepsilon _a}{2}) \vert T_s \vert ^3 T_s  \Bigr] \end{pmatrix}.
$$
Therefore, $ F \begin{pmatrix} T_a = 2^{-1/4} T_s\\ T_s \end{pmatrix}$ is vertical and its second component is positive if and only if $T_s \in (0, T_{s,1}^*) \cup (T_{s,2}^*, T_{s,3}^*)$, and negative on $(T_{s,1}^*, T_{s,2}^*) \cup (T_{s,3}^*, +\infty)$.

On the curve $\mathcal C_2$ we have
$$
F \begin{pmatrix} T_a \\ T_s \end{pmatrix} 
= \begin{pmatrix} \frac{1}{\gamma _a} \Bigl[\varepsilon _a \sigma _B \vert T_s \vert ^3 T_s  - 2 \varepsilon _a \sigma _B \vert T_a \vert ^3 T_a  \Bigr] \\ 0 \end{pmatrix}.
$$
Hence in this case $ F \begin{pmatrix} T_a \\ T_s \end{pmatrix}$ is horizontal. It heads right if and only if the point $(T_a,T_s)$ is above the line $\mathcal C_1$ whereas it heads left if and only if $(T_a,T_s)$ is below $\mathcal C_1$. This completes the explanation of Figure \ref{fig8:region}.

\subsubsection{Monotonicity and convergence} \hfill

The easiest cases to analyse are the following
\begin{itemize}
\item $(T_a ^{(0)}, T_s ^{(0)}) \in \mathcal Q_1$: the solution does not leave $Q_1$, the functions $T_a$ and $T_s$ are increasing and the solution $(T_a, T_s)$ converges to the equilibrium point $P^{eq} _1$; 
\item $(T_a ^{(0)}, T_s ^{(0)}) \in \mathcal Q_1 '$: the solution does not leave $\mathcal Q_1 '$, the functions $T_a$ and $T_s$ are decreasing and the solution $(T_a, T_s)$ converges to the equilibrium point $P^{eq} _1$; 
\item $(T_a ^{(0)}, T_s ^{(0)}) \in \mathcal Q_3 '$: the solution does not leave $\mathcal Q_3^{'}$, the functions $T_a$ and $T_s$ are increasing and the solution $(T_a, T_s)$ converges to the equilibrium point $P^{eq} _3$; 
\item $(T_a ^{(0)}, T_s ^{(0)}) \in \mathcal Q_3$: the solution does not leave $\mathcal Q_3$, the functions $T_a$ and $T_s$ are decreasing and the solution $(T_a, T_s)$ converges to the equilibrium point $P^{eq} _3$.
\end{itemize}
It remains to study the behaviour of the solution when the initial condition belongs either to $\mathcal Q_2$ or  $\mathcal Q_4$. Let us define
\begin{equation}
\label{Cleft}
\mathcal C _{\text{left}} := 
\partial \mathcal Q_2 \cap \mathcal Q ,
\quad \text{ and } \quad 
\mathcal C _{\text{right}} := \partial \mathcal Q_4 \cap \mathcal Q .
\end{equation}
Consider an initial condition in the region $\mathcal Q_4$. As long as the solution is in $Q_4$, $T_a$ decreases and $T_s$ increases. Therefore, either the solution does not leave $Q_4$, and it converges monotonically to some equilibrium point, or it reaches $\mathcal C _{\text{right}}$. The latter case, however, cannot happen at an equilibrium point. Hence, either the solution attains $\partial \mathcal Q_1$, enters $\mathcal{Q}_1$ and converges increasingly to $P_1 ^{eq}$, or it gets to $\partial \mathcal Q_1 '$, enters $\mathcal Q_1 '$ and it converges decreasingly to $P_1 ^{eq}$, or it move to $\partial \mathcal Q_3 '$ enters $\mathcal Q_3 '$ and it converges increasingly to $P_3 ^{eq}$, or it reaches $\partial \mathcal Q_3$, enters $\mathcal Q_3$ and it converges decreasingly to $P_3 ^{eq}$.

Analogously, if the initial condition is of the form $(T_a ^{(0)}, 0)$ with $T_a ^{(0)} >0$, then the solution immediately enters $\mathcal Q_4$ and it behaves as explained above.

A similar argument can be adopted for initial conditions lying on $\mathcal Q_2$, or of the form $(0,T_s ^{(0)})$ with $ T_s ^{(0)} \geq 0$. Hence, we have a complete description of the behaviour of the solution with initial condition in $\overline{\mathcal Q}$.

Now, we study what happens backward in time. We are going to prove the following

\begin{Lemma}
\label{lem-derdeder?}
Given an initial condition in $\mathcal C _{\text{right}}$ that is not an equilibrium point, there exists $\tau >0$ such that $T_s (-\tau) =0$ and $T_a (-\tau) >0$. Therefore, the solution at time $-\tau$ is on the horizontal axis.
\end{Lemma}

\noindent {\it Proof of Lemma \ref{lem-derdeder?}.} Consider an initial condition lying in $\mathcal C _{\text{right}}$ that is not an equilibrium point, and let the time go backward. Then, the vector field $F$ drives the solution towards $\mathcal Q_4$ so that the solution cannot reach again $\mathcal C _{\text{right}}$. Moreover, when time goes backward, $T_a$ increases and $T_s$ decreases, hence the solution moves to the south-east direction. Therefore, either there exists $\tau >0$ such that $T_s(-\tau) =0$, or $T_s$ remains positive. We shall prove that the latter possibility cannot happen.

By contradiction, assume that $T_s$ remains positive. If $T_a$ is bounded, then the solution has to converge to some equilibrium point. However, such behaviour is not compatible with the monotonicity of the solution that moves in the "south-east direction". Therefore $T_a$ must be unbounded and diverge to $+\infty$. In this case \eqref{2layer-pbm-EDO} implies that also $T_s '$ diverges to infinity.
Moreover,
$$ \frac{\gamma _s}{\gamma _a} \frac{T_s '}{T_a '} = \frac{- \sigma _B \vert T_s \vert ^3 T_s  + \varepsilon _a \sigma _B \vert T_a \vert ^3 T_a + \mathcal R_s}{\varepsilon _a \sigma _B \vert T_s \vert ^3 T_s  - 2 \varepsilon _a \sigma _B \vert T_a \vert ^3 T_a} \to -\frac{1}{2} \quad \text{ as } T_a \to +\infty ,$$
therefore
$$ - \frac{3}{4} \frac{\gamma _a}{\gamma _s} T_a ' \leq T_s ' \leq - \frac{1}{4} \frac{\gamma _a}{\gamma _s} T_a ' $$
for $T_a$ large enough. By integration, we deduce that $T_s\to-\infty$, which contradicts the assumption.
%
%
Therefore, we conclude that there exists $\tau >0$ such that $T_s(-\tau) =0$, and Lemma \ref{lem-derdeder?} is proved. \qed

We use the above result to prove the following
\begin{Lemma}
\label{lem-derdeder}

There exists a unique value $T_{a, \text{threshold}} >0$ such that
\begin{itemize} 
\item if $T_a ^{(0)} \in (0, T_{a, \text{threshold}})$, the solution starting from $(T_a ^{(0)},0)$ converges to $P_1 ^{eq}$,
\item if $T_a ^{(0)} = T_{a, \text{threshold}}$, the solution starting from $(T_a ^{(0)},0)$ converges to $P_2 ^{eq}$,
\item if $T_a ^{(0)} > T_{a, \text{threshold}}$, the solution starting from $(T_a ^{(0)},0)$ converges to $P_3 ^{eq}$.
\end{itemize}

\end{Lemma}

\noindent {\it Proof of Lemma \ref{lem-derdeder}.}
Consider
the subsets of initial conditions:
\begin{equation*}
\mathcal I _1 := \{ T_a ^{(0)} >0\,:\, \text{ the solution starting from $(T_a ^{(0)},0)$} \text{ converges to } P^{eq} _1 \},
\end{equation*}
\begin{equation*}
\mathcal I _3 := \{ T_a ^{(0)} >0\,:\, \text{ the solution starting from $(T_a ^{(0)},0)$} \text{ converges to } P^{eq} _3 \} .
\end{equation*}
Since two solutions with different initial conditions on the horizontal axis cannot cross each other because of the uniqueness, $\mathcal I _1$ and 
$\mathcal I _3$ are intervals. Moreover, we claim that $\mathcal I _1$ and $\mathcal I _3$ are open because $P_1 ^{eq}$ and $P_3 ^{eq}$ are asymptotically exponentially stable. Let us prove the latter property for $\mathcal{I}_1$. There exists $\eta _1 >0$ such that any solution starting from an initial condition in the ball $B(P_1 ^{eq}, \eta_1)$, with center $P_1 ^{eq}$ and radius $\eta_1$, converges to $P_1 ^{eq}$. If $T_a ^{(0)} \in \mathcal I _1$, there exists $\tau_0 >0$ such that the solution $(T_a,T_s)$ with initial condition $(T_a ^{(0)},0)$ is, at time $\tau _0$, in the ball $B(P_1 ^{eq}, \frac{\eta_1}{3})$. By continuity with respect to the initial condition (Gronwall's lemma), there exists $\eta _0 >0$ such that, if $\vert T_a ^{(0)} - \tilde T_a ^{(0)} \vert < \eta _0$, then the solution $(\tilde T_a, \tilde T_s)$ starting from $(\tilde T_a ^{(0)},0)$ satisfies
$$ \Vert (T_a (\tau_0), T_s (\tau_0)) - (\tilde T_a (\tau_0), \tilde T_s (\tau_0)) \Vert < \frac{\eta_1}{3};$$ This implies that $(\tilde T_a (\tau_0), \tilde T_s (\tau_0)) \in B(P_1 ^{eq}, \frac{2\eta_1}{3})$, and the solution $(\tilde T_a, \tilde T_s)$ converges to $P_1 ^{eq}$. This proves that $\mathcal I _1$ is open.

We observe that since $\mathcal I _1$ and $\mathcal I _3$ are open intervals contained in $(0,+\infty)$,  we cannot have $\mathcal I _1 \cup \mathcal I _3 = (0,+\infty)$.
This implies that there exists $T_{a, threshold}$ such that the solution starting from $(T_{a, \text{threshold}},0)$
converges neither to $P^{eq} _1$ nor to $P^{eq} _3$. The only possibility that remains is that this solutions does not leave $Q_4$, and therefore goes monotonically to $P^{eq} _2$. Since $P^{eq} _2$ is unstable
then we have that $\mathcal I _1 = (0, T_{a, \text{threshold}})$, and $\mathcal I _3 = (T_{a, \text{threshold}}, +\infty)$. This concludes the proof of Lemma \ref{lem-derdeder}. 
\qed

We observe that a similar argument applies for solutions starting from the vertical axis. Hence, also in this case there exists a threshold value $T_{s, \text{threshold}}$ such that the solution with initial condition of the form $(0,T^{(0)}_s)$, with $T^{(0)}_s<T_{s,\text{threshold}}$ converges to $P^1_{eq}$, whereas if $T^{(0)}_s>T_{s,\text{threshold}}$ the solution converges to $P^3_{eq}$. And finally, the solution starting in $(0,T_{s,\text{threshold}})$ converges to the unstable equilibrium $P^2_{eq}$.

Therefore, the trajectories of the solutions starting 
from $(T_{a, \text{threshold}},0)$ and  $(0, T_{s, \text{threshold}})$ merge in $P^{eq} _2$ and separate the quadrant $\mathcal Q$ into two subdomains that are the attraction basins of $P^{eq} _1$ and of $P^{eq} _3$, see Figure \ref{fig:basins}. 


\subsection{Weaker assumptions on the coalbedo $\beta_s$} \hfill
\label{sec-weaker-beta}

In this section we discuss the results obtained till now in presence of a more general function $\beta_s$. Propositions \ref{prop-wellposed} and \ref{prop-comp-lambda-eq-cv} are stated and proved just requiring global lipschitz continuity and positivity on $\beta _s$. The piecewise linear assumption \eqref{beta} has been used in Proposition \ref{prop-comp-lambda-eq-gene} mainly to determine the number of equilibria of system \eqref{2layer-pbm-EDO}, see section \ref{sub-sec-eq123}.

\subsubsection{Influence of the assumptions on $\beta _s$ on the number of equilibrium points} \hfill

Under weaker assumptions on $\beta_s$, equation \eqref{equil-eqT0} can have more than three solutions.
For instance, if we assume $\beta _s$ to be positive and to have to some limit for  $T_s \to +\infty$, then \eqref{equil-eqT0} can have several solutions on $[T_{s,+},+\infty)$. There can be even a continuum of solutions if $q\beta _s$ and $T_s \mapsto \sigma _B (1-\frac{\varepsilon _a}{2}) T_s ^4$ coincide on some compact interval. However, there are some quite general assumptions for $\beta_s$ that lead to a finite number of equilibrium points for problem \eqref{2layer-pbm-EDO}.

If we assume that $\beta _s$ is positive and converges to some finite limit as $T_s \to +\infty$ and additionally that $\beta _s$ is {\it analytic} on $[T_{s,+},+\infty)$, then equation \eqref{equil-eqT0} cannot have an infinite number of solutions on $[T_{s,+},+\infty)$. Indeed, if we consider the difference
\begin{equation}\label{difference}
    T_s \mapsto \sigma _B (1-\frac{\varepsilon _a}{2}) T_s ^4 - q \beta _s (T_s),
\end{equation}
which is an analytic function as well, it would have an infinite number of zeros contained in a compact set. This would imply that the above function equals $0$ on $[T_{s,+},+\infty)$, that is a contradiction with the behaviour at $+\infty$.

If we assume $\beta _s$ to be {\it concave} on $[T_{s,+},+\infty)$ in place of the analyticity on the same interval, then \eqref{difference} is strictly convex on $[T_{s,+},+\infty)$ and can have at most two zeros on $[T_{s,+},+\infty)$.

If instead we assume that $\beta _s$ is {\it analytic} and {\it positive} on $[0,T_{s,-}]$ and {\it concave} on $[T_{s,-},+\infty)$, then \eqref{difference} has a finite number of zeros. This fact can be proved combining the two previous arguments: the function has a finite number of zeros on $[0,T_{s,-}]$ by analyticity (otherwise it would be equal to zero, which is not true for $T_s=0$) and has a finite number of zeros on $[T_{s,-}, +\infty)$ by convexity.

Finally, let $\beta _s$ be  {\it analytic} and {\it positive} on $[0,T_{s,-}]$, {\it concave} on $[T_{s,-},T_{s,+}]$, {\it analytic} on $[T_{s,+},+\infty)$ and converge to some limit as $T_s \to +\infty$. Then \eqref{difference} has at most a finite number of zeros on $[0,T_{s,-}]$ and $[T_{s,+},+\infty)$, and at most 2 zeros on $[T_{s,-},T_{s,+}]$, hence a finite number of zeros on $[0,+\infty)$. These assumptions cover the case of a piece-wise linear function $\beta _s$ that we mainly use in the paper.


\subsubsection{The influence on the asymptotic analysis} \hfill

\begin{figure}[h!]
\centering
\subfloat{
\begin{tikzpicture}[scale=.65]
\draw[-stealth] (-1,0) -- (7,0);
\draw[-stealth] (0,-1) -- (0,7);
\draw (0,.5) -- (2.5,.5);
\draw (2.5,.5) -- (3,2.5);
\draw(3,2.5) to[out=70,in=-120](4.5,3);
\draw(4.5,3) to[out=60,in=-170](6,5.5);
\draw(6,5.5) -- (6.5,5.6);
\draw[domain=0:6.2, smooth, variable=\x,green!50!black] plot ({\x}, {.17*\x*\x});
\node[green!50!black] at (4.2,6.2) {$\scriptstyle{T\mapsto \sigma_B\left(1-\frac{\eps_a}{2}\right)T^4}$};
\draw[dashed] (1.7,.5) -- (1.7,0)node[below]{$\scriptstyle{T^*_{s,1}}$};
\draw[dashed] (2.7,1.2) -- (2.7,0)node[below]{$\quad \scriptstyle{T^*_{s,2}}$};
\draw[dashed] (4.05,2.75) -- (4.05,0)node[below]{$\quad \scriptstyle{T^*_{s,3}}$};
\draw[dashed] (4.9,4.1) -- (4.9,0)node[below]{$\quad \scriptstyle{T^*_{s,4}}$};
\draw[dashed] (5.65,5.4) -- (5.65,0)node[below]{$\quad \scriptstyle{T^*_{s,5}}$};
\end{tikzpicture}}\\
\subfloat{
\begin{tikzpicture}[scale=.3]
\draw[-stealth] (-1,0) -- (16,0)node[below]{$\scriptstyle{T_a}$};
\draw[-stealth] (0,-1) -- (0,15)node[left]{$\scriptstyle{T_s}$};
\draw[blue] (0,0) -- (12.5,12.5)node[left]{$\scriptstyle{T_s = 2^{\frac{1}{4}}T_a}\,\,$};
\draw[red,-stealth] (0,1.4) to[in=180] (1.75,1.75);
\draw[red,stealth-stealth] (1.75,1.75) to[out=0,in=-90] (4.25,4.25);
\draw[red,stealth-stealth] (4.25,4.25) to[out=90,in=190] (6.75,6.75);
\draw[red,stealth-stealth] (6.75,6.75) to[out=0,in=-90] (9,9);
\draw[red,stealth-stealth] (9,9) to[out=90,in=190] (11.25,11.25);
\draw[red] (11.25,11.25) -- (13.5,{.176*(2.25)+11.25});
\draw[red,-stealth] (.265,1.6) -- (.765,1.6);
\draw[blue,-stealth] (1,1) -- (1,1.5);
\draw[blue,-stealth] (.25,.25) -- (.25,.75);
\draw[red,-stealth] (4,3.125) -- (3.5,3.125);
\draw[blue,-stealth] (3.75,3.75) -- (3.75,3.25);
\draw[blue,-stealth] (2.75,2.75) -- (2.75,2.25);
\draw[red,-stealth] (2.72,1.93) -- (2.22,1.93);
\draw[blue,-stealth] (4.5,4.5) -- (4.5,5);
\draw[blue,-stealth] (5.75,5.75) -- (5.75,6.25);
\draw[blue,-stealth] (7.65,7.65) -- (7.65,7.15);
\draw[blue,-stealth] (8.15,8.15) -- (8.15,7.65);
\draw[red,-stealth] (4.46,5.25) -- (4.96,5.25);
\draw[red,-stealth] (4.89,5.85) -- (5.39,5.85);
\draw[red,-stealth] (7.8,7) -- (7.3,7);
\draw[red,-stealth] (8.8,8) -- (8.3,8);
\draw[red,-stealth] (9.25,10) -- (9.75,10);
\draw[red,-stealth] (10,10.8) -- (10.5,10.8);
\draw[blue,-stealth] (10,10) -- (10,10.5);
\draw[blue,-stealth] (9.3,9.3) -- (9.3,9.8);
\draw[blue,-stealth] (12,12) -- (12,11.5);
\draw[red,-stealth] (13.05,11.58) -- (12.55,11.58);
\draw[red,dashed] (0,1.75)node[left]{$\scriptstyle{T^*_{s,1}}$} -- (1.75,1.75);
\draw[red,dashed] (0,4.25)node[left]{$\scriptstyle{T^*_{s,2}}$} -- (4.25,4.25);
\draw[red,dashed] (0,6.75)node[left]{$\scriptstyle{T^*_{s,3}}$} -- (6.75,6.75);
\draw[red,dashed] (0,9)node[left]{$\scriptstyle{T^*_{s,4}}$} -- (9,9);
\draw[red,dashed] (0,11.25)node[left]{$\scriptstyle{T^*_{s,5}}$} -- (11.25,11.25);
\draw[red] (16,10.5)node{$\scriptstyle{\eps_a\sigma_B T_a^4 = \sigma_B T_s^4 - \beta_s(T_s)}$};
\end{tikzpicture}}
\caption{On the top a non-piecewise coalbedo which intersects five times the curve $T\mapsto \sigma_B\left(1-\frac{\eps_a}{2}\right)T^4$, and on the bottom the associated phase space.}
\label{fig:5eq-pts}
\end{figure}
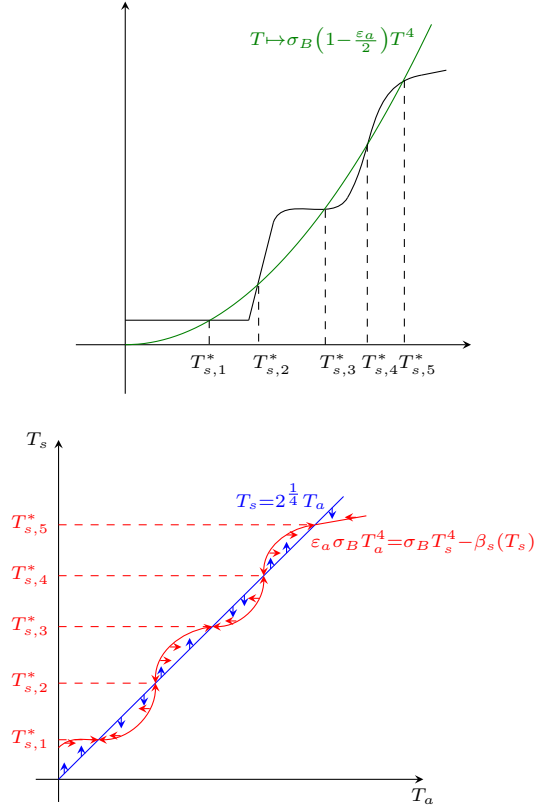

Once the number of equilibrium points have being determined, the asymptotic analysis remains essentially the same. Assume that we are in one of the situations described in the previous section: there is a finite number $N$ of solutions $T_{s,i} ^*$ of \eqref{equil-eqT0},
$\beta _s$ is $C^1$ at any $T_{s,i} ^*$ and
$$ \forall\, i \in \{1, \cdots, N \}, \quad \frac{d}{dT_s} _{\vert T_s = T_{s,i} ^*}  \Bigl( \sigma _B (1-\frac{\varepsilon _a}{2}) T_s ^4 \Bigr) \neq q \beta ' _s (T_{s,i} ^*) ,$$
that is, curves $T_s \mapsto \sigma _B (1-\frac{\varepsilon _a}{2}) T_s ^4$ and $T_s \mapsto q \beta _s (T_s)$ are not tangent on $T^*_{s,i}$ for all $i=1,\dots,N$. Then, under such assumptions, the sign of \eqref{difference} is alternatively negative and positive: negative on $[0,T_{s,1} ^*]$, positive on $[T_{s,1} ^*, T_{s,2} ^*]$, $\cdots$, positive on $[T_{s,N} ^*, +\infty )$.
This first implies that $N$ has to be odd. Furthermore, it forces $\mathcal C_2$ to be first on the left of $\mathcal C_1$, then, after the first equilibrium point, on the right, then, after the second equilibrium point, again on the left and so on till the last equilibrium point where it is on the right of $\mathcal{C}_1$. Finally, the assumptions on $\beta_s$ determine also the direction of the vector field $F$ on $\mathcal C_1$ and $\mathcal C_2$. $F$ is horizontal on $\mathcal C_2$ because on this curve its vertical component is equal to $0$. Moreover, if  $ T_s \in (T_{s,1} ^*, T_{s,2} ^*)$ the sign of the first component of $F$ is the same as the one of
\begin{equation*}
\varepsilon _a \sigma _B T_s ^4 - 2 \varepsilon _a \sigma _B T_a ^4 = 
\varepsilon _a \sigma _B T_s ^4 - 2 \Bigl( \sigma _B T_s ^4 - q \beta _s (T_s) \Bigr) = -2 \Bigl( \sigma _B (1-\frac{\varepsilon _a}{2}) T_s ^4 - q \beta _s (T_s) \Bigr) ,
\end{equation*}
which is negative on $(T_{s,1} ^*, T_{s,2} ^*)$. $F$ is vertical on $\mathcal C_1$ and if $ T_s \in (T_{s,1} ^*, T_{s,2} ^*)$ the sign of the second component of $F$ is the same of
$$ - \sigma _B T_s ^4 + \varepsilon _a \sigma _B T_a ^4 + q \beta _s (T_s) = - \sigma _B (1-\frac{\varepsilon _a}{2}) T_s ^4 + q \beta _s (T_s) $$
which is negative on $(T_{s,1} ^*, T_{s,2} ^*)$. Hence, if one consider, for instance, a coalbedo functions as the black curve in left picture of Figure \ref{fig:5eq-pts}, a similar argument leads to the phase space description (for $N=5$) on the right of Figure \ref{fig:5eq-pts}.

\section{Asymptotic behaviour for $\lambda>0$ and $\mathcal{R}_a=0$}
\label{sec7}

Equilibria of problem \eqref{2layer-pbm-EDO} are solutions $(T_a,T_s)$ of the following system:
\begin{equation}
\label{2layer-pbm-EDO-eq*}
\begin{cases}
-\lambda (T_a - T_s) + \varepsilon _a \sigma _B \vert T_s \vert ^3 T_s  - 2 \varepsilon _a \sigma _B \vert T_a \vert ^3 T_a =0 , \\
-\lambda (T_s - T_a) - \sigma _B \vert T_s \vert ^3 T_s  + \varepsilon _a \sigma _B \vert T_a \vert ^3 T_a + q \beta_s (T_s) = 0 .
\end{cases}
\end{equation}
Let us prove several properties of such points.


\subsection{Equilibrium points: uniform bounds}\label{sec-uniform-bounds}\hfill

Assume that $(T_a,T_s)$ is an equilibrium point of problem \eqref{2layer-pbm-EDO}, that is, a solution of \eqref{2layer-pbm-EDO-eq*}. Therefore $(T_a,T_s)$ solves
\begin{equation}
\label{2layer-pbm-EDO-eq-lambda}
\begin{cases}
\lambda T_s + \varepsilon _a \sigma _B \vert T_s \vert ^3 T_s   & = \lambda T_a + 2 \varepsilon _a \sigma _B \vert T_a \vert ^3 T_a, \\
\lambda T_s + \sigma _B \vert T_s \vert ^3 T_s
 &= \lambda T_a + \varepsilon _a \sigma _B \vert T_a \vert ^3 T_a + q\beta _s (T_s)
\end{cases}
\end{equation}
Observe that, from the first equation in \eqref{2layer-pbm-EDO-eq-lambda}, $T_a=0$ if and only if $T_s=0$. However, the pair $(T_a,T_s)=(0,0)$ does not satisfy the second equation. Thus, the equilibrium points have positive components.

Since $T_a >0$, we have
$$ \lambda T_a +  \varepsilon _a \sigma _B  T_a ^4 < \lambda T_a + 2 \varepsilon _a \sigma _B T_a ^4 < \lambda 2^{1/4} T_a + 2 \varepsilon _a \sigma _B T_a ^4 .$$
Therefore, from the first equation in \eqref{2layer-pbm-EDO-eq-lambda} we deduce that
$$ \lambda T_a +  \varepsilon _a \sigma _B  T_a ^4 < \lambda T_s + \varepsilon _a \sigma _B  T_s ^4 
< \lambda (2^{1/4} T_a) +  \varepsilon _a \sigma _B (2^{1/4}T_a )^4.$$
The map $T \mapsto \lambda T + \varepsilon _a \sigma _B  T ^4$ is increasing on $(0,+\infty)$, and this implies that
\begin{equation}
\label{lambda-posi}
T_a < T_s < 2^{1/4} T_a .
\end{equation}
We want to prove that the equilibrium points belong to a compact set independent of $\lambda$. To this purpose, we compute the difference of the two equations in \eqref{2layer-pbm-EDO-eq-lambda}, which gives
\begin{equation}
\label{lambda-diff}
\sigma _B (\varepsilon _a -1) T_s ^4 = \varepsilon _a \sigma _B T_a ^4 - q\beta _s (T_s) .
\end{equation}
If $\varepsilon _a \in (0,1]$, then
$ \varepsilon _a \sigma _B T_a ^4 - q\beta _s (T_s) \leq 0$,
and therefore 
\begin{equation}
\label{lambda-01-sup}
T_a \leq \Bigl( \frac{q\beta _{s,+}}{\varepsilon _a \sigma _B} \Bigr) ^{1/4}.
\end{equation}
Since $T_a < T_s$, we also have
\begin{equation*}
    q\beta _s (T_s ) 
= \varepsilon _a \sigma _B T_a ^4 
+ \sigma _B (1-\varepsilon _a) T_s ^4  <\varepsilon _a \sigma _B T_s ^4 
+ \sigma _B (1-\varepsilon _a) T_s ^4  
= \sigma _B T_s ^4 ,
\end{equation*} 
from which we get
\begin{equation}
\label{lambda-01-inf}
T_s > \Bigl( \frac{q\beta _{s,-}}{\sigma _B} \Bigr) ^{1/4}.\end{equation}
Thus, from \eqref{lambda-posi}, \eqref{lambda-01-sup} and \eqref{lambda-01-inf} we obtain
\begin{equation*}
    \Bigl( \frac{q\beta _{s,-}}{2\sigma _B} \Bigr) ^{1/4} < T_a \leq \Bigl( \frac{q\beta _{s,+}}{\varepsilon _a \sigma _B} \Bigr) ^{1/4}, \quad\text{ and } \quad 
\Bigl( \frac{q\beta _{s,-}}{\sigma _B} \Bigr) ^{1/4} < T_s \leq 
\Bigl( \frac{2q\beta _{s,+}}{\varepsilon _a \sigma _B} \Bigr) ^{1/4} ,
\end{equation*}
that is, uniform bounds of the equilibria independent of $\lambda >0$.

If $\varepsilon _a \in (1,2)$, from \eqref{lambda-diff} and the left hand side of \eqref{lambda-posi} we deduce that
$$ \varepsilon _a \sigma _B T_a ^4
= \sigma _B (\varepsilon _a -1) T_s ^4 + q\beta _s (T_s)
> \sigma _B (\varepsilon _a -1) T_a ^4 + q\beta _s (T_s) .
$$
Hence it holds that
$$ \sigma _B T_a ^4 > q\beta _s (T_s) ,$$
which gives
\begin{equation}
\label{lambda-12-inf}
T_a > \Bigl( \frac{q\beta _{s,-}}{\sigma _B} \Bigr) ^{1/4}.\end{equation}
Using again \eqref{lambda-diff} and \eqref{lambda-posi}, we get
$$ \varepsilon _a \sigma _B T_a ^4
= \sigma _B (\varepsilon _a -1) T_s ^4 + q\beta _s (T_s)
< \sigma _B (\varepsilon _a -1) 2 T_a ^4 + q\beta _s (T_s) .$$
Thus we have
$$ (2 - \varepsilon _a) \sigma _B T_a ^4 < q\beta _s (T_s) $$
which implies that
\begin{equation}
\label{lambda-12-sup}
T_a \leq \Bigl( \frac{q\beta _{s,+}}{(2-\varepsilon _a )\sigma _B} \Bigr) ^{1/4}.
\end{equation}
From \eqref{lambda-posi}, \eqref{lambda-12-inf} and \eqref{lambda-12-sup} we obtain
\begin{equation*}
\Bigl( \frac{q\beta _{s,-}}{\sigma _B} \Bigr) ^{1/4} < T_a \leq \Bigl( \frac{q\beta _{s,+}}{(2-\varepsilon _a )\sigma _B} \Bigr) ^{1/4}, \quad \text{ and } \quad 
\Bigl( \frac{q\beta _{s,-}}{\sigma _B} \Bigr) ^{1/4} < T_s \leq \Bigl( \frac{2q\beta _{s,+}}{(2-\varepsilon _a )\sigma _B} \Bigr) ^{1/4}.
\end{equation*}
Therefore, also in the case $\eps_a\in(1,2)$ we have found uniform bounds, independent of $\lambda >0$, for the equilibrium points. This will be useful later (see in particular the proof of Corollary \ref{cor-greenhouse}).\qed


\subsection{Equilibrium points: existence} \hfill

As proved in Proposition \ref{prop-comp-lambda-eq-cv}, problem \eqref{2layer-pbm-EDO} has at least one equilibrium point which follows from the convergence of the solutions. We can also directly prove their existence as follows.
\begin{Lemma}
\label{lem-nbre-eq-1}
Given $\lambda >0$ and $\varepsilon _a \in (0,2)$, problem \eqref{2layer-pbm-EDO} has at least one equilibrium point.
\end{Lemma}

We are going to give two short proofs of this result, each one having its own interest.

\subsubsection{Geometrical proof of Lemma \ref{lem-nbre-eq-1}}\label{sec-geom-pf} \hfill

Consider 
\begin{equation*}
    \mathcal C_1 := \{(T_a,T_s) \in \overline{\mathcal Q}\,:\,
    -\lambda (T_a - T_s) + \varepsilon _a \sigma _B \vert T_s \vert ^3 T_s  - 2 \varepsilon _a \sigma _B \vert T_a \vert ^3 T_a =0 \}
\end{equation*}
and
\begin{equation*}
    \mathcal C_2 := \{ (T_a,T_s)\in \overline{\mathcal Q}\,:\,
    -\lambda (T_s - T_a) - \sigma _B \vert T_s \vert ^3 T_s  + \varepsilon _a \sigma _B \vert T_a \vert ^3 T_a + q \beta _s (T_s) = 0 \} .
\end{equation*}
We first analyse $\mathcal C_1$. Given $T_s \geq 0$, there exists a unique value $T_a$, that we denote $T_a ^{(1)} (T_s)$, such that 
$$ \lambda T_a + 2 \varepsilon _a \sigma _B \vert T_a \vert ^3 T_a = \lambda T_s + \varepsilon _a \sigma _B \vert T_s \vert ^3 T_s ,$$
that is, such that $(T_a,T_s) \in \mathcal C_1$. By the implicit function theorem we deduce that $\mathcal C_1$ is a curve. Moreover, it contains $(0,0)$ and 
$$ T_a ^{(1)} (T_s) \sim 2^{-1/4} T_s \quad \text{ as } T_s \to +\infty .$$
Let us now study $\mathcal C_2$. It contains points of the form  $(0,T_s)$. Let $T_{s,\text{max}}$ be the largest value $T_s$ such that $(0, T_s)\in \mathcal C_2$, that is, the largest value such that
$$ \lambda T_{s} + \sigma _B  T _{s} ^4 =  q \beta (T_{s}) .$$
Then 
$$ T_s > T_{s,\text{max}} \quad \implies \quad \lambda T_{s} + \sigma _B  T_s ^4 - q \beta _s (T_s) >  0 ,$$
and therefore, for all $T_s \geq T_{s,\text{max}}$, there exists a unique value of $T_a$, denoted by $T_a ^{(2)} (T_s)$, such that
$$ \lambda T_a + \varepsilon _a \sigma _B T_a ^4 = \lambda T_{s} + \sigma _B  T_s ^4 - q \beta _s (T_s) .$$
Note that 
$$ T_a ^{(2)} (T_s)  \sim \varepsilon _a ^{-1/4} T_s
\quad \text{ as } T_s \to +\infty .$$
Consider now the function
$$ T_s \in [T_{s,\text{max}}, +\infty ) \to \mathbb R, \quad T_s \mapsto T_a ^{(1)} (T_s) - T_a ^{(2)} (T_s) .$$
This map is continuous, positive for $T_s = T_{s,\text{max}}$ (since $T_a ^{(2)} (T_{s,\text{max}}) =0$) and negative for $T_s$ large. This implies that the sets 
$$ \mathcal C_{1,\text{max}} := \{(T_a ^{(1)} (T_s),T_s), T_s \geq T_{s,\text{max}} \}$$
and
$$
\mathcal C_{2,\text{max}} := \{ (T_a ^{(2)} (T_s),T_s), T_s \geq T_{s,\text{max}} \} $$
have to intersect at least once. This implies that problem \eqref{2layer-pbm-EDO} has at least an equilibrium point. \qed


\subsubsection{Analytical proof of Lemma \ref{lem-nbre-eq-1}} \label{subsec-analytic-proof}\hfill

Equilibrium points $(T_a, T_s)$ of problem \eqref{2layer-pbm-EDO} satisfy \eqref{2layer-pbm-EDO-eq-lambda}. Therefore, for such points it holds
\begin{equation*}
\begin{split}
q \beta _s (T_s) &= \Bigl( \lambda T_s + \sigma _B T_s ^4 \Bigr)  - \Bigl( \lambda T_a + \varepsilon _a \sigma _B T_a ^4 \Bigr)
\\
&= \Bigl( \lambda T_s + \sigma _B T_s ^4 \Bigr) - \Bigl( \frac{\lambda}{2} T_a + \frac{1}{2} \Bigr[ \lambda T_a + 2 \varepsilon _a \sigma _B T_a ^4 \Bigr] \Bigr)
\\
&= \Bigl( \lambda T_s + \sigma _B T_s ^4 \Bigr) - \Bigl( \frac{\lambda}{2} T_a + \frac{1}{2} \Bigr[ \lambda T_s + \varepsilon _a \sigma _B T_s ^4 \Bigr] \Bigr)
\\
&= \frac{\lambda}{2} (T_s - T_a ) + \sigma _B (1 - \frac{\varepsilon _a}{2}) T_s ^4 .
\end{split}
\end{equation*}
Thus, for nonnegative $T_a$ and $T_s$, $(T_a, T_s)$ solves \eqref{2layer-pbm-EDO-eq-lambda} if and only if it verifies the equivalent system
\begin{equation}
\label{2layer-pbm-EDO-eq-lambda-equiv}
\begin{cases}
\lambda T_s + \varepsilon _a \sigma _B  T_s ^4   & = \lambda T_a + 2 \varepsilon _a \sigma _B T_a ^4, \\
\frac{\lambda}{2} (T_s - T_a ) + \sigma _B (1 - \frac{\varepsilon _a}{2}) T_s ^4 &  = q \beta _s (T_s) .
\end{cases}
\end{equation}
Now we prove that \eqref{2layer-pbm-EDO-eq-lambda-equiv} has at least one solution. Indeed, given $T_s \geq 0$, the first equation in \eqref{2layer-pbm-EDO-eq-lambda-equiv} has a unique solution $T_a=T_a ^{(1)} (T_s)$. Furthermore, it holds that
$$ T_a ^{(1)} (0)=0 \quad \text{ and } \quad T_a ^{(1)} (T_s) \sim _{T_s \to 0} T_s,$$
$$ T_a ^{(1)} (T_s) \to _{T_s \to +\infty} +\infty \quad \text{ and } \quad T_a ^{(1)} (T_s) \sim _{T_s \to +\infty} 2^{-1/4} T_s .$$
Moreover, as already observed in \eqref{lambda-posi}, we have
\begin{equation}
\label{lambda-posi2}
\forall\, T_s >0, \quad T_a ^{(1)} (T_s) < T_s < 2^{1/4} T_a ^{(1)} (T_s) .
\end{equation}
Finally, since $T_s \mapsto T_a ^{(1)} (T_s)$ is a smooth function, we obtain that
$$ \frac{dT_a ^{(1)}}{dT_s}(T_s) = \frac{\lambda + 4 \varepsilon _a \sigma _B  T_s ^3}{\lambda + 8\varepsilon _a \sigma _B [ T_a ^{(1)} (T_s)] ^3} .$$
Therefore
\begin{equation*}
\begin{split}
 1 - \frac{dT_a ^{(1)}}{dT_s} (T_s) 
&= 1 - \frac{\lambda + 4 \varepsilon _a \sigma _B  T_s ^3}{\lambda + 8\varepsilon _a \sigma _B  [ T_a ^{(1)} (T_s)] ^3}
\\
&= \frac{8\varepsilon _a \sigma _B  [ T_a ^{(1)} (T_s)] ^3 - 4 \varepsilon _a \sigma _B  T_s ^3}{\lambda + 8\varepsilon _a \sigma _B  [ T_a ^{(1)} (T_s)] ^3} 
\\
&= 4 \varepsilon _a \sigma _B \frac{2 [ T_a ^{(1)} (T_s)] ^3 -  T_s ^3}{\lambda + 8\varepsilon _a \sigma _B  [ T_a ^{(1)} (T_s)] ^3}  .    
\end{split}
\end{equation*}
Thanks to \eqref{lambda-posi2}, we get that
\begin{equation*}
     2 [ T_a ^{(1)} (T_s)] ^3 -  T_s ^3 > 2 ^{3/4} [ T_a ^{(1)} (T_s)] ^3 -  T_s ^3 = [ 2 ^{1/4} T_a ^{(1)} (T_s)] ^3 -  T_s ^3 >0 ,
\end{equation*}
and so
$$  1 - \frac{dT_a ^{(1)}}{dT_s} (T_s)  >0 .$$
This implies that the function
\begin{equation}\label{def-Phi1*}
\Phi_1: T_s \mapsto \frac{\lambda}{2} (T_s - T_a ^{(1)} (T_s) )
\end{equation}
is strictly increasing on $[0,+\infty)$. Therefore the following function
\begin{equation}\label{def-Phi*}
    \Phi: T_s \mapsto \frac{\lambda}{2} (T_s - T_a ^{(1)} (T_s) ) + \sigma _B (1 - \frac{\varepsilon _a}{2}) T_s ^4 
\end{equation}
is also strictly increasing on $[0,+\infty)$. Moreover, $\Phi (0)=0$, and $\Phi (T_s) \sim \sigma _B (1 - \frac{\varepsilon _a}{2}) T_s ^4$ as $T_s \to \infty$. Hence, there exists at least one positive value of $T_s$, denoted by $T_{s,*}$, such that $\Phi (T_{s,*})=q \beta _s (T_{s,*})$. Therefore $(T_a ^{(1)} (T_{s,*}), T_{s,*})$ solves \eqref{2layer-pbm-EDO-eq-lambda-equiv} or, equivalently, it satisfies \eqref{2layer-pbm-EDO-eq-lambda}.
We have therefore analytically proved Lemma \ref{lem-nbre-eq-1}.


\subsection{Number of equilibria}\label{sec-number-eq-pts}\hfill

\subsubsection{Existence of at most one warm equilibrium}\hfill

By using the argument Section \ref{subsec-analytic-proof}, we deduce that there exists at most one value $T_s \in (0,T_{s,-}]$ 
and at most one value $T_s \geq T_{s,+}$ such that $\Phi (T_s)= q \beta _s (T_s)$, that is, at most one "cold" and one "warm" equilibrium.

Indeed, the existence two different warm equilibria would imply the existence of two different values $T_s , \tilde T_s \geq T_{s,+}$ such that
$$ \Phi (T_s) = q \beta _s (T_s) , \quad \text{ and } \quad \Phi (\tilde T_s) = q \beta _s (\tilde T_s) .$$
However, since $T_s, \tilde T_s \geq T_{s,+}$ then $q \beta _s (T_s) = q \beta _s (\tilde T_s)$. Therefore, it would hold that $\Phi (T_s) = \Phi (\tilde T_s)$, which implies that $T_s= \tilde T_s$ (because we recall that $\Phi$ is strictly increasing). This proves the third point of Proposition \ref{prop-comp-lambda-eq-gene}. \qed

\subsubsection{Existence of at most a finite number of equilibria}\hfill

Now, let us prove that the number of equilibrium is finite  (for the case of a piecewise linear function such as $\beta _s$).
\begin{Lemma}
\label{lem-nbre-eq-1*}
Given $\lambda >0$ and $\varepsilon _a \in (0,2)$, problem \eqref{2layer-pbm-EDO} has
at most a finite number of equilibrium points.
\end{Lemma}

\noindent {\it Proof of Lemma \ref{lem-nbre-eq-1*}.} We have already proved that there are at most one cold and one warm equilibrium point. If the number of equilibria is infinite, there exists a sequence $(T_{s,n})_n$ of distinct values belonging to $[T_{s,-}, T_{s,+}]$ satisfying 
$$ \Phi (T_{s,n}) = q \beta _s  (T_{s,n}).$$
Such sequence admits a converging subsequence. Moreover, it is possible to extract a strictly monotone subsequence. Indeed, consider a convergent subsequence
$(T_{s,\varphi (n)})_n$ and let $T_{s,\infty}$ be the limit value. Define the set
$$ \mathcal A := \{n\in\mathbb{N}\,:\, T_{s,\varphi (n)} < T_{s,\infty} \}.$$
It is clear that $\mathcal A \cup (\mathbb N \setminus \mathcal A) = \mathbb N$, therefore either $\mathcal A$, $\mathbb N \setminus \mathcal A$, or both sets have a infinite number of elements. We can then construct a strictly monotone subsequence, $(T_{s,\psi (n)})_n$, converging to $T_{s,\infty}$. Since it holds that
$$ (\Phi - q \beta _s) (T_{s,\psi (n)}) = \Phi (T_{s,\psi (n)}) - q \beta _s (T_{s,\psi (n)}) = 0, $$
by applying Rolle's Theorem to the function $\Phi - q \beta _s$, we obtain that if $(T_{s,\psi (n)})_n$
is increasing there exists 
$\tilde T_{s,n} \in (T_{s,\psi (n)}, T_{s,\psi (n+1)})$  (or $\tilde T_{s,n} \in (T_{s,\psi (n+1)}, T_{s,\psi (n)})$ if $(T_{s,\psi (n)})_n$
is decreasing), such that
$$ \Phi ' (\tilde T_{s,n}) - q \beta ' _s  (\tilde T_{s,n}) = 0,$$
or, equivalently,
$$ \Phi ' (\tilde T_{s,n}) = q \beta ' _s  (\tilde T_{s,n})  = q \frac{\beta _{s,+} - \beta _{s,-}}{T _{s,+} - T _{s,-}} .$$
Since $\Phi$ is analytic on $(0,+\infty)$ \cite[Proposition 2.20 p. 39]{Conways} also $\Phi '$ is analytic. However, if $\Phi '$ is equal to the same value an infinite number of times on a compact set, then it must be constant. This means that $\Phi$ must grow at most linearly at infinity, which turns out to be false. Therefore, we conclude that problem \eqref{2layer-pbm-EDO} has
at most a finite number of equilibrium points. \qed

\subsubsection{Convexity of function $\Phi$ and number of equilibrium points} \hfill

Once we have established that the number of equilibrium points is finite, it is natural to try to have an estimate of
such number. We recall that equilibrium points are of the form $(T_a = T_a ^{(1)} (T_s), T_s)$ where $T_s$ satisfies
\begin{equation}
\label{eq-eq-Phi}
\Phi (T_s) = q \beta _s (T_s), 
\end{equation}
and $\Phi$ is defined in \eqref{def-Phi*}.
In the case $\lambda =0$, function $\Phi$ is convex, which implies that there exist at most 3 equilibrium points. When $\lambda >0$ we have obtained the following result concerning the physically relevant case $\varepsilon _a \leq1$.

\begin{Lemma}
\label{lem-nbre-eq-1*conv-phys}

Assume that $\varepsilon _a \leq1$ and let $\lambda >0$. Then $\Phi $, defined in \eqref{def-Phi*}, is strictly increasing and strictly convex on $[0,+\infty )$.
Therefore, problem \eqref{2layer-pbm-EDO} has at most three equilibrium points.
\end{Lemma}
For what concerns the case $\varepsilon _a \in (1,2)$, we have obtained some information by using numerical tests. We describe such results in the remark that follows.
\begin{Remark}
\label{lem-nbre-eq-1*conv-simpl}

Let $\lambda >0$. Then, we observe thanks to some numerical tests that there exists a unique universal constant $\varepsilon _{a,0}\in (1.99, 1.991)$, independent of $\lambda >0$, such that
\begin{itemize}
\item if $\varepsilon _a \in (0,\varepsilon _{a,0})$, $\Phi $ is strictly convex on $[0,+\infty )$;
\item if $\varepsilon _a \in (\varepsilon _{a,0},2)$, $\Phi $ is successively convex, concave and convex on $[0,+\infty )$. 
\end{itemize}
From the above numerical results we deduce that
\begin{itemize}
\item if $\varepsilon _a \in (0,\varepsilon _{a,0})$, problem \eqref{2layer-pbm-EDO} has at most three equilibrium points;
\item if $\varepsilon _a \in (\varepsilon _{a,0},2)$, problem \eqref{2layer-pbm-EDO} has at most five equilibrium points. 
\end{itemize}
\end{Remark}
The proof of Lemma \ref{lem-nbre-eq-1*conv-phys}can be found in Appendix \ref{sec-conv-Phi1}. It is based on a careful analysis of $\Phi ''$.
We also explain in Appendix \ref{sec-conv-Phi1} the results stated in Remark \ref{lem-nbre-eq-1*conv-simpl}, which take advantage of some computations from Lemma \ref{lem-nbre-eq-1*conv-phys} and from some numerical tests.

\subsection{ Local stability of the equilibrium points}\hfill
\label{sec-stab-lambda}

\subsubsection{ Local stability of a warm (resp. cold) equilibrium point}\hfill

Assume that there exists a warm equilibrium: $(T_a,T_s)$ satisfying \eqref{2layer-pbm-EDO-eq-lambda-equiv}
with $T_s > T_{s,+}$. Then, we have
\begin{equation}
\label{DF-lambda} 
D F\begin{pmatrix} T_a \\ T_s \end{pmatrix}  
= 
\begin{pmatrix}
\frac{1}{\gamma _a}[-\lambda -8 \varepsilon _a \sigma _B T_a ^3 ]
& 
\frac{1}{\gamma _a} [\lambda + 4 \varepsilon _a \sigma _B T_s ^3] 
\\
\frac{1}{\gamma _s}[\lambda + 4 \varepsilon _a \sigma _B T_a ^3 ]
& 
\frac{1}{\gamma _s}[-\lambda -4 \sigma _B T_s ^3 ]
\end{pmatrix} ,
\end{equation}
where $F$ is defined in \eqref{def-Fnonlin}.
It is easy to check that the trace of the above matrix is negative. Moreover, from the characteristic polynomial associated to the matrix, we deduce that it has two different real eigenvalues. Moreover, a direct computation of the determinant of $DF\begin{pmatrix}T_a \\T_s\end{pmatrix}$ matrix gives
\begin{equation*}
\begin{split}
\gamma _a \gamma _s \, \det D F \begin{pmatrix} T_a \\ T_s \end{pmatrix}&= [\lambda +8 \varepsilon _a \sigma _B T_a ^3 ][\lambda +4 \sigma _B T_s ^3 ]- [\lambda + 4 \varepsilon _a \sigma _B T_s ^3][\lambda + 4 \varepsilon _a \sigma _B T_a ^3 ]
\\
&= 4 \lambda \sigma _B \Bigl[ (1-\varepsilon _a) T_s ^3 + \varepsilon _a T_a ^3 \Bigr] + 16 \varepsilon _a (2 - \varepsilon _a) \sigma _B ^2 T_a ^3 T_s ^3 .
\end{split}
\end{equation*}
The bound in \eqref{lambda-posi} yields
\begin{equation*}
    \Bigl[ (1-\varepsilon _a) T_s ^3 + \varepsilon _a T_a ^3 \Bigr]
\geq \Bigl[ (1-\varepsilon _a) T_s ^3 + \varepsilon _a 2^{-3/4} T_s ^3 \Bigr]= \Bigl( 1-\varepsilon _a + 2^{-3/4} \varepsilon _a \Bigr) T_s ^3 
.
\end{equation*}
Since $\varepsilon _a \in (0,2)$ and $2^{-3/4} -1 <0$, we obtain
\begin{equation*}
    \Bigl( 1 + (2^{-3/4} -1) \varepsilon _a \Bigr) T_s ^3
\geq \Bigl( 1 + 2 (2^{-3/4} -1) \Bigr) T_s ^3
= (2^{1/4} -1 ) T_s ^3 >0 ,
\end{equation*}
which implies that
\begin{equation*}
    \gamma _a \gamma _s \, \det D F \begin{pmatrix} T_a \\ T_s \end{pmatrix}
\geq 4 \lambda \sigma _B \, (2^{1/4} -1 ) T_s ^3 + 16 \varepsilon _a (2 - \varepsilon _a) \sigma _B ^2 T_a ^3 T_s ^3 >0 .
\end{equation*}
Therefore, the two eigenvalues of $D F \begin{pmatrix} T_a \\ T_s \end{pmatrix} $ are negative and the warm equilibrium is asymptotically exponentially stable. The same holds true for a cold equilibrium ($T_s < T_{s,-}$).

\subsubsection{ Local stability of an intermediate equilibrium point}\hfill
\label{sec-unstable-interm}

Assume that there exists an intermediate equilibrium point $(\tilde T_a, \tilde T_s)$ satisfying \eqref{2layer-pbm-EDO-eq-lambda-equiv}
with $\tilde T_s \in ( T_{s,-},  T_{s,+})$. We compute the Jacobian matrix of $F$ at t $(\tilde T_a, \tilde T_s)$:
\begin{equation}
\label{DF-lambda-q} 
D F\begin{pmatrix} \tilde T_a \\ \tilde T_s \end{pmatrix}  
= 
\begin{pmatrix}
\frac{1}{\gamma _a}[-\lambda -8 \varepsilon _a \sigma _B \tilde T_a ^3 ]
& 
\frac{1}{\gamma _a} [\lambda + 4 \varepsilon _a \sigma _B \tilde T_s ^3] 
\\
\frac{1}{\gamma _s}[\lambda + 4 \varepsilon _a \sigma _B \tilde T_a ^3 ]
& 
\frac{1}{\gamma _s}[-\lambda -4 \sigma _B \tilde T_s ^3 + q \beta ' _s (\tilde T_s) ]
\end{pmatrix} ,
\end{equation}
where $F$ is defined in \eqref{def-Fnonlin}. Thus, we get
\begin{multline}
\label{eq-stab-interm}
\gamma _a \gamma _s \, \det D F \begin{pmatrix} \tilde T_a \\ \tilde T_s \end{pmatrix}=[\lambda +8 \varepsilon _a \sigma _B \tilde T_a ^3 ]
\Bigl[ \lambda +4 \sigma _B \tilde T_s ^3 - q \beta ' _s (\tilde T_s)
\\- \frac{\lambda + 4 \varepsilon _a \sigma _B \tilde T_s ^3}{\lambda +8 \varepsilon _a \sigma _B \tilde T_a ^3 } \Bigl( \lambda + 4 \varepsilon _a \sigma _B \tilde T_a ^3 \Bigr)  \Bigr]  .
\end{multline}
Let us compute the derivative of the function $\Phi$
defined in \ref{def-Phi*} at $\tilde T_s$:
\begin{multline*} \Phi '(\tilde T_s) = \frac{\lambda}{2} (1- \frac{dT_a ^{(1)}}{dT_s}(\tilde T_s)) + 4 \sigma _B (1 - \frac{\varepsilon _a}{2}) \tilde T_s ^3
\\
= \frac{\lambda}{2} (1- \frac{\lambda + 4 \varepsilon _a \sigma _B  \tilde T_s ^3}{\lambda + 8\varepsilon _a \sigma _B [ T_a ^{(1)} (\tilde T_s)] ^3}) + 4 \sigma _B (1 - \frac{\varepsilon _a}{2}) \tilde T_s ^3 ,
\end{multline*}
from which we obtain that
\begin{equation}
\label{eq-interm}
\Phi '(\tilde T_s)
= \lambda +4 \sigma _B \tilde T_s ^3 
- \frac{\lambda + 4 \varepsilon _a \sigma _B \tilde T_s ^3}{\lambda +8 \varepsilon _a \sigma _B \tilde T_a ^3 } \Bigl( \lambda + 4 \varepsilon _a \sigma _B \tilde T_a ^3 \Bigr) .
\end{equation}
Indeed, explicit computations show that
\begin{equation*}
\begin{split}
\Phi '(\tilde T_s)&- \Bigl[ \lambda +4 \sigma _B \tilde T_s ^3 
- \frac{\lambda + 4 \varepsilon _a \sigma _B \tilde T_s ^3}{\lambda +8 \varepsilon _a \sigma _B \tilde T_a ^3 } \Bigl( \lambda + 4 \varepsilon _a \sigma _B \tilde T_a ^3 \Bigr) \Bigr]
\\
\\
&= \Bigl[ \frac{\lambda}{2} (1- \frac{\lambda + 4 \varepsilon _a \sigma _B  \tilde T_s ^3}{\lambda + 8\varepsilon _a \sigma _B \tilde T_a ^3}) - 2 \varepsilon _a \sigma _B \tilde T_s ^3 \Bigr]- \Bigl[ \lambda 
- \frac{\lambda + 4 \varepsilon _a \sigma _B \tilde T_s ^3}{\lambda +8 \varepsilon _a \sigma _B \tilde T_a ^3 } \Bigl( \lambda + 4 \varepsilon _a \sigma _B \tilde T_a ^3 \Bigr) \Bigr]
\\
&= -\frac{\lambda}{2} - 2 \varepsilon _a \sigma _B \tilde T_s ^3 + \frac{\lambda + 4 \varepsilon _a \sigma _B \tilde T_s ^3}{\lambda +8 \varepsilon _a \sigma _B \tilde T_a ^3 } \Bigl( -\frac{\lambda}{2} + \lambda + 4 \varepsilon _a \sigma _B \tilde T_a ^3 \Bigr)
\\
&= -\frac{\lambda}{2} - 2 \varepsilon _a \sigma _B \tilde T_s ^3 + \frac{1}{2} \Bigl[ \lambda + 4 \varepsilon _a \sigma _B \tilde T_s ^3 \Bigr] = 0 .
\end{split}
\end{equation*}
By using \eqref{eq-interm} in \eqref{eq-stab-interm} we obtain
\begin{equation}
\label{eq-interm-stab2}
\gamma _a \gamma _s \, \det D F \begin{pmatrix} \tilde T_a \\ \tilde T_s \end{pmatrix} =
[\lambda +8 \varepsilon _a \sigma _B \tilde T_a ^3 ] 
\Bigl[ \Phi '(\tilde T_s) - q \beta ' _s (\tilde T_s) \Bigr] .
\end{equation}
Observe that the above formula is the generalization of that obtained in section \ref{sec-local-stab} (for the case $\lambda =0$).

Now if $\varepsilon _a <1$, then there exist at most three equilibrium points (see Lemma \ref{lem-nbre-eq-1*conv-phys}).
Assume that there are exactly three equilibrium points: one cold, one intermediate $(\tilde T_a, \tilde T_s)$ and one warm. Function $\Phi$ is strictly convex, and we have
$$ \Phi '(\tilde T_s) - q \beta ' _s (\tilde T_s) < 0 ,$$
(see Section \ref{sec-local-stab}), otherwise the convexity of $\Phi$ would prevent the existence of the warm equilibrium. Therefore, we deduce that the intermediate equilibrium point is unstable, as in the case $\lambda =0$.


\subsection{Phase space analysis for $\lambda >0$} \hfill
\label{sec-asympt-lambda>0}

In this section we extend to the case $\lambda >0$ the study of the phase space we have already done for the case $\lambda =0$. Consider again $\mathcal C_1$ and $\mathcal C_2$ defined in section \ref{sec-geom-pf}. We observe that $\mathcal C_2$ has at most three intersection points with the vertical axis $(0,T_s)$. Indeed, the solutions of equation $\lambda T_s + \sigma _B T_s ^4 = q\beta _s (T_s)$ are at most three because $T_s \mapsto \lambda T_s + \sigma _B T_s ^4$ is strictly convex and increasing. We have already observed in Section \ref{sec-number-eq-pts} that equilibrium points are related to the solutions of the equation $\Phi (T_s)=q \beta _s (T_s)$, where $\Phi$ is defined in \eqref{def-Phi*}.
Assume that $\Phi$ is strictly convex  ( which indeed happens at least for $\varepsilon _a < 1$, see Lemma \ref{lem-nbre-eq-1*conv-phys}), and that there exist three equilibrium points $T_{s,1} < T_{s,2} < T_{s,3}$ (that is, we are in the same situation of Figure \ref{fig2:region}, case (A)). Then 
$q\beta _s (T_s) - \Phi (T_s)$ is positive on $(0,T_{s,1})$, negative on $(T_{s,1}, T_{s,2})$,
positive on $(T_{s,2}, T_{s,3})$ and negative on $(T_{s,3}, +\infty)$. We use the following identity
$$ \lambda T_a ^{(1)} (T_s) + 2 \varepsilon _a \sigma _B T_a ^{(1)} (T_s) ^4 = \lambda T_s + \varepsilon _a \sigma _B T_s ^4 ,$$
inside the expression of $\Phi(T_s)$ and we get
\begin{equation*}
\begin{split}
\Phi (T_s) &= \frac{\lambda}{2} (T_s - T_a ^{(1)} (T_s) ) + \sigma _B (1 - \frac{\varepsilon _a}{2}) T_s ^4
\\
&= \frac{\lambda}{2} T_s + \sigma _B T_s ^4 - \frac{\varepsilon _a \sigma _B}{2} T_s ^4 - \frac{\lambda}{2} T_a ^{(1)} (T_s)
\\
&= \Bigl( \lambda T_s + \sigma _B T_s ^4 \Bigr) - \frac{\lambda}{2} T_s  - \frac{\varepsilon _a \sigma _B}{2} T_s ^4 - \frac{\lambda}{2} T_a ^{(1)} (T_s)
\\
&= \Bigl( \lambda T_s + \sigma _B T_s ^4 \Bigr) - \frac{1}{2} \Bigl( \lambda T_s + \varepsilon _a \sigma _B T_s ^4 \Bigr) - \frac{\lambda}{2} T_a ^{(1)} (T_s)
\\
&= \Bigl( \lambda T_s + \sigma _B T_s ^4 \Bigr) - \frac{1}{2} \Bigl( \lambda T_a ^{(1)} (T_s) + 2 \varepsilon _a \sigma _B T_a ^{(1)} (T_s) ^4 \Bigr) \\
&\qquad\qquad\qquad\qquad\qquad\qquad\qquad\qquad- \frac{\lambda}{2} T_a ^{(1)} (T_s)
\\
&= \Bigl( \lambda T_s + \sigma _B T_s ^4 \Bigr) - \Bigl( \lambda T_a ^{(1)} (T_s) + \varepsilon _a \sigma _B T_a ^{(1)} (T_s) ^4 \Bigr) .
\end{split}
\end{equation*}
We note that the function
$$\psi ^{(\lambda)}: \mathbb R \to \mathbb R, \quad T_a \mapsto \lambda T_a + \varepsilon _a \sigma _B \vert T_a \vert ^3 T_a $$
is strictly increasing. Therefore, given $T_s$ there exists one and only one value $T_a ^{(2)} (T_s)$ such that
$$\psi ^{(\lambda)} (T_a ^{(2)} (T_s))= \lambda T_s + \sigma _B \vert T_s \vert ^3 T_s - q \beta _s (T_s) .$$
Thus, we obtain that
\begin{equation*}
    \Phi (T_s)
= \Bigl( \psi ^{(\lambda)} (T_a ^{(2)} (T_s)) + q \beta _s (T_s) \Bigr)- \Bigl( \lambda T_a ^{(1)} (T_s) + \varepsilon _a \sigma _B T_a ^{(1)} (T_s) ^4 \Bigr) ,
\end{equation*}
which implies
$$ \Phi (T_s) - q \beta _s (T_s) = \psi ^{(\lambda)} (T_a ^{(2)} (T_s)) - \psi ^{(\lambda)} (T_a ^{(1)} (T_s)).$$
Therefore, $T_a ^{(1)} (T_s) - T_a ^{(2)} (T_s)$ and $q \beta _s (T_s) - \Phi (T_s)$ have the same sign, that is,
$T_a ^{(1)} (T_s) - T_a ^{(2)} (T_s)$ is positive on $(0,T_{s,1})$, negative on $(T_{s,1}, T_{s,2})$,
positive on $(T_{s,2}, T_{s,3})$ and negative on $(T_{s,3}, +\infty)$. This describes the relative position of $\mathcal C_1$
with respect to $\mathcal C_2$ which is similar to that represented in Figure \ref{fig8:region} (the only difference is that, for the current case $\lambda>0$, $\mathcal C_1$ is no more a line, but the curve of the strictly increasing function $T_a \mapsto T_s ^{(1)} (T_a)$). The direction of vector field $F$ on the curves $\mathcal C_1$ and $\mathcal C_2$ is similar to the described by arrows in Figure \ref{fig8:region}, therefore we will have same phase plane analysis as in the case $\lambda=0$. Note that the existence of the separatrix between the basins of attraction of the stable equilibrium points
follows exactly in the same way from the property that the intermediate equilibrium point is unstable (which is proved in section \ref{sec-unstable-interm}), as in the case $\lambda=0$.
This concludes the proof of Proposition \ref{prop-comp-lambda-eq-gene}. \qed


\section{Sensitivity of the equilibria to parameters}
\label{sec8}

\subsection{Proof of Proposition \ref{prop-comp-lambda>0}: monotonicity of the equilibrium points with respect to $\lambda$} \hfill

Consider the function $G: \mathbb R \times \mathbb R^2 \to \mathbb R^2$ defined by
\begin{equation}
\label{def-G-lambda}
G (\lambda, \begin{pmatrix} T_a \\ T_s \end{pmatrix} )
= \begin{pmatrix} \frac{1}{\gamma _a} \Bigl[ -\lambda (T_a - T_s) + \varepsilon _a \sigma _B \vert T_s \vert ^3 T_s  - 2 \varepsilon _a \sigma _B \vert T_a \vert ^3 T_a \Bigr] \\ \frac{1}{\gamma _s} \Bigl[ -\lambda (T_s - T_a) - \sigma _B \vert T_s \vert ^3 T_s  + \varepsilon _a \sigma _B \vert T_a \vert ^3 T_a + \mathcal R_s \Bigr] \end{pmatrix} .
\end{equation}
If $(T_a,T_s)$ is an equilibrium point of \eqref{2layer-pbm-EDO} with parameter $\lambda$, that is, $(T_a,T_s)$ solves \eqref{2layer-pbm-EDO-eq-lambda}, then $G (\lambda,  \begin{pmatrix} T_a \\ T_s \end{pmatrix} )=0$.

Fixed $\lambda >0$ and a point $(T_a,T_s)\in \mathcal{Q}$ of differentiability for $\beta _s$, we differentiate with respect to the second variable of $G$:
\begin{equation}
\label{D2G} 
D_2 G(\lambda, \begin{pmatrix} T_a \\ T_s \end{pmatrix} )
= 
\begin{pmatrix}
\frac{1}{\gamma _a}[-\lambda -8 \varepsilon _a \sigma _B T_a ^3 ]
& 
\frac{1}{\gamma _a} [\lambda + 4 \varepsilon _a \sigma _B T_s ^3] 
\\
\frac{1}{\gamma _s}[\lambda + 4 \varepsilon _a \sigma _B T_a ^3 ]
& 
\frac{1}{\gamma _s}[-\lambda -4 \sigma _B T_s ^3 + q \beta ' _s (T_s) ]
\end{pmatrix} .
\end{equation}
The stability of an equilibrium point $(T_a,T_s)$ of \eqref{2layer-pbm-EDO} is related to the sign of determinant and of the trace of the matrix $D_2 G(\lambda, \begin{pmatrix} T_a \\ T_s \end{pmatrix} )$.

Now, assume that $(T_a ^{eq,\lambda ^*}, T_s ^{eq,\lambda ^*})$ is an equilibrium for problem \eqref{2layer-pbm-EDO} with parameter $\lambda = \lambda ^*$ for which it holds that
$$ \begin{cases}
\det D_2 G(\lambda ^*, \begin{pmatrix} T_a ^{eq,\lambda ^*}\\ T_s ^{eq,\lambda ^*} \end{pmatrix} ) >0 ,
\\ 
\text{Tr } D_2 G(\lambda ^*, \begin{pmatrix} T_a ^{eq,\lambda ^*}
\\ T_s ^{eq,\lambda ^*} \end{pmatrix} ) <0 ,
\\
T_s ^{eq,\lambda ^*} \notin [T_{s,-}, T_{s,+}] .
 \end{cases} $$
Then, by applying the implicit function theorem to $G$, we deduce that there exists a neighborhood $\mathcal V^*$ of $\lambda ^*$ in $\mathbb R$, a neighborhood $\mathcal V ^{eq,\lambda ^*}$ of $(T_a ^{eq,\lambda ^*}, T_s ^{eq,\lambda ^*})$ in $\mathbb R^2$, and a $C^1$ function
$$ \Psi ^{(\varepsilon _a)} : \mathcal V^* \to \mathcal V ^{eq,\lambda ^*}, \quad \lambda \mapsto \Psi ^{(\varepsilon _a)} (\lambda) $$
such that
$$ \begin{cases}
G (\lambda, \begin{pmatrix} T_a \\ T_s \end{pmatrix} )= \begin{pmatrix} 0 \\ 0 \end{pmatrix},
\\
\begin{pmatrix} T_a \\ T_s \end{pmatrix} \in \mathcal V ^{eq,\lambda ^*} \end{cases}
\,\, \iff \,\,\begin{cases} \begin{pmatrix} T_a \\ T_s \end{pmatrix} = \Psi ^{(\varepsilon _a)} (\lambda) , \\ \lambda \in \mathcal V^* . \end{cases}$$
This fact ensures the existence and uniqueness of an equilibrium of problem \eqref{2layer-pbm-EDO} for $\lambda$ close to $\lambda ^*$. Such equilibrium, that we will denote by $(T_a ^{eq,\lambda}, T_s ^{eq,\lambda})$, will be close to $(T_a ^{eq,\lambda ^*}, T_s ^{eq,\lambda ^*})$.
Furthermore, for $\lambda$ close to $\lambda ^*$, we still have
$T_s ^{eq,\lambda} \notin [T_{s,-}, T_{s,+}]$
and thus $(T_a ^{eq,\lambda}, T_s ^{eq,\lambda})$ is asymptotically stable (by continuity of $\det$ and $\text{Tr}$). The analyticity of the functions $\lambda \mapsto T_s ^{eq,\lambda}$
and $\lambda \mapsto T_a ^{eq,\lambda}$ can be deduced by the analytic version of the implicit function theorem, see e.g. \cite[Proposition 6.1, p. 138]{Cartan}.

It is interesting to determine the monotonicity 
of $\lambda \mapsto T_a ^{eq,\lambda}$ and $\lambda \mapsto T_s ^{eq,\lambda}$. We differentiate with respect to $\lambda$
the equation
$$ 0 = G(\lambda, \Psi ^{(\varepsilon _a)} (\lambda) ) = G (\lambda, \begin{pmatrix} T_a ^{eq,\lambda} \\ T_s ^{eq,\lambda} \end{pmatrix} ),$$
we get
$$ 0 = D_1 G (\lambda, \begin{pmatrix} T_a ^{eq,\lambda} \\ T_s ^{eq,\lambda} \end{pmatrix} )
+ D_2 G (\lambda, \begin{pmatrix} T_a ^{eq,\lambda} \\ T_s ^{eq,\lambda} \end{pmatrix} )  \cdot\begin{pmatrix} \frac{\partial T_a ^{eq,\lambda}}{\partial \lambda} \\ \frac{\partial T_s ^{eq,\lambda}}{\partial \lambda} \end{pmatrix} .$$
Therefore, we deduce that
$$ \begin{pmatrix} \frac{\partial T_a ^{eq,\lambda}}{\partial \lambda} \\ \frac{\partial T_s ^{eq,\lambda}}{\partial \lambda} \end{pmatrix}
= - D_2 G (\lambda, \begin{pmatrix} T_a ^{eq,\lambda} \\ T_s ^{eq,\lambda} \end{pmatrix} ) ^{-1} \, D_1 G (\lambda, \begin{pmatrix} T_a ^{eq,\lambda} \\ T_s ^{eq,\lambda} \end{pmatrix} ) .$$
Since we are interested in warm and cold equilibria, we have that, for such points, $\beta _s ' (T_s ^{eq,\lambda}) = 0$. We compute the inverse of $D_2G(\lambda, \begin{pmatrix} T_a ^{eq,\lambda} \\ T_s ^{eq,\lambda} \end{pmatrix} )$
\begin{multline*}
D_2 G (\lambda, \begin{pmatrix} T_a ^{eq,\lambda} \\ T_s ^{eq,\lambda} \end{pmatrix} ) ^{-1}=\frac{1}{\det D_2 G (\lambda, \begin{pmatrix} T_a ^{eq,\lambda} \\ T_s ^{eq,\lambda} \end{pmatrix} )}\cdot
\\ 
\begin{pmatrix}
\frac{1}{\gamma _s}[-\lambda -4 \sigma _B \vert T_s ^{eq,\lambda} \vert ^3 ]
& 
\frac{-1}{\gamma _a} [\lambda + 4 \varepsilon _a \sigma _B \vert T_s ^{eq,\lambda} \vert ^3] 
\\
\frac{-1}{\gamma _s}[\lambda + 4 \varepsilon _a \sigma _B \vert T_a ^{eq,\lambda} \vert ^3 ]
& 
\frac{1}{\gamma _a}[-\lambda -8 \varepsilon _a \sigma _B \vert T_a ^{eq,\lambda} \vert ^3 ]
\end{pmatrix} .
\end{multline*}
Since 
\begin{equation*}
     D_1 G (\lambda, \begin{pmatrix} T_a ^{eq,\lambda} \\ T_s ^{eq,\lambda} \end{pmatrix} ) =  
\begin{pmatrix} \frac{-1}{\gamma _a} (T_a ^{eq,\lambda} - T_s ^{eq,\lambda})\\ \frac{-1}{\gamma _s} (T_s ^{eq,\lambda} - T_a ^{eq,\lambda})  \end{pmatrix}
\\= (T_s ^{eq,\lambda} - T_a ^{eq,\lambda}) \begin{pmatrix} \frac{1}{\gamma _a}  -\frac{1}{\gamma _s}  \end{pmatrix} ,
\end{equation*}
we obtain that
\begin{equation*}
    \begin{pmatrix} \frac{\partial T_a ^{eq,\lambda}}{\partial \lambda} \\ \frac{\partial T_s ^{eq,\lambda}}{\partial \lambda} \end{pmatrix}
= \frac{T_s ^{eq,\lambda} - T_a ^{eq,\lambda}}{\gamma _a \, \gamma _s \, \det D_2 G (\lambda, \begin{pmatrix} T_a ^{eq,\lambda} \\ T_s ^{eq,\lambda} \end{pmatrix})}
\begin{pmatrix} 4 \sigma _B (1 - \varepsilon _a) \vert T_s ^{eq,\lambda} \vert ^3 \\
-4 \varepsilon _a \sigma _B \vert T_a ^{eq,\lambda} \vert ^3 \end{pmatrix} .
\end{equation*}
We observe that we have already proved that $T_s ^{eq,\lambda} - T_a ^{eq,\lambda} >0$ and therefore we conclude that $T_s ^{eq,\lambda}$ is decreasing with respect to $\lambda$ because the second component of the above vector is negative. Moreover, from the first component we deduce that $\frac{\partial T_a ^{eq,\lambda}}{\partial \lambda}$ has the same sign as $1-\varepsilon _a$. In particular, in the physical case, that is, for $\varepsilon _a <1$, $T_a ^{eq,\lambda}$ is increasing with respect to $\lambda$.
This concludes the proof of Proposition \ref{prop-comp-lambda>0}. \qed

\subsection{Proof of Proposition \ref{prop-comp-lambda=0}: monotonicity of the equilibrium points with respect to $\varepsilon _a$} \hfill

Given $\lambda >0$, consider the function $H: \mathbb R \times \mathbb R^2 \to \mathbb R^2$ defined by
\begin{equation}
\label{def-H-lambda}
H (\varepsilon _a, \begin{pmatrix} T_a \\ T_s \end{pmatrix} )
= \begin{pmatrix} \frac{1}{\gamma _a} \Bigl[ -\lambda (T_a - T_s) + \varepsilon _a \sigma _B \vert T_s \vert ^3 T_s  - 2 \varepsilon _a \sigma _B \vert T_a \vert ^3 T_a \Bigr] \\ \frac{1}{\gamma _s} \Bigl[ -\lambda (T_s - T_a) - \sigma _B \vert T_s \vert ^3 T_s  + \varepsilon _a \sigma _B \vert T_a \vert ^3 T_a + \mathcal R_s \Bigr] \end{pmatrix} .
\end{equation}
We note that $H (\varepsilon _a, \left( \begin{array}{c} T_a \\ T_s \end{array} \right) )=0$ if $(T_a,T_s)$ solves \eqref{2layer-pbm-EDO-eq-lambda}, that is, if $(T_a,T_s)$ is an equilibrium point of problem \eqref{2layer-pbm-EDO} with parameter $\lambda$.

Assume that $(T_a ^{eq,\varepsilon _a ^*}, T_s ^{eq,\varepsilon _a ^*})$ is an equilibrium point for problem \eqref{2layer-pbm-EDO} with parameter $\varepsilon _a = \varepsilon _a ^*$ such that
$$ \begin{cases}
\det D_2 H(\varepsilon _a ^*, \begin{pmatrix} T_a ^{eq,\varepsilon _a ^*}\\ T_s ^{eq,\varepsilon _a ^*} \end{pmatrix} ) >0 ,
\\ 
\text{Tr } D_2 H(\varepsilon _a ^*, \begin{pmatrix} T_a ^{eq,\varepsilon _a ^*}\\ T_s ^{eq,\varepsilon _a ^*} \end{pmatrix} ) <0 ,
\\
T_s ^{eq,\varepsilon _a ^*} \notin [T_{s,-}, T_{s,+}] .
 \end{cases} $$
Then, the implicit function theorem applied to $H$ implies that there exists a neighborhood $\mathcal V^*$ of $\varepsilon _a ^*$ in $\mathbb R$, a neighborhood $\mathcal V ^{eq,\varepsilon _a ^*}$ of $(T_a ^{eq,\lambda ^*}, T_s ^{eq,\lambda ^*})$ in $\mathbb R^2$, and a $C^1$ function
$$ \tilde \Psi ^{(\lambda)} : \mathcal V^* \to \mathcal V ^{eq,\varepsilon _a ^*}, \quad \varepsilon _a \mapsto \tilde \Psi ^{(\lambda)} (\varepsilon _a) $$
such that
$$ \begin{cases}
H (\varepsilon _a, \begin{pmatrix} T_a \\ T_s \end{pmatrix} )= \begin{pmatrix} 0 \\ 0 \end{pmatrix},
\\
\begin{pmatrix} T_a \\ T_s \end{pmatrix} \in \mathcal V ^{eq,\varepsilon _a ^*} \end{cases}
\,\, \iff \,\, \begin{cases} \begin{pmatrix} T_a \\ T_s \end{pmatrix} = \tilde \Psi ^{(\lambda)} (\varepsilon _a) , \\ \varepsilon _a \in \mathcal V^* . \end{cases}$$
Thus, there exists a unique equilibrium for problem \eqref{2layer-pbm-EDO} close to $(T_a ^{eq,\varepsilon _a ^*}, T_s ^{eq,\varepsilon _a ^*})$ with $\varepsilon _a$ close to $\varepsilon _a ^*$. We denote such equilibrium by $(T_a ^{eq,\varepsilon _a}, T_s ^{eq,\varepsilon _a})$. Moreover, $T_s ^{eq,\varepsilon _a} \notin [T_{s,-}, T_{s,+}]$ for $\varepsilon _a$ close to $\varepsilon _a ^*$ and $(T_a ^{eq,\varepsilon _a}, T_s ^{eq,\varepsilon _a})$ is asymptotically stable (by continuity of $\det$ and $\text{Tr}$). The analyticity of the functions $\varepsilon _a \mapsto T_s ^{eq,\varepsilon _a}$
and $\varepsilon _a \mapsto T_a ^{eq,\varepsilon _a}$ can be deduced by the analytic version of the implicit function theorem, see for instance \cite[Proposition 6.1, p. 138]{Cartan}.

To study the monotonicity 
of $\varepsilon _a \mapsto T_s ^{eq,\varepsilon _a}$ we compute the derivative of $H$ with respect to $\eps_a$:
$$ 0 = \frac{d}{d\varepsilon _a} H (\varepsilon _a,\tilde \Psi ^{(\lambda)} (\varepsilon _a) ) = 
\frac{d}{d\varepsilon _a} \Bigl[ H (\varepsilon _a, \begin{pmatrix} T_a ^{eq,\varepsilon _a} \\ T_s ^{eq,\varepsilon _a} \end{pmatrix} ) \Bigr] ,$$
and we obtain that
$$ \begin{pmatrix} \frac{\partial T_a ^{eq,\varepsilon _a}}{\partial \varepsilon _a} \\ \frac{\partial T_s ^{eq,\varepsilon _a}}{\partial \varepsilon _a} \end{pmatrix}
= - D_2 H (\varepsilon _a, \begin{pmatrix} T_a ^{eq,\varepsilon _a} \\ T_s ^{eq,\varepsilon _a} \end{pmatrix} ) ^{-1} \, D_1 H (\varepsilon _a, \begin{pmatrix} T_a ^{eq,\varepsilon _a} \\ T_s ^{eq,\varepsilon _a} \end{pmatrix} ) .$$
For what concerns the second component, we get that
\begin{equation*}
\begin{split}
\Bigl[& \gamma _a \gamma _s \det D_2 H(\varepsilon _a, \begin{pmatrix} T_a ^{eq,\varepsilon _a} \\ T_s ^{eq,\varepsilon _a} \end{pmatrix} ) \Bigr] \, \frac{\partial T_s ^{eq,\varepsilon _a}}{\partial \varepsilon _a}=
\\
& =  \Bigl( \lambda +4\varepsilon _a \sigma_B \vert T_a ^{eq,\varepsilon _a} \vert ^3 \Bigr) \Bigl( \sigma _B \vert T_s ^{eq,\varepsilon _a}\vert ^4 - 2 \sigma _B \vert T_a ^{eq,\varepsilon _a} \vert ^4 \Bigr)  \\
&\qquad+ \Bigl( \lambda +8\varepsilon _a \sigma_B \vert T_a ^{eq,\varepsilon _a} \vert ^3 \Bigr) \Bigl( \sigma _B \vert T_a ^{eq,\varepsilon _a} \vert ^4 \Bigr) 
\\ 
&=  \lambda \sigma _B \Bigl(  \vert T_s ^{eq,\varepsilon _a}\vert ^4 - \vert T_a ^{eq,\varepsilon _a} \vert ^4 \Bigr) + 4\varepsilon _a \sigma_B ^2 \, \vert T_a ^{eq,\varepsilon _a} \vert ^3 \, \vert T_s ^{eq,\varepsilon _a}\vert ^4  ,
\end{split}
\end{equation*}
which is positive. Therefore, $T_s ^{eq,\varepsilon _a}$ is increasing with respect to $\varepsilon _a$.

With easy computations on the first component, we obtain that
\begin{equation*}
\begin{split}
\Bigl[ &\frac{\gamma _a \gamma _s}{\sigma _B} \det D_2 H(\varepsilon _a, \begin{pmatrix} T_a ^{eq,\varepsilon _a} \\ T_s ^{eq,\varepsilon _a} \end{pmatrix} ) \Bigr] \, \frac{\partial T_a ^{eq,\varepsilon _a}}{\partial \varepsilon _a}=
\\
&= 
\Bigl( \lambda +4 \sigma_B \vert T_s ^{eq,\varepsilon _a} \vert ^3 \Bigr) \Bigl( \vert T_s ^{eq,\varepsilon _a} \vert ^4
- 2 \vert T_a ^{eq,\varepsilon _a}\vert) ^4 \Bigr)  + \Bigl( \lambda + 4 \varepsilon _a \sigma_B \vert T_s ^{eq,\varepsilon _a} \vert ^3 \Bigr) \vert T_a ^{eq,\varepsilon _a} 
\vert ^4 
\\
&= 
4 (\varepsilon _a -1) \sigma _B \vert T_s ^{eq,\varepsilon _a} \vert ^3 \vert T_a ^{eq,\varepsilon _a}\vert ^4 \Bigr)+ (\lambda + 4 \sigma _B \vert T_s ^{eq,\varepsilon _a}\vert ^3) \Bigl( \vert T_s ^{eq,\varepsilon _a}\vert ^4  - \vert T_a ^{eq,\varepsilon _a}\vert ^4 \Bigr) ,
\end{split}
\end{equation*}
and this quantity is positive if $\varepsilon _a \geq 1$.  This concludes the proof of Proposition \ref{prop-comp-lambda=0}. \qed

\noindent {\it Remark. We were not able to determine the sign of $\frac{\partial T_a ^{eq,\varepsilon _a}}{\partial \varepsilon _a}$ if $\varepsilon _a \in (0,1)$ which is not straightforward from the previous expression. It is possible to rewrite the above expression as follows
\begin{equation*}
    \lambda \Bigl( \vert T_s ^{eq,\varepsilon _a}\vert ^4  - \vert T_a ^{eq,\varepsilon _a}\vert ^4 \Bigr)
    + 4 \sigma _B \vert T_s ^{eq,\varepsilon _a}\vert ^3 \Bigl( \vert T_s ^{eq,\varepsilon _a}\vert ^4  - (2-\varepsilon _a) \vert T_a ^{eq,\varepsilon _a}\vert ^4 \Bigr), 
\end{equation*} 
that however does not give a hint of the sign for $\varepsilon _a \in (0,1)$ (observe that this quantity is positive when $\lambda =0$). The knowledge of the sign would be interesting, since $\varepsilon _a \in (0,1)$ is the relevant physical case.
}

\subsection{Proof of Corollary \ref{cor-greenhouse}}\hfill

Let $(T_a ^{eq,\varepsilon _a ^*}, T_s ^{eq,\varepsilon _a ^*})$ be a warm equilibrium, therefore $T_s ^{eq,\varepsilon _a ^*} > T_{s,+}$. Then, thanks to Proposition \ref{prop-comp-lambda=0}, there exists $\delta^*>0$ such that for any $\varepsilon _a\in(\varepsilon _a ^* - \delta ^*,
\varepsilon _a ^* + \delta ^*)$ there exists a unique warm equilibrium 
$(T_a ^{eq,\varepsilon _a }, T_s ^{eq,\varepsilon _a})$. Moreover, $\varepsilon _a \mapsto T_s ^{eq,\varepsilon _a}$ is increasing and $T_s ^{eq,\varepsilon _a} >  T_{s,+}$
for all $\varepsilon _a \in [\varepsilon _a ^* ,
\varepsilon _a ^* + \delta ^*)$.

One can iterate the procedure by applying Proposition \ref{prop-comp-lambda=0} for any $(T_a ^{eq,\varepsilon _a }, T_s ^{eq,\varepsilon _a})$ with $\eps_a\in[\varepsilon _a ^* ,
\varepsilon _a ^* + \delta ^*)$. Our aim is to prove that there exists a unique warm equilibrium for any $\eps_a\in [\varepsilon _a ^*,2)$

Define
\begin{equation*}
\varepsilon _{a,max} := \sup \{ \tilde \varepsilon _a \in [\varepsilon _a ^*,2)\,:\,
\exists \text{ a warm equilibrium for all } \varepsilon _a \in [\varepsilon _a ^*, \tilde \varepsilon _a ]\} .
\end{equation*}
Then, there exists a unique warm equilibrium for all $\varepsilon _a \in [\varepsilon _a ^*, \varepsilon _{a,max})$ and function $\varepsilon _a \in [\varepsilon _a ^*, \varepsilon _{a,max}) \mapsto T_s ^{eq,\varepsilon _a}$ is increasing. Thus, either $\varepsilon _{a,max}=2$ and the proof is completed or $\varepsilon _{a,max} <2$. Assume that $\varepsilon _{a,max} <2$. Observe that, as we have proved in Section \ref{sec-uniform-bounds}, function $\varepsilon _a \in [\varepsilon _a ^*, \varepsilon _{a,max}) \mapsto T_s ^{eq,\varepsilon _a}$ is bounded and therefore $T^{eq,\eps_s}_s\to T_{s,max}$ as $\varepsilon _a \to \varepsilon _{a,max}$. Moreover, function $\varepsilon _a \in [\varepsilon _a ^*, \varepsilon _{a,max}) \mapsto T_a ^{eq,\varepsilon _a}$ is also bounded and so there exists a subsequence $\varepsilon _{a,n} \to \varepsilon _{a,max}$ such that $T_a ^{eq,\varepsilon _{a,n}}\to T_{a,max}$. Taking the limit $\varepsilon _{a,n} \to \varepsilon _{a,max}$ in \eqref{2layer-pbm-EDO-eq*}, we deduce that $(T_{a,max}, T_{s,max})$ is an equilibrium point associated to $\eps_{a,max}$. Therefore, we have just proved that there exists a warm equilibrium associated to the parameter $\varepsilon _{a,max}$.
We can then apply Proposition \ref{prop-comp-lambda=0} and deduce that there exists $\delta>0$ such that, for every $\eps_a\in(\eps_{a,max}-\delta,\eps_{a,max}+\delta)$, there exists a unique warm equilibrium point. However, this contradicts the definition of $\varepsilon _{a,max}$. Therefore $\varepsilon _{a,max} =2$. 

Finally, the analyticity of the functions $\varepsilon _a \in [\varepsilon _a ^*, 2) \mapsto T_s ^{eq,\varepsilon _a}$
and $\varepsilon _a \in [\varepsilon _a ^*, 2) \mapsto T_a ^{eq,\varepsilon _a}$ can be deduced by the analytic version of the implicit function theorem, see, e.g., \cite[Proposition 6.1, p. 138]{Cartan}. This proves the first item of Corollary \ref{cor-greenhouse}.

To prove the second item, let the former warm equilibrium point $(T_a ^{eq, \varepsilon _a ^*}, T_s ^{eq, \varepsilon _a ^*})$ be the initial condition of our system with parameter $\varepsilon _a ^+ > \varepsilon _a ^*$. 
Since $(T_a ^{eq, \varepsilon _a ^*}, T_s ^{eq, \varepsilon _a ^*})$ is an equilibrium of problem \eqref{2layer-pbm-EDO} with parameter $\varepsilon _a ^*$, we have
\begin{equation}
\label{eq-eq} 
\left\{ 
\begin{split}
- \lambda (T_a ^{eq, \varepsilon _a ^*} - T_s ^{eq, \varepsilon _a ^*} )
&= 2 \varepsilon _a ^* \sigma _B \vert T_a ^{eq, \varepsilon _a ^*} \vert ^4 - \varepsilon _a ^* \sigma _B \vert T_s ^{eq, \varepsilon _a ^*} \vert ^4 ,
\\
- \lambda (T_s ^{eq, \varepsilon _a ^*} - T_a ^{eq, \varepsilon _a ^*} )
&= \sigma _B \vert T_s ^{eq, \varepsilon _a ^*} \vert ^4 - \varepsilon _a ^* \sigma _B \vert T_a ^{eq, \varepsilon _a ^*} \vert ^4 - q \beta _s (T_s ^{eq, \varepsilon _a ^*}) .
\end{split}
\right.
\end{equation}
Thanks to the first identity of \eqref{eq-eq}, we obtain
\begin{equation} 
\label{eq-eq_a}
\begin{split}
\gamma _a T_a '(0)&= - \lambda (T_a ^{eq, \varepsilon _a ^*} - T_s ^{eq, \varepsilon _a ^*} ) + \varepsilon _a ^+ \sigma _B \vert T_s ^{eq, \varepsilon _a ^*} \vert ^4 - 2 \varepsilon _a ^+ \sigma _B \vert T_a ^{eq, \varepsilon _a ^*} \vert ^4
\\
&= 2 \varepsilon _a ^* \sigma _B \vert T_a ^{eq, \varepsilon _a ^*} \vert ^4 - \varepsilon _a ^* \sigma _B \vert T_s ^{eq, \varepsilon _a ^*} \vert ^4 + \varepsilon _a ^+ \sigma _B \vert T_s ^{eq, \varepsilon _a ^*} \vert ^4 - 2 \varepsilon _a ^+ \sigma _B \vert T_a ^{eq, \varepsilon _a ^*} \vert ^4
\\
&= ( \varepsilon _a ^+ -\varepsilon _a ^*) \sigma _B \vert T_s ^{eq, \varepsilon _a ^*} \vert ^4
- 2 (\varepsilon _a ^+ - \varepsilon _a ^*) \sigma _B \vert T_a ^{eq, \varepsilon _a ^*} \vert ^4
\\
&= \sigma _B ( \varepsilon _a ^+ - \varepsilon _a ^*) \Bigl( \vert T_s ^{eq, \varepsilon _a ^*} \vert ^4 - 2 \vert T_a ^{eq, \varepsilon _a ^*} \vert ^4\Bigr) < 0 .
\end{split}
\end{equation}
And from the second identity of \eqref{eq-eq}, we deduce that
\begin{equation}
\label{eq-eq_s}
\begin{split}
\gamma _s T_s '(0)&= - \lambda (T_s ^{eq, \varepsilon _a ^*} - T_a ^{eq, \varepsilon _a ^*} ) - \sigma _B \vert T_s ^{eq, \varepsilon _a ^*} \vert ^4
+ \varepsilon _a ^+ \sigma _B \vert T_a ^{eq, \varepsilon _a ^*} \vert ^4+ q \beta _s (T_s ^{eq, \varepsilon _a ^*})
\\
&= \sigma _B \vert T_s ^{eq, \varepsilon _a ^*} \vert ^4 - \varepsilon _a ^* \sigma _B \vert T_a ^{eq, \varepsilon _a ^*} \vert ^4 - q \beta _{s,+} - \sigma _B \vert T_s ^{eq, \varepsilon _a ^*} \vert ^4 + \varepsilon _a ^+ \sigma _B \vert T_a ^{eq, \varepsilon _a ^*} \vert ^4 + q \beta _{s,+}
\\ 
&= ( \varepsilon _a ^+ - \varepsilon _a ^*) \sigma _B \vert T_a ^{eq, \varepsilon _a ^*} \vert ^4 >0 .
\end{split}
\end{equation}
Let us define the subsets of $\mathcal{Q}$
\begin{multline*}
    \mathcal Q_1 ^{(\lambda, \varepsilon _a ^+)}:= \{ (T_a, T_s) \in \mathcal Q,
T_s \in (0, T_{s,1}^{eq,\varepsilon _a ^+})\,:\, \\
\begin{cases}
-\lambda (T_a-T_s) + \varepsilon _a ^+ T_s ^4 -2  \varepsilon _a ^+  T_a ^4 >0 , \\
-\lambda (T_s-T_a)- \sigma _B T_s ^4  + \varepsilon _a ^+ \sigma _B T_a ^4 + q\beta _s(T_s) > 0 ,
\end{cases}
 \} ,
\end{multline*}
\begin{multline*}
\mathcal Q_1  ^{',(\lambda, \varepsilon _a ^+)} := \{ (T_a, T_s) \in \mathcal Q, 
T_s \in (T_{s,1}^{eq,\varepsilon _a ^+}, T_{s,2}^{eq,\varepsilon _a ^+})\,:\,\\
\begin{cases}
-\lambda (T_a-T_s) + \varepsilon _a ^+ T_s ^4 -2  \varepsilon _a ^+  T_a ^4 <0 , \\
-\lambda (T_s-T_a)- \sigma _B T_s ^4  + \varepsilon _a ^+ \sigma _B T_a ^4 + q\beta _s(T_s) < 0 ,
\end{cases}
 \} ,
\end{multline*}
\begin{multline*}
\mathcal Q_2  ^{(\lambda, \varepsilon _a ^+)}:= \{ (T_a, T_s) \in \mathcal Q\,:\,\\
\begin{cases}
-\lambda (T_a-T_s) + \varepsilon _a ^+ T_s ^4 -2  \varepsilon _a ^+  T_a ^4 >0 , \\
-\lambda (T_s-T_a)- \sigma _B T_s ^4  + \varepsilon _a ^+ \sigma _B T_a ^4 + q\beta _s(T_s) < 0 ,
\end{cases}
 \} ,
\end{multline*}
\begin{multline*}
\mathcal Q_3  ^{(\lambda, \varepsilon _a ^+)}:= \{ (T_a, T_s) \in \mathcal Q, 
T_s >T_{s,3}^{eq,\varepsilon _a ^+}\,:\,\\
\begin{cases}
-\lambda (T_a-T_s) + \varepsilon _a ^+ T_s ^4 -2  \varepsilon _a ^+  T_a ^4 <0 , \\
-\lambda (T_s-T_a)- \sigma _B T_s ^4  + \varepsilon _a ^+ \sigma _B T_a ^4 + q\beta _s(T_s) < 0 ,
\end{cases}
\} ,
\end{multline*}
\begin{multline*}
\mathcal Q_3  ^{',(\lambda, \varepsilon _a ^+)} := \{ (T_a, T_s) \in \mathcal Q, 
T_s \in (T_{s,2}^{eq,\varepsilon _a ^+}, T_{s,3}^{eq,\varepsilon _a ^+} )\,:\,\\
\begin{cases}
-\lambda (T_a-T_s) + \varepsilon _a ^+ T_s ^4 -2  \varepsilon _a ^+  T_a ^4 >0 , \\
-\lambda (T_s-T_a)- \sigma _B T_s ^4  + \varepsilon _a ^+ \sigma _B T_a ^4 + q\beta _s(T_s) > 0 ,
\end{cases}
\} ,
\end{multline*}
\begin{multline*}
\mathcal Q_4  ^{(\lambda, \varepsilon _a ^+)} := \{ (T_a, T_s) \in \mathcal Q\,:\,\\
\begin{cases}
-\lambda (T_a-T_s) + \varepsilon _a ^+ T_s ^4 -2  \varepsilon _a ^+  T_a ^4 <0 , \\
-\lambda (T_s-T_a)- \sigma _B T_s ^4  + \varepsilon _a ^+ \sigma _B T_a ^4 + q\beta _s(T_s) > 0 ,
\end{cases}
 \} .
\end{multline*}
Observe that from \eqref{eq-eq_a}-\eqref{eq-eq_s} we deduce that
$(T_a ^{eq, \varepsilon _a ^*}, T_s ^{eq, \varepsilon _a ^*})$ belongs to the subset $\mathcal Q_4  ^{(\lambda, \varepsilon _a ^+)}$
(that generalizes $Q_4$ in Figure \ref{fig8:region}). Since $T_s (0)=T_s ^{eq, \varepsilon _a ^*} > T_{s,+}$, then as long as $(T_a,T_s)$ remains in $\mathcal Q_4  ^{(\lambda, \varepsilon _a ^+)}$, $T_s$ increases and therefore it will remain 
bigger than $T_{s,+}$. Then, several possibilities can happen
\begin{itemize} 
\item either $(T_a,T_s)$ remains forever in $\mathcal Q_4  ^{(\lambda, \varepsilon _a ^+)}$. In this case $T_s$ is increasing and $(T_a,T_s)$ converges to the new warm equilibrium,

\item or it enters into $\mathcal Q_3  ^{'(\lambda, \varepsilon _a ^+)}$. In this case the solution increasingly attains the new warm equilibrium,

\item or it enters into $\mathcal Q_3  ^{(\lambda, \varepsilon _a ^+)}$. Hence, the solution decreasingly converges to the new warm equilibrium.
\end{itemize}
This concludes the proof of Corollary \ref{cor-greenhouse}. \qed

\section{Concluding remarks}
\label{conclusions}
First, we observe that the mathematical analysis of the problem we present in this paper can be completed with some physics-informed comments. Indeed, showing that having larger values of $\eps_a$ leads to a warmer temperature at surfaces amounts to mathematically proving the greenhouse effect for this simple model. 
In the case of a very opaque atmosphere,  if one wants to proceed along the lines of Eq. \ref{2layer-pbm}, in order to respect the basic laws of thermodynamics, one needs to consider multiple layers that behave radiatively as black bodies stacked on top of each other. 

Let us nonetheless assume that if one considers a model  with a single atmospheric layer having (an unphysical) $\epsilon_a\geq1$, amounts to representing by and large a very opaque atmosphere. Then the fact that the sensitivity of the atmospheric temperature with respect to $\eps_a$ is positive for $1<\eps_a<2$, meaning that larger values of $\eps_a$ lead to higher atmospheric temperature, can be roughly interpreted as a  signature of the runaway greenhouse effect, which indeed manifests itself when $\eps_a\geq2$  leading to a blow-up in finite time of the solution. 

Moreover, larger values of $\lambda$ correspond to having a stronger coupling between the surface layer and the atmosphere (the two temperatures being identical, ceteris paribus, in the $\lambda\rightarrow\infty$ limit). Since the incoming solar radiation is primarily absorbed at surface, larger values of $\lambda$ allow for a more efficient upward transfer of energy to the atmosphere, thus decreasing the surface temperature, as a result of enhanced vertical sensible or latent heat transport.

Finally, we would like to point out that the results of this paper will be used in the study of the 1D two-layer energy balance model \eqref{2layer-pbm}, which we will develop in future works. The key mathematical questions are the same: global existence, long time behaviour, influence of the different parameters. However, problem \eqref{2layer-pbm} is composed of two nonlinear and degenerate parabolic equations, which brings many mathematical challenges. For such partial differential equations, comparison results will probably derive bounds from the ODE system studied in the present work. Furthermore, problem \eqref{2layer-pbm} is particularly interesting when the function $q$ is not constant, and depends on $x$ (the insolation being not the same at the poles or at the equator), and then the study of the equilibrium points is one of the interesting challenges.



\appendix

\section{Proof of Proposition \ref{prop-blowup}}
\label{sec-blowup-proof}

Let us prove that blow up in finite time occurs for problem \eqref{2layer-pbm-EDO} when $\varepsilon _a >2$, first under quite restrictive hypothesis (part a) of Proposition \ref{prop-blowup}), then under weaker assumptions (part b) of Proposition \ref{prop-blowup}).

\subsection{Proposition \ref{prop-blowup}, part a)}\hfill

Consider the following problem
\begin{equation}
\label{2layer-pbm-EDO-blow}
\begin{cases}
\gamma _a T_a' = \varepsilon _a \sigma _B \vert T_s \vert ^3 T_s  - 2 \varepsilon _a \sigma _B \vert T_a \vert ^3 T_a , \\
\gamma _s T_s' = - \sigma _B \vert T_s \vert ^3 T_s  + \varepsilon _a \sigma _B \vert T_a \vert ^3 T_a + \mathcal R_s (T_s), \\
T_a (0)=T_a ^{(0)} ,\\
T_s (0)=T_s ^{(0)} ,
\end{cases}
\end{equation}
when $\beta_s$ is given by \eqref{beta}. Observe that, under assumptions  \eqref{global-ci},\eqref{global-coeffs-bl}-\eqref{global-eps-bl} one can still carry out the same proof of Proposition \ref{prop-wellposed} to show that a unique maximal solution of problem \eqref{2layer-pbm-EDO} exists for $t\in [0, T^*)$. Moreover, thanks to Lemma \ref{lem-postivity}, such solution is positive on $(0, T^*)$. 

On the other hand, the boundedness of solutions established in Section \ref{sec-boundedness} was based on the fact that $\varepsilon _a < 2$. Thus, this property cannot be deduced when $\varepsilon _a \geq 2$.

\subsubsection{Lack of equilibrium points} \hfill

Equilibrium points are stationary solutions of \eqref{2layer-pbm-EDO-blow}: they solve the system
\begin{equation}
\label{equil-bl} \begin{cases}
\varepsilon _a \sigma _B \vert T_s \vert ^3 T_s  - 2 \varepsilon _a \sigma _B \vert T_a \vert ^3 T_a = 0 , \\
- \sigma _B \vert T_s \vert ^3 T_s  + \varepsilon _a \sigma _B \vert T_a \vert ^3 T_a + \mathcal R_s = 0 .
\end{cases}
\end{equation}
Since $T_s \geq 0$ and $T_a \geq 0$, from the first equation of \eqref{equil-bl} we get
$$ T_s= 2^{1/4} T_a,$$
and, plugging such identity into the second equation of \eqref{equil-bl}, we obtain
$$ (\varepsilon _a -2) \sigma _B T_a ^4 + q \beta _s (T_s)=0,$$
which has no solutions because $\varepsilon _a \geq 2$ and $q\beta _s >0$.

{\it Remark. Observe that, by replacing
$T_a$ by $2^{-1/4} T_s$ in the ODE \eqref{2layer-pbm-EDO-blow} satisfied by $T_s$, we get
$$ \gamma _s T_s' = (\frac{\varepsilon _a}{2} -1) \sigma _B  T_s ^4  + \mathcal R_s (T_s) ,$$
where blow up in finite time occurs when $\frac{\varepsilon _a}{2} -1>0$. Of course, this procedure
is not rigorous, but anyhow it gives an insight of blow up in finite time for $\varepsilon _a >2$.

}

\subsubsection{Monotonicity of the solution} \hfill

Now let us consider the sets
$$ \mathcal C _1 = \{(T_a,T_s) \in \overline{\mathcal Q}, \varepsilon _a \sigma _B \vert T_s \vert ^3 T_s  - 2 \varepsilon _a \sigma _B \vert T_a \vert ^3 T_a = 0 \},$$
and 
$$ \mathcal C _2 = \{(T_a,T_s)\in \overline{\mathcal Q}, - \sigma _B \vert T_s \vert ^3 T_s  + \varepsilon _a \sigma _B \vert T_a \vert ^3 T_a + \mathcal R_s (T_s) =0 \}.$$
Observe that $(T_a,T_s) \in \mathcal C _1$ if and only if $T_s=2^{1/4} T_a$ with $T_a \geq 0$. Given $T_s \geq 0$, we will denote $T_a ^{(1)} (T_s)$ the value for which $(T_a ^{(1)} (T_s), T_s) \in \mathcal C _1$. Hence, we have $T_a ^{(1)} (T_s) =  2^{-1/4} T_s$.

Now, let us analyse $\mathcal C _2$. Since there are no equilibrium points in $\overline{\mathcal Q}$, $\mathcal C _2$ does not intersect $\mathcal C _1$ in $\overline{\mathcal Q}$. 

Moreover, $\mathcal C _2$ contains points of the form $(0,T_s)$. Indeed, such points satisfy
\begin{equation}
\label{blowup-bord}
 - \sigma _B \vert T_s \vert ^3 T_s + \mathcal R_s (T_s) = 0, 
\end{equation}
which has at least one solution (function $T_s \mapsto - \sigma _B \vert T_s \vert ^3 T_s + \mathcal R_s (T_s)$ is positive for $T_s$ close to $0$ and negative for $T_s$ large). More precisely, we have proved in Section \ref{sec-comp-asympt-lambda=0} that, when $\beta _s$ is given by \eqref{beta}, \eqref{blowup-bord} can have one, two or three solutions $T_s$ depending on the values of the parameters appearing in $\beta_{s}$, on $q$ and $\sigma _B$.

We further observe that if 
$- \sigma _B \vert T_s \vert ^3 T_s + \mathcal R_s (T_s) > 0$, there is no $T_a \geq 0$ such that 
$(T_a, T_s) \in \mathcal C _2$. On the other hand, if $- \sigma _B \vert T_s \vert ^3 T_s + \mathcal R_s (T_s) \leq 0$, there exists a unique value $T_a \geq 0$ such that $(T_a, T_s) \in \mathcal C _2$. We denote such value $T_a ^{(2)} (T_s)$.
Therefore 
$$ \mathcal C _2 = \{ (T_a ^{(2)} (T_s), T_s), \sigma _B \vert T_s \vert ^3 T_s \geq \mathcal R_s (T_s) \}.$$
Note that the set 
$$ \mathcal J := \{ T_s, \sigma _B \vert T_s \vert ^3 T_s \geq \mathcal R_s (T_s) \}$$ 
is a union of intervals where the extremal points of such intervals are zeros of the map $T_s \mapsto \sigma _B \vert T_s \vert ^3 T_s - \mathcal R_s (T_s)$.

Finally, we note that
$$ (T_a,T_s) \in \mathcal C _2 \quad \implies \quad T_s \sim \varepsilon _a ^{1/4} T_a \quad \text{as}\quad T_s \to \infty .$$
This gives the asymptotic shape of $\mathcal C _2$.

We claim that
\begin{equation}
\label{mathcalJ}
\forall\, T_s \in \mathcal J, \quad T_a ^{(2)} (T_s) < T_a ^{(1)} (T_s) .
\end{equation}
Indeed, first we observe that on a compact connected component of $\mathcal J$ it holds that
$$ \forall\, T_s \in \partial \mathcal J, \quad T_a ^{(2)} (T_s) < T_a ^{(1)} (T_s) $$
because $ T_a ^{(2)} (T_s)=0$ on $\partial \mathcal J$. Moreover, there is no $T_s$ such that
$T_a ^{(2)} (T_s)= T_a ^{(1)} (T_s)$. Therefore, by continuity, \eqref{mathcalJ} holds true on all the compact connected components of $\mathcal J$.

\begin{figure}[h!]
\centering
\subfloat[]{
\begin{tikzpicture}[scale=.4]
\draw[-stealth] (-1,0) -- (7,0)node[below]{$T_a$};
\draw[-stealth] (0,-1) -- (0,7)node[left]{$T_s$};
\draw[blue] (0,0) -- (6,7)node[right]{$\mathcal{C}_1$};
\draw[domain=0:1.9, smooth, variable=\x,red] plot ({\x}, {1.1*\x*\x+3});
\node[red] at (0.5,7.5)[right] {$\mathcal{C}_2$};
\draw[red,-stealth] (.68,3.5) -- (1.18,3.5);
\draw[red,-stealth] (1.36,5) -- (1.86,5);
\draw[red,-stealth] (1.8,6.5) -- (2.3,6.5);
\draw[blue,-stealth] (2,2.32) -- (2,2.82);
\draw[blue,-stealth] (3,3.5) -- (3,4);
\draw[blue,-stealth] (4,4.68) -- (4,5.18);
\node at (0,3)[left] {$T_{1,*}$};
\node at (3,5) {$\mathcal{E}$};
\end{tikzpicture}
}\quad
\subfloat[]{
\begin{tikzpicture}[scale=.4]
\draw[-stealth] (-1,0) -- (7,0)node[below]{$T_a$};
\draw[-stealth] (0,-1) -- (0,7)node[left]{$T_s$};
\draw[blue] (0,0) -- (6,7)node[right]{$\mathcal{C}_1$};
\draw[domain=0:1.9, smooth, variable=\x,red] plot ({\x}, {1.1*\x*\x+3});
\node[red] at (0.5,7.5)[right] {$\mathcal{C}_2$};
\draw[red,-stealth] (.68,3.5) -- (1.18,3.5);
\draw[red,-stealth] (1.36,5) -- (1.86,5);
\draw[red,-stealth] (1.8,6.5) -- (2.3,6.5);
\draw[blue,-stealth] (2,2.32) -- (2,2.82);
\draw[blue,-stealth] (3,3.5) -- (3,4);
\draw[blue,-stealth] (4,4.68) -- (4,5.18);
\node at (0,3)[left] {$T_{3,*}$};
\node at (3,5) {$\mathcal{E}$};
\draw[red] (0,1) arc (-90: 90: .7);
\node at (0,1)[left] {$T_{1,*}$};
\node at (0,2.2)[left] {$T_{2,*}$};
\end{tikzpicture}
}	
\caption{In the phase space we represent the sets $\mathcal{C}_1=\{(T_1,T_2)\in\bar{Q}\,:\, \eps_a\sigma_B|T_s|^3T_s-2\eps_a\sigma_B|T_a|^3T_a=0\}$, $\mathcal{C}_2=\{(T_a,T_s)\in\bar{Q}\,:\,-\sigma_B|T_s|^3T_s+\eps_a\sigma_B|T_a|^3T_a+\mathcal{R}_s(T_s)=0\}$ and $\mathcal{E}=\{(T_a,T_s)\in\bar{Q}\,:\,\sigma_BT_s^4-\mathcal{R}_s(T_s)<\eps_a\sigma_BT^4_a<\frac{1}{2}\eps_a\sigma_BT^4_s\}$. In (a) we show the case in which $-\sigma_B|T_s|^3T_s+\mathcal{R}_s(T_s)=0$ has a unique solution and in (b) the case of three solutions.}
\label{fig1:region-blow}
\end{figure}
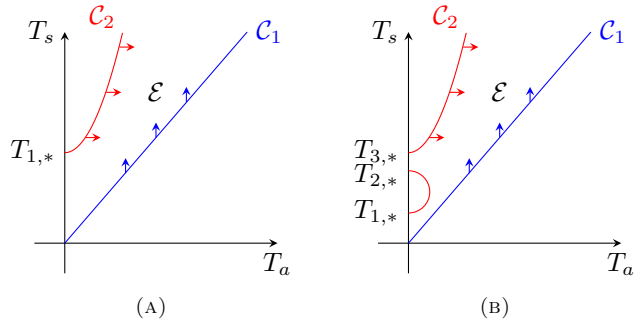

Furthermore, on the unbounded connected component, since $\varepsilon _a >2$, the asymptotics of $T_a ^{(2)} (T_s)$ and $ T_a ^{(1)} (T_s)$ give that
$$\text{ $T_s$ large enough} \quad \implies \quad T_a ^{(2)} (T_s) < T_a ^{(1)} (T_s).$$
Since there is no $T_s$ such that
$T_a ^{(2)} (T_s)= T_a ^{(1)} (T_s)$, the sign of $T_a ^{(2)} (T_s) - T_a ^{(1)} (T_s)$ cannot change
in the unbounded connected component (by continuity), and therefore \eqref{mathcalJ} holds true.
This implies that $\mathcal C_2$ remains strictly above $\mathcal C_1$.
In the following we will consider the set
$$\mathcal E := \{ (T_a, T_s) \in \overline{\mathcal Q}, \sigma _B T_s ^4 - \mathcal R_s (T_s) < \varepsilon _a \sigma _B T_a ^4 < \frac{1}{2} \varepsilon _a \sigma _B T_s ^4 \},$$
see Figure \ref{fig1:region-blow}.

Now, let us rewrite the vector field $F$ defined in \eqref{def-Fnonlin} under the current assumptions
\begin{equation}
\label{def-Fnonlin-blow}
F \begin{pmatrix} T_a \\ T_s \end{pmatrix}
= \begin{pmatrix} \frac{1}{\gamma _a} \Bigl[ \varepsilon _a \sigma _B \vert T_s \vert ^3 T_s  - 2 \varepsilon _a \sigma _B \vert T_a \vert ^3 T_a \Bigr] \\ \frac{1}{\gamma _s} \Bigl[ - \sigma _B \vert T_s \vert ^3 T_s  + \varepsilon _a \sigma _B \vert T_a \vert ^3 T_a + \mathcal R_s (T_s)\Bigr] \end{pmatrix} .
\end{equation}
Assume that $T_s ^{(0)} < 2^{1/4} T_a ^{(0)}$. Then, as long as
the solution $(T_a, T_s)$ remains below $\mathcal C _1$, we have $T_a ' \leq 0$ and $T_s ' >0$,
hence, in the phase plane, the solution goes monotonically in the \lq\lq north-west" direction, and attains $\mathcal C _1$ in finite time (otherwise, it would have to converge to some equilibrium point which, however, does not exist). Hence, under assumption $T_s ^{(0)} < 2^{1/4} T_a ^{(0)}$, the solution enters $\mathcal E$.

Analogously, if the initial condition is above $\mathcal C_1$ and not in $\mathcal E$, then it enters $\mathcal E$ in finite time. And if the initial condition belongs to $\mathcal E$, then the solution cannot leave $\mathcal E$, since the vector field goes inward $\mathcal E$.

Therefore, no matter where the initial condition is located in $\bar{Q}$, the solution enters $\mathcal E$ (in a monotonical way). When it is in $\mathcal E$, then 
$T_a ' >0$ and $T_s ' >0$, that is, $T_a$ and $T_s$ are increasing. Therefore, since the temperatures cannot converge to some equilibrium point, they both go to $+\infty$ as $t\to \tau _{a,s} ^+$, which is the maximum existence time (see Section \ref{sec-local-existence}) that, at this stage of the proof, can be finite or infinite.


\subsubsection{Asymptotic behaviour in the phase plane}\hfill

As already noted, the solution, once entered in $\mathcal E$, remains in such set for all time of existence. Hence, there exists $\tau _0$ such that
$$ \forall  \, t \in [\tau_0, \tau _{a,s} ^+), \quad T_a ^{(2)} (T_s (t))  < T _a (t) < T_a ^{(1)} (T_s (t)) .$$
Since $ T_a ^{(1)} (T_s (t)) = 2^{-1/4} T_s (t)$, and $ T_a ^{(2)} (T_s (t))  \sim \varepsilon _a ^{-1/4} T_s (t)$
as $t \to \tau _{a,s} ^+$, then the quotient $\frac{T_s (t)}{T_a(t)}$ remains bounded between two positive constants.

Assume for the moment that the behaviour is perfectly linear, that is, there exists $\mu _*>0$ such that
\begin{equation}
\label{hope-mu*}
\forall\, t \in [\tau_0, \tau _{a,s} ^+), \quad T_s (t) = \mu _* T_a (t) .
\end{equation}
Then, since $(T_a,T_s) \in \mathcal E$, this linear behaviour would imply that $\mu _* \in [2^{1/4}, \varepsilon _a ^{1/4}]$, and moreover
$$\mu _* =  \frac{T_s '(t)}{T_a '(t)}.$$
From \eqref{2layer-pbm-EDO-blow} we would have
\begin{equation*}
\begin{split}
\frac{\gamma _s}{\gamma _a} \frac{T_s '(t)}{T_a '(t)}
&= \frac{- \sigma _B \vert T_s \vert ^3 T_s  + \varepsilon _a \sigma _B \vert T_a \vert ^3 T_a + \mathcal R_s}{\varepsilon _a \sigma _B \vert T_s \vert ^3 T_s  - 2 \varepsilon _a \sigma _B \vert T_a \vert ^3 T_a} 
\\
&= \frac{- \sigma _B \mu _* ^4 T_a ^4  + \varepsilon _a \sigma _B T_a ^4 + \mathcal R_s}{\varepsilon _a \sigma _B \mu _* ^4  T_a ^4  - 2 \varepsilon _a \sigma _B T_a ^4} \\
&= \frac{\varepsilon _a - \mu _* ^4 + \frac{\mathcal R_s}{\sigma _B T_a ^4}}{\varepsilon _a (\mu _* ^4 - 2 )} 
\to \frac{\varepsilon _a - \mu _* ^4 }{\varepsilon _a (\mu _* ^4 - 2 )},\quad\text{as}\quad T_a \to +\infty.
\end{split}
\end{equation*}
Thus, $\mu_*$ would solve the equation
$$ \frac{\gamma _s}{\gamma _a} \mu_* = \frac{\varepsilon _a - \mu _* ^4 }{\varepsilon _a (\mu _* ^4 - 2 )}, 
\quad \text{ with } \quad \mu _* \in [2^{1/4}, \varepsilon _a ^{1/4}].$$
We observe that the map $y\in [2, \varepsilon _a] \mapsto \frac{\gamma _s}{\gamma _a} y^{1/4}$ is strictly increasing,
whereas the map $y\in (2, \varepsilon _a) \mapsto \frac{\varepsilon _a - y }{\varepsilon _a (y - 2 )}$
is decreasing, it goes to $+\infty$ as $y\to 2^+$, and to $0$ as $y\to \varepsilon _a$. Therefore, there exists a unique value $y_* \in (2, \varepsilon _a)$ such that 
$$ \frac{\gamma _s}{\gamma _a} y_* ^{1/4} = \frac{\varepsilon _a - y_* }{\varepsilon _a (y_* - 2 )},$$
and so a unique value $\mu_* \in (2^{1/4}, \varepsilon _a ^{1/4})$ such that
\begin{equation}
\label{value-mu*}
\frac{\gamma _s}{\gamma _a} \mu_* = \frac{\varepsilon _a - \mu _* ^4 }{\varepsilon _a (\mu _* ^4 - 2 )} .
\end{equation}
Such value is $\mu_* = y_* ^{1/4}$. To sum up, if \eqref{hope-mu*} would hold, $\mu _*$ would have to be the unique value in $(2^{1/4}, \varepsilon _a ^{1/4})$ solving \eqref{value-mu*}.

Of course, the linear behavior \eqref{hope-mu*} can not be ensured. However, from the above argument we are able to prove -- in the general case -- the following bound from below of the ratio $\frac{T_s(t)}{T_a(t)}$.
 
\begin{Lemma}
\label{lem-comp-lin}
Let $\tau _0$ be any value for which $(T_a (\tau_0), T_s (\tau_0)) \in \mathcal E$. Then, 
there exists $\mu \in (2^{1/4}, \mu_*)$ such that
\begin{equation}
\label{eq-complin}
\forall\, t\geq \tau _0, \quad T_s (t) \geq \mu T_a (t) .
\end{equation}
\end{Lemma}
What is important in Lemma \ref{lem-comp-lin} is the fact that \eqref{eq-complin} holds for some $\mu > 2^{1/4}$. This fact is crucial to prove that blow up in finite time occurs.

\noindent {\it Proof of Lemma \ref{lem-comp-lin}.} Let us consider system \eqref{2layer-pbm-EDO-blow} as a first order non-autonomous ODE in $T_s$, with $T_a$ the new variable. Indeed, since the map
$$ \varphi: [\tau_0, \tau _{a,s} ^+) \to [T_a (\tau _0), +\infty), \quad t \mapsto \varphi (t) = T_a (t) $$
is increasing, we can consider its inverse 
$$ \varphi ^{-1}: [T_a (\tau _0), +\infty) \to [\tau_0, \tau _{a,s} ^+) , \quad z \mapsto \varphi ^{-1} (z) ,$$
and then consider 
$$ u: [T_a (\tau _0), +\infty) \to [T_s (\tau _0), + \infty ), \quad u(z) = T_s (\varphi ^{-1} (z)).$$
We have
\begin{equation*}
\begin{split}
&u'(z) = T_s ' (\varphi ^{-1} (z)) \ (\varphi ^{-1})' (z)
=\frac{ T_s ' (\varphi ^{-1} (z)) }{\varphi ' (\varphi ^{-1} (z))}
= \frac{ T_s ' (\varphi ^{-1} (z)) }{T_a ' (\varphi ^{-1} (z))}
\\ 
&= \frac{\gamma _a}{\gamma _s}  \, \frac{- \sigma _B  T_s (\varphi ^{-1} (z)) ^4  + \varepsilon _a \sigma _B T_a (\varphi ^{-1} (z)) ^4 + \mathcal R_s (T_s (\varphi ^{-1} (z)))}{\varepsilon _a \sigma _B T_s (\varphi ^{-1} (z)) ^4  - 2 \varepsilon _a \sigma _B T_a (\varphi ^{-1} (z)) ^4 }
\\ 
&= \frac{\gamma _a}{\gamma _s}  \, \frac{- \sigma _B  u(z) ^4  + \varepsilon _a \sigma _B z ^4 + \mathcal R_s (u(z)) }{\varepsilon _a \sigma _B u(z) ^4  - 2 \varepsilon _a \sigma _B z^4} \\
&= \frac{\gamma _a}{\gamma _s}  \, \frac{\varepsilon _a z^4 -  u(z) ^4   + \frac{\mathcal R_s (u(z))}{\sigma _B } }{\varepsilon _a  (u(z) ^4 -2 z^4)} .
\end{split}
\end{equation*}
Note that 
$$ \frac{ u(z) }{z} = \frac{ T_s (\varphi ^{-1} (z))}{ \varphi(\varphi ^{-1} (z))}
= \frac{ T_s (\varphi ^{-1} (z))}{T_a (\varphi ^{-1} (z))} > 2 ^{1/4} $$
since $(T_a, T_s)$ is above $\mathcal C_1$. Therefore the quantity $u(z) ^4-2 z^4$ is positive. 
Let us now consider the set
$$ U:= \{ (z,x), z\in [T_a (\tau _0), +\infty), x> 2^{1/4}z \},$$
and the function 
$$ f: U \to \mathbb{R}, \quad f(z,x) := \frac{\gamma _a}{\gamma _s}  \,  \frac{\varepsilon _a z^4 -  x ^4   + \frac{\mathcal R_s (x)}{\sigma _B } }{\varepsilon _a  (x ^4 -2 z^4)} .$$
Thus, $u$ satisfies the following differential equation
\begin{equation}
\label{eq-1storder}
\begin{cases} u'(z) = f(z,u(z)), \quad z\in [T_a (\tau _0), +\infty), \\
u(T_a (\tau _0)) = T_s (\tau _0) .\end{cases} 
\end{equation}
Now, since $T_s (\tau _0) > 2^{1/4} T_a (\tau _0)$, we can choose $\mu \in (2^{1/4}, \mu_*)$ -- where $\mu_*$ is the value for which \eqref{value-mu*} holds -- such that $T_s (\tau _0) > \mu T_a (\tau _0)$. Graphically, we are choosing $\mu \in (2^{1/4}, \mu_*)$ so that
the point $(T_a (\tau _0),T_s (\tau _0))$ is above the line $T_s = \mu T_a$. Then, the function
$$ v: [T_a (\tau _0), +\infty) \to \mathbb R, \quad v(z) := \mu z $$
is a subsolution of the first order non-autonomous equation \eqref{eq-1storder}. Indeed,
\begin{equation*}
\begin{split}
 f(z,v(z))&= \frac{\gamma _a}{\gamma _s}  \, \frac{\varepsilon _a z^4 - v(z) ^4   + \frac{\mathcal R_s (v(z))}{\sigma _B } }{\varepsilon _a  (v(z) ^4 -2 z^4)}\\
&= \frac{\gamma _a}{\gamma _s}  \, \frac{(\varepsilon _a - \mu ^4)z^4  + \frac{\mathcal R_s (v(z))}{\sigma _B } }{\varepsilon _a  (\mu ^4 -2) z^4}
\\
&\geq \frac{\gamma _a}{\gamma _s}  \, \frac{(\varepsilon _a - \mu ^4)z^4  }{\varepsilon _a  (\mu ^4 -2) z^4}\\
&= \frac{\gamma _a}{\gamma _s}  \, \frac{\varepsilon _a - \mu ^4}{\varepsilon _a  (\mu ^4 -2)}
\\
&> \frac{\gamma _a}{\gamma _s}  \, \frac{\varepsilon _a - \mu_* ^4}{\varepsilon _a  (\mu _* ^4 -2)} = \mu _* > \mu = v'(z) .
\end{split}
\end{equation*}
Moreover, 
$$ v(T_a (\tau _0)) = \mu T_a (\tau _0) < T_s (\tau _0) = u(T_a (\tau _0)) .$$
To sum up, we have
$$ \begin{cases} u'(z)=f(z,u(z)), \\ u(T_a (\tau _0)) = T_s (\tau_0)
\end{cases} \,\, \text{ and } \,\, 
\begin{cases} v'(z) < f(z,v(z)), \\ v(T_a (\tau _0)) < u(T_a (\tau _0))
\end{cases} .$$
We claim that
$$ \forall z \geq T_a (\tau _0), \quad v(z) < u(z) .$$
Indeed, if there exists $z_1$ such that $v(z_1)=u(z_1)$, we would have
$$ v'(z_1) < f(z_1, v(z_1)) = f(z_1, u(z_1)) = u'(z_1) .$$ 
Such inequality implies that $v-u$ is strictly decreasing near $z_1$ and equal to $0$ at $z_1$. Therefore, $v-u$ would have to be positive before $z_1$, which contradicts the minimality of $z_1$.
Thus, we have that
\begin{equation*}
    \forall\, z \geq T_a (\tau _0)\qquad T_s (\varphi ^{-1} (z)) = u(z) > \mu z = \mu  \varphi(\varphi ^{-1} (z)) 
= \mu T_a (\varphi ^{-1} (z)) ,
\end{equation*}
which implies that
$$ \forall\, t \in [\tau _0, \tau _{a,s} ^+),  \quad T_s (t) > \mu T_a (t) .$$
The proof of Lemma \ref{lem-comp-lin} is therefore complete. \qed


\subsubsection{Blow up in finite time} \hfill

In this section we prove that $\tau _{a,s} ^+ <\infty$. From the first equation in \eqref{2layer-pbm-EDO-blow} we get that
\begin{equation*}
    \gamma _a \, \frac{T_a '(t)}{T_a (t)^4} = \frac{\varepsilon _a \sigma _B \vert T_s \vert ^3 T_s  - 2 \varepsilon _a \sigma _B \vert T_a \vert ^3 T_a}{T_a (t)^4} 
    = \varepsilon _a \sigma _B  \frac{T_s (t) ^4}{T_a (t)^4} - 2 \varepsilon _a \sigma _B .
\end{equation*}
Using Lemma \ref{lem-comp-lin}, we deduce that there exists $\mu\in(2^{1/4},\mu_*)$ such that
$$ \forall\, t \in [\tau _0, \tau _{a,s} ^+), \quad \gamma _a \, \frac{T_a '(t)}{T_a (t)^4} \geq \varepsilon _a \sigma _B  (\mu ^4 -2) .$$
Then, for all $\tau \in [\tau _0, \tau _{a,s} ^+)$, we have
$$ \gamma _a \, \int _{\tau _0} ^\tau \frac{T_a '(t)}{T_a (t)^4} \, dt \geq \int _{\tau _0} ^\tau \varepsilon _a \sigma _B  (\mu ^4 -2) \, dt .$$
By computing the integrals we get
$$ \gamma _a \, \Bigl[ \frac{-1}{3 T_a (t) ^3} \Bigr] _{\tau _0} ^\tau  \geq \varepsilon _a \sigma _B  (\mu ^4 -2) (\tau - \tau _0) ,\quad \forall\, \tau \in [\tau _0, \tau _{a,s} ^+)$$
which gives
\begin{equation*}
     \varepsilon _a \sigma _B  (\mu ^4 -2) (\tau - \tau _0) \leq \frac{\gamma _a}{3 T_a (\tau _0) ^3} - \frac{\gamma _a}{3 T_a (\tau ) ^3}
\leq \frac{\gamma_a}{3 T_a (\tau _0) ^3},\quad\forall\, \tau \in [\tau _0, \tau _{a,s} ^+) .
\end{equation*}
Since $\mu > 2^{1/4}$, we obtain that
$$ \tau - \tau _0 \leq \frac{\gamma_a}{3 \varepsilon _a \sigma _B  (\mu ^4 -2) T_a (\tau _0) ^3} ,\quad\forall\, \tau \in [\tau _0, \tau _{a,s} ^+)$$
which implies that $\tau _{a,s} ^+$ is finite:
$$ \tau _{a,s} ^+ \leq \tau _0 + \frac{\gamma _a}{3 \, \varepsilon _a \, \sigma _B  \, (\mu ^4 -2) \, T_a (\tau _0) ^3} .$$
This concludes the proof of Proposition \ref{prop-blowup}, part a).

\subsection{Proposition \ref{prop-blowup}, part b)} \hfill

In order to deal with the case $\lambda>0$ and $\mathcal{R}_a\geq0$ we consider the set of points
\begin{equation*}
    \mathcal C _1 ^{(g)}= \{(T_a,T_s),\,
    -\lambda (T_a-T_s) +\varepsilon _a \sigma _B \vert T_s \vert ^3 T_s  - 2 \varepsilon _a \sigma _B \vert T_a \vert ^3 T_a + \mathcal R_a (T_a)= 0 \},
\end{equation*}
and 
\begin{equation*}
    \mathcal C _2 ^{(g)}= \{(T_a,T_s),\, -\lambda (T_s -T_a) - \sigma _B \vert T_s \vert ^3 T_s  + \varepsilon _a \sigma _B \vert T_a \vert ^3 T_a + \mathcal R_s (T_s) =0 \}.
\end{equation*}
Observe that for $T_a$ large, there exists a unique value $T_s$ such that 
$$ \lambda T_s +\varepsilon _a \sigma _B \vert T_s \vert ^3 T_s = \lambda T_a  + 2 \varepsilon _a \sigma _B \vert T_a \vert ^3 T_a - \mathcal R_a (T_a) ,$$
that is, $(T_a,T_s)\in \mathcal C_1 ^{(g)}$. Such value, denoted by $T_s ^{(1,g)} (T_a) $, satisfies
$$ T_s ^{(1,g)} (T_a) \sim 2^{1/4} T_a \quad \text{ as } T_a \to \infty .$$
Analogously, given $T_s$ large enough, there is one and only value of $T_a$ such that
$$ \lambda T_a  + \varepsilon _a \sigma _B \vert T_a \vert ^3 T_a = \lambda T_s + \sigma _B \vert T_s \vert ^3 T_s - \mathcal R_s (T_s) ,$$
namely, such that $(T_a,T_s)\in \mathcal C_2 ^{(g)}$. We denote this value $T_a ^{(2,g)} (T_s)$, and we claim that it satisfies
$$ T_a ^{(2,g)} (T_s) \sim \varepsilon _a ^{-1/4} T_s \quad \text{ as } T_s \to \infty .$$
Therefore, we have
\begin{equation*}
    T_a ^{(2,g)} (T_s ^{(1,g)} (T_a)) \sim \varepsilon _a ^{-1/4} T_s ^{(1,g)} (T_a)
\sim \Bigl( \frac{2}{\varepsilon _a} \Bigr) ^{1/4} T_a\qquad \text{ as } T_a \to \infty .
\end{equation*}
Hence, there exists $T_{a,*}$ such that
$$ T_a  \geq T_{a,*} \quad \implies \quad T_a ^{(2,g)} (T_s ^{(1,g)} (T_a  )) < T_a .$$
We now consider the set
$$  \mathcal E ^{(g)} := \bigcup _{T_a \geq T_{a,*}}  (T_a ^{(2,g)} (T_s ^{(1,g)} (T_a  )), T_a) \times \{T_s ^{(1,g)} (T_a ) \},$$
or, in other words, the set of $(T_a, T_s) \in \mathcal Q$ such that $T_a \geq T_{a,*}$ and such that
$$
\begin{cases} 
-\lambda (T_a-T_s) +\varepsilon _a \sigma _B \vert T_s \vert ^3 T_s  - 2 \varepsilon _a \sigma _B \vert T_a \vert ^3 T_a + \mathcal R_a (T_a) >0 ,\\
-\lambda (T_s -T_a) - \sigma _B \vert T_s \vert ^3 T_s  + \varepsilon _a \sigma _B \vert T_a \vert ^3 T_a + \mathcal R_s (T_s) >0 .
\end{cases} $$
Then, we deduce that if the initial condition $(T_a ^{(0)}, T_s ^{(0)})$ belongs to $\mathcal E^{(g)}$, the solution never leaves $\mathcal E^{(g)}$. Furthermore, the components $T_a$ and $T_s$ are increasing and blow up in finite time. \qed

\section{Proof of Lemma \ref{lem-nbre-eq-1*conv-phys} and Remark \ref{lem-nbre-eq-1*conv-simpl}}
\label{sec-conv-Phi1}

\subsection{Preliminary computations}\hfill

When studying the convexity of $\Phi$, a useful tool is a suitable expression of $\Phi ''$.
This is the goal of this section. 

Let us first consider function $\Phi_1$, defined in \eqref{def-Phi1*}.
By computing the second derivative, we get
$$\Phi_1 '' ( T_s) = - \frac{\lambda}{2} \frac{d^2 T_a ^{(1)}}{d T_s ^2} (T_s) .$$
We introduce the following function
\begin{equation}
\label{eq-rho-app}
\rho (T_s) := \frac{T_a ^{(1)} (T_s)}{T_s}.
\end{equation}
Thanks to \eqref{lambda-posi2}, we know that $\rho (T_s) \in (2^{-1/4},1)$.
We replace $T_a ^{(1)} (T_s)$ by $\rho (T_s) T_s$ in the first equation of
\eqref{2layer-pbm-EDO-eq-lambda-equiv} and we obtain the identity
$$ \lambda (1 - \rho (T_s)) T_s = \varepsilon _a \sigma _B (2 \rho (T_s) ^4 - 1) T_s ^4 ,$$
that can be rewritten as follows
\begin{equation}
\label{def-inv-rho}
\frac{2 \rho (T_s) ^4 - 1}{1 - \rho (T_s)} = \frac{\lambda}{\varepsilon _a \sigma _B T_s ^3} .
\end{equation}
Let us introduce the map
\begin{equation}
\label{def-L-K}
L: x\in [2^{-1/4},1) \mapsto L(x) := \frac{2 x ^4 - 1}{1 - x} , \,\,\text{ and } \,\, K_{\text{ph}} := \frac{\lambda}{\varepsilon _a \sigma _B } .
\end{equation}
Observe that $L$ is strictly increasing on $[2^{-1/4},1)$ and it holds that
\begin{equation*}
\label{def-rho}
\rho (T_s) = L^{-1} \Bigl( \frac{ K_{\text{ph}} }{ T_s ^3 }  \Bigr)
\end{equation*}
from which we deduce the following identity
\begin{equation*}
T_a ^{(1)} (T_s) = T_s \, \rho (T_s) = T_s \, L^{-1} \Bigl( \frac{ K_{\text{ph}} }{ T_s ^3 }  \Bigr) .
\end{equation*}
We compute the second derivative of $T^{(1)}_a$ with respect to $T_s$ and we get
$$
\frac{d \, T_a ^{(1)}}{dT_s} (T_s) 
= L^{-1} \Bigl( \frac{ K_{\text{ph}} }{ T_s ^3 }  \Bigr) + T_s \, \frac{d }{d T_s} \Bigl( L^{-1} \Bigl( \frac{ K_{\text{ph}} }{ T_s ^3 }  \Bigr) \Bigr) ,$$
and
$$
\frac{d^2 \, T_a ^{(1)}}{dT_s ^2} (T_s) 
= 2 \frac{d }{d T_s} \Bigl( L^{-1} \Bigl( \frac{ K_{\text{ph}} }{ T_s ^3 }  \Bigr) \Bigr) + T_s \, \frac{d^2 }{d T_s ^2} \Bigl( L^{-1} \Bigl( \frac{ K_{\text{ph}} }{ T_s ^3 }  \Bigr) \Bigr) .$$
The derivation rule of an inverse function gives
$$
\frac{d }{d T_s} \Bigl( L^{-1} \Bigl( \frac{ K_{\text{ph}} }{ T_s ^3 }  \Bigr) \Bigr) 
= \frac{-3 K_{\text{ph}}}{T_s ^4} \,
\frac{1}{L' \Bigl( L^{-1} \Bigl( \frac{ K_{\text{ph}} }{ T_s ^3 }  \Bigr) \Bigr)} ,
$$
and therefore we obtain
\begin{equation*}
\frac{d^2 }{d T_s ^2} \Bigl( L^{-1} \Bigl( \frac{ K_{\text{ph}} }{ T_s ^3 }  \Bigr) \Bigr)
= \frac{12 K_{\text{ph}}}{T_s ^5} \,
\frac{1}{L' \Bigl( L^{-1} \Bigl( \frac{ K_{\text{ph}} }{ T_s ^3 }  \Bigr) \Bigr)}
+ \frac{3 K_{\text{ph}}}{T_s ^4} \, \frac{\frac{d }{d T_s} \Bigl[ L' \Bigl( L^{-1} \Bigl( \frac{ K_{\text{ph}} }{ T_s ^3 }  \Bigr) \Bigr) \Bigr] }{\Bigl[ L' \Bigl( L^{-1} \Bigl( \frac{ K_{\text{ph}} }{ T_s ^3 }  \Bigr) \Bigr)\Bigr] ^2} .    \end{equation*}
Using the latter expression inside the second derivative of $T^{(1)}_a$ we get
\begin{equation*}
\begin{split}
&\frac{d^2 \, T_a ^{(1)}}{dT_s ^2} (T_s)  
= \frac{-6 K_{\text{ph}}}{T_s ^4} \,
\frac{1}{L' \Bigl( L^{-1} \Bigl( \frac{ K_{\text{ph}} }{ T_s ^3 }  \Bigr) \Bigr)}
\\
&\qquad+ \frac{12 K_{\text{ph}}}{T_s ^4} \,
\frac{1}{L' \Bigl( L^{-1} \Bigl( \frac{ K_{\text{ph}} }{ T_s ^3 }  \Bigr) \Bigr)}
+ \frac{3 K_{\text{ph}}}{T_s ^3} \, \frac{\frac{d }{d T_s} \Bigl[ L' \Bigl( L^{-1} \Bigl( \frac{ K_{\text{ph}} }{ T_s ^3 }  \Bigr) \Bigr) \Bigr] }{\Bigl[ L' \Bigl( L^{-1} \Bigl( \frac{ K_{\text{ph}} }{ T_s ^3 }  \Bigr) \Bigr)\Bigr] ^2} 
\\
&\quad= \frac{6 K_{\text{ph}}}{T_s ^4} \,
\frac{1}{L' \Bigl( L^{-1} \Bigl( \frac{ K_{\text{ph}} }{ T_s ^3 }  \Bigr) \Bigr)}
\, \left\{ 1+ \frac{T_s}{2} \frac{\frac{d }{d T_s} \Bigl[ L' \Bigl( L^{-1} \Bigl( \frac{ K_{\text{ph}} }{ T_s ^3 }  \Bigr) \Bigr) \Bigr] }{\Bigl[ L' \Bigl( L^{-1} \Bigl( \frac{ K_{\text{ph}} }{ T_s ^3 }  \Bigr) \Bigr)\Bigr] }\right\} .
\end{split}
\end{equation*}
We observe that
\begin{equation*}
\begin{split}
1+ \frac{T_s}{2} \frac{\frac{d }{d T_s} \Bigl[ L' \Bigl( L^{-1} \Bigl( \frac{ K_{\text{ph}} }{ T_s ^3 }  \Bigr) \Bigr) \Bigr] }{\Bigl[ L' \Bigl( L^{-1} \Bigl( \frac{ K_{\text{ph}} }{ T_s ^3 }  \Bigr) \Bigr)\Bigr] }&= 1+ \frac{T_s}{2} \, \frac{d }{d T_s} \Bigl\{ \ln  L' \Bigl( L^{-1} \Bigl( \frac{ K_{\text{ph}} }{ T_s ^3 }  \Bigr) \Bigr)    \Bigr \}
\\
&= 1+ \frac{T_s}{2} \, \frac{d }{d T_s} \Bigl\{ \ln \Bigl( L' (\rho (T_s)) \Bigr)    \Bigr \} .
\end{split}
\end{equation*}
On the other hand we have that
\begin{equation*}
\ln \left( L' (\rho) \right) = \ln \left( \frac{8\rho ^3 - 6 \rho ^4 -1}{(1-\rho)^2} \right)= \ln \left(8\rho ^3 - 6 \rho ^4 -1 \right) - 2 \ln \left(1-\rho \right) ,
\end{equation*}
hence
\begin{equation*}
\begin{split}
\frac{d }{d T_s} \left\{ \ln \left( L' (\rho (T_s)) \right)    \right \}
&= \frac{d }{d T_s} \left\{ \ln \left(  8\rho (T_s) ^3 - 6 \rho (T_s) ^4 -1 \right) - 2 \ln \left( 1-\rho (T_s) \right) \right \}
\\
&= \left[ \frac{24\rho (T_s)^2 \, (1-\rho (T_s))}{8\rho (T_s) ^3 - 6 \rho (T_s) ^4 -1}+ \frac{2}{1-\rho (T_s)} \right] \, \rho ' (T_s)
\\
&= \left[ \frac{24\rho (T_s)^2 \, (1-\rho (T_s))}{8\rho (T_s) ^3 - 6 \rho (T_s) ^4 -1}+ \frac{2}{1-\rho (T_s)} \right] \,
\frac{-3 K_{\text{ph}}}{T_s ^4} \, \frac{1}{L' (\rho (T_s)) } .
\end{split}
\end{equation*}
Using the latter expression, we have that
\begin{equation*}
\begin{split}
& 1+ \frac{T_s}{2} \frac{\frac{d }{d T_s} \Bigl[ L' \Bigl( L^{-1} \Bigl( \frac{ K_{\text{ph}} }{ T_s ^3 }  \Bigr) \Bigr) \Bigr] }{\Bigl[ L' \Bigl( L^{-1} \Bigl( \frac{ K_{\text{ph}} }{ T_s ^3 }  \Bigr) \Bigr)\Bigr] }=
\\
&\,\,= 1 - \frac{3 K_{\text{ph}}}{2 T_s ^3}\, \left[ \frac{24\rho (T_s)^2 \, (1-\rho (T_s))}{8\rho (T_s) ^3 - 6 \rho (T_s) ^4 -1}+ \frac{2}{1-\rho (T_s)} \right] \frac{1}{L' (\rho (T_s)) }
\\
&\,\,= 1 - \frac{3 }{2 }\, \left[ \frac{24\rho (T_s)^2 \, (1-\rho (T_s))}{8\rho (T_s) ^3 - 6 \rho (T_s) ^4 -1}+ \frac{2}{1-\rho (T_s)} \right] \cdot \frac{L (\rho (T_s))}{L' (\rho (T_s)) } 
\\
&\,\,= 1 - \frac{3 }{2 }\, \left[ \frac{24\rho (T_s)^2 \, (1-\rho (T_s))}{8\rho (T_s) ^3 - 6 \rho (T_s) ^4 -1}+ \frac{2}{1-\rho (T_s)} \right] \frac{(2 \rho (T_s) ^4 -1) (1-\rho (T_s))}{8 \rho (T_s) ^3 - 6 \rho (T_s) ^4 - 1 } .
\end{split}
\end{equation*}
Let us introduce the following function on $[2^{-1/4},1)$
\begin{equation}
\label{def-N}
\begin{split}
N: \rho &\mapsto 1 - \frac{3 }{2 }\, \Bigl[ \frac{24\rho ^2 \, (1-\rho )}{8\rho  ^3 - 6 \rho  ^4 -1}+ \frac{2}{1-\rho } \Bigr] \, \frac{(2 \rho ^4 -1) (1-\rho )}{8 \rho ^3 - 6 \rho ^4 - 1 } 
\\
&\,\,= 1 - \frac{3 }{2 }\, \Bigl[ \frac{24\rho ^2 \, (2 \rho ^4 -1) \, (1-\rho )^2}{(8\rho  ^3 - 6 \rho  ^4 -1)^2}+ \frac{2 (2 \rho ^4 -1)}{8 \rho ^3 - 6 \rho ^4 - 1 } \Bigr] .
\end{split}
\end{equation}
This gives the following expressions of $\Phi _1 ''$
\begin{equation}
\label{phi_1''}
\Phi _1 '' (T_s) = -\frac{3 \lambda K_{\text{ph}}}{T_s ^4 \, L' (\rho (T_s))} \, N (\rho (T_s)) .
\end{equation}
We use the above identity for $\Phi ''$: thanks to \eqref{def-L-K} we have that $L(\rho (T_s))= \frac{K_{ph}}{T_s ^3}$ and therefore we get
\begin{equation*}
\begin{split}
\Phi '' (T_s)& = - \frac{\lambda}{2} \frac{d^2 \, T_a ^{(1)}}{dT_s ^2} (T_s) + 12 \sigma _B (1-\frac{\varepsilon _a}{2}) T_s ^2
\\
&= - \frac{3 \lambda }{K_{\text{ph}} ^{1/3}} \, \frac{L(\rho (T_s)) ^{4/3} }{L' (\rho (T_s))} \,
N(\rho (T_s))
+ 12 \sigma _B (1-\frac{\varepsilon _a}{2}) \Bigl( \frac{ K_{\text{ph}} }{ L(\rho (T_s)) } \Bigr) ^{2/3} 
\\
&= 3 \lambda \frac{(1-\rho (T_s)) ^{2/3} \, (2 \rho (T_s) ^4 -1) ^{4/3} }{ K_{\text{ph}} ^{1/3}}\\
&\qquad
\cdot\Bigl[ -\frac{N(\rho (T_s))}{8 \rho (T_s) ^3 - 6 \rho (T_s) ^4 - 1} 
+ 4 (\frac{1}{\varepsilon _a} -\frac{1}{2}) \, \frac{1}{(2\rho (T_s) ^4 -1)^2} \Bigr] .
\end{split}
\end{equation*}
Let us define
$$ \rho _s := \rho (T_s) = L^{-1} \Bigl( \frac{ K_{\text{ph}} }{ T_s ^3 }  \Bigr) ,$$
and the following function on $(2^{-1/4}, 1)$
\begin{equation}
\label{def-N*}
N^*: \rho\mapsto \Bigl[ - \frac{N(\rho )}{8 \rho  ^3 - 6 \rho  ^4 - 1} 
+ 4 (\frac{1}{\varepsilon _a} -\frac{1}{2}) \, \frac{1}{(2\rho  ^4 -1)^2} \Bigr] .
\end{equation}
Thus, we rewrite $\Phi''$ as follows
\begin{equation}
\label{phi''}
\Phi  '' (T_s) = \Bigl( \frac{3 \lambda}{K_{\text{ph}} ^{1/3}} \, (1-\rho _s) ^{2/3} \, (2 \rho _s ^4 -1) ^{4/3} \Bigr) \, N^* (\rho _s) .
\end{equation}

\subsection{Proof of Lemma \ref{lem-nbre-eq-1*conv-phys}} \hfill

Our aim is to prove that
\begin{equation}
\label{phi''0}
\varepsilon _a < 1 \quad \implies \quad \begin{cases} N^* (\rho) >0 \\
\forall \rho \in (2^{-1/4}, 1) .
\end{cases}
\end{equation}
We note that
\begin{equation*}
\begin{split}
&N^* (\rho)
= 4 (\frac{1}{\varepsilon _a} -\frac{1}{2}) \, \frac{1}{(2\rho  ^4 -1)^2}
- \frac{1}{8 \rho  ^3 - 6 \rho  ^4 - 1} \\
&\qquad\,\,\cdot\Bigl( 1 - \frac{3 }{2 }\, \Bigl[ \frac{24\rho ^2 \, (2 \rho ^4 -1) \, (1-\rho )^2}{(8\rho  ^3 - 6 \rho  ^4 -1)^2}+ \frac{2 (2 \rho ^4 -1)}{8 \rho ^3 - 6 \rho ^4 - 1 } \Bigr] \Bigr)
\\
&= 4 (\frac{1}{\varepsilon _a} -\frac{1}{2}) \, \frac{1}{(2\rho  ^4 -1)^2} - \frac{1}{8 \rho  ^3 - 6 \rho  ^4 - 1}+ \frac{3/2}{8 \rho  ^3 - 6 \rho  ^4 - 1}
\\
&\qquad\qquad\qquad \cdot\Bigl[ \frac{24\rho ^2 \, (2 \rho ^4 -1) \, (1-\rho )^2}{(8\rho  ^3 - 6 \rho  ^4 -1)^2}+ \frac{2 (2 \rho ^4 -1)}{8 \rho ^3 - 6 \rho ^4 - 1 } \Bigr] .
\end{split}
\end{equation*}
The function $\rho \mapsto 8 \rho ^3 - 6 \rho ^4 - 1$ is increasing and positive on $(2^{-1/4}, 1)$.
Therefore, for all $\rho \in (2^{-1/4}, 1)$, we have
\begin{multline*}
N^* (\rho)
\geq 4 (\frac{1}{\varepsilon _a} -\frac{1}{2}) \, \frac{1}{(2\rho  ^4 -1)^2} - \frac{1}{8 \rho  ^3 - 6 \rho  ^4 - 1}
\\
= \frac{1}{2\rho  ^4 -1} \Bigl( 4 (\frac{1}{\varepsilon _a} -\frac{1}{2}) \, \frac{1}{2\rho  ^4 -1} - \frac{2\rho  ^4 -1}{8 \rho  ^3 - 6 \rho  ^4 - 1} \Bigr) .
\end{multline*}
We observe that the function $\rho \mapsto \frac{1}{2\rho  ^4 -1}$ is decreasing on $(2^{-1/4}, 1)$ and is equal to $1$ for $\rho =1$. Whereas the function
$$ \tilde N: \rho \in (2^{-1/4}, 1) \mapsto \tilde N (\rho) := \frac{2 \rho ^4 -1}{8 \rho  ^3 - 6 \rho  ^4 - 1} $$
is strictly increasing on $(2^{-1/4}, 1)$.
Indeed,
\begin{equation*}
\begin{split}
\tilde N ' (\rho)&= \frac{8\rho ^3  (8 \rho  ^3 - 6 \rho  ^4 - 1) - (2 \rho ^4 -1)(24 \rho ^2 - 24 \rho ^3)}{(8 \rho  ^3 - 6 \rho  ^4 - 1)^2}
\\
&= \frac{8\rho ^2}{(8 \rho  ^3 - 6 \rho  ^4 - 1)^2}\Bigl( (8 \rho  ^4 - 6 \rho  ^5 - \rho ) - (2 \rho ^4 -1)(3 - 3 \rho)\Bigr) \\
&=\frac{8\rho ^2}{(8 \rho  ^3 - 6 \rho  ^4 - 1)^2}\Bigl( 2 \rho ^4 - 4 \rho +3\Bigr)\geq0\quad\text{on}\quad (2^{-1/4}, 1).
\end{split}
\end{equation*}
Therefore $\tilde N$ is increasing on $(2^{-1/4}, 1)$.
Hence, 
$$ (2 \rho ^4 -1) N^* (\rho) \geq 4 (\frac{1}{\varepsilon _a} -\frac{1}{2}) - 1 
= \frac{4-3\varepsilon _a}{\varepsilon _a}.$$
Therefore, if $\varepsilon _a \leq1$ (actually the same holds for $\varepsilon _a<\frac{4}{3}$), $N^*$ is positive on $(2^{-1/4}, 1)$ and so $\Phi ''$ is positive on $(0,+\infty)$, thanks to \eqref{phi''}.
This was the main part of Lemma \ref{lem-nbre-eq-1*conv-phys}.

Now, if $\varepsilon _a \leq1$, $\Phi$ is strictly convex and strictly increasing on $[0,+\infty)$ (see the proof of Lemma \ref{lem-nbre-eq-1}). Thus, equation $\Phi (T_s) = q \beta _s (T_s)$ can have at most one solution on $[0,T_{s,-}]$, one on $[T_{s,+}, +\infty)$ and
two on $[T_{s,-}, T_{s,+}]$. If there are two solutions on $[T_{s,-}, T_{s,+}]$, the strict convexity of $\Phi$
implies that there cannot be other solutions on $[T_{s,+}, +\infty)$. Therefore, in this case, there are at most three solutions.
This concludes the proof of Lemma \ref{lem-nbre-eq-1*conv-phys}. \qed

\subsection{On Remark \ref{lem-nbre-eq-1*conv-simpl}}\hfill

\begin{figure}[h!]
        \centering
       \includegraphics[scale=.25]{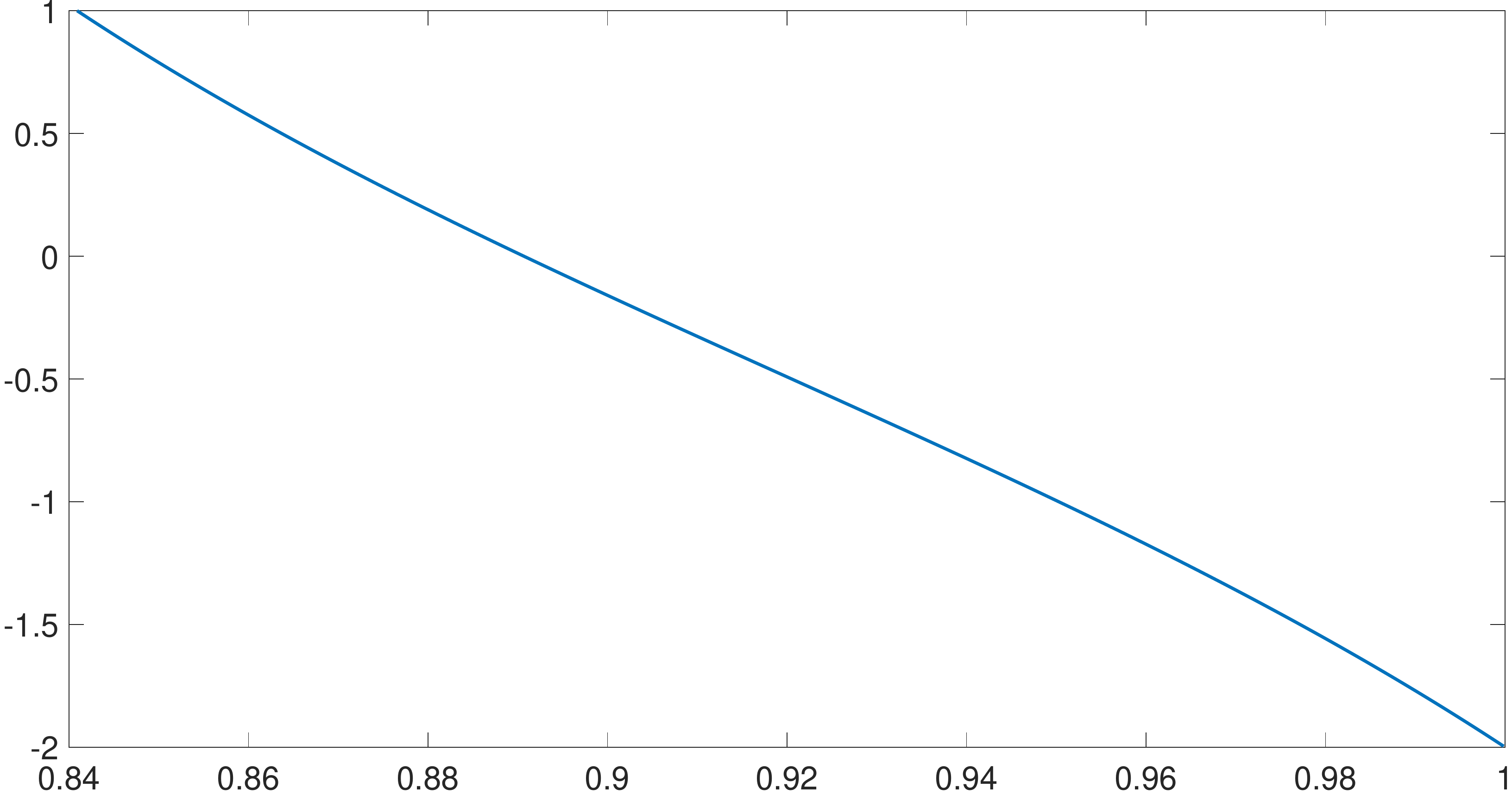}
        \caption{Graph of function $N$, defined in \eqref{def-N}.}\label{AC-figure2.1}
\end{figure}

For $\varepsilon _a \in [1,2)$, we used numerical tools have an idea the convexity of $\Phi$. 
First we note that $N$ is equal to $1$ for  $\rho = 2^{-1/4}$, converges to $-2$ for $\rho \to 1^-$, hence has at least one zero. By plotting such function we observe that it is decreasing on
$[2^{-1/4},1)$ and it has a unique zero, $\rho_0$, whose value is approximately $0,89$ (see Figure \ref{AC-figure2.1}).

This gives information on the convexity of $T_a ^{(1)}$ and of $\Phi _1$.

\begin{figure}[h!]
\centering
\subfloat[]{
\includegraphics[scale=.26]{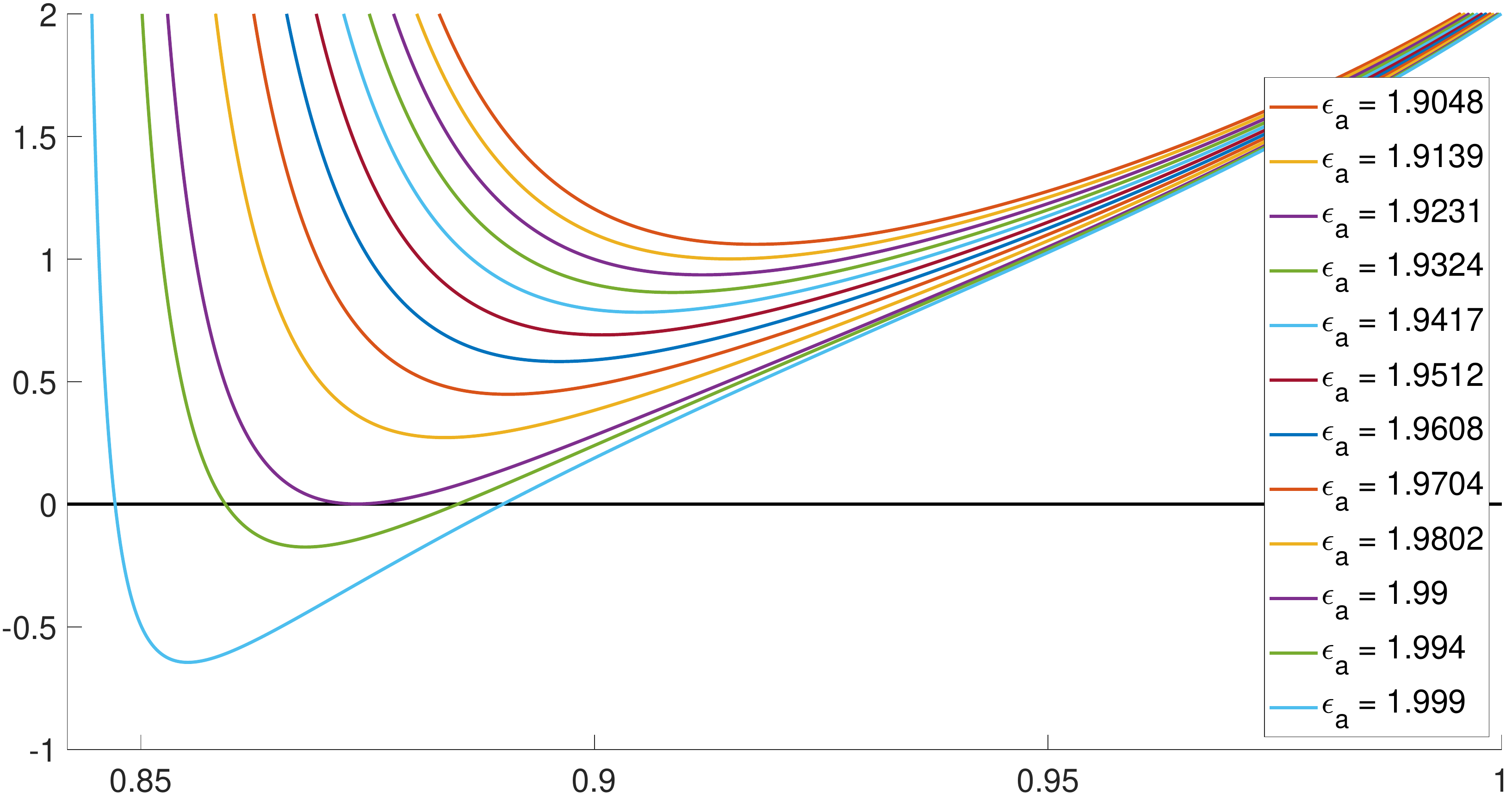}
}
\caption{The graph $N^*$ for $\varepsilon _a$ close to $2$.}
\label{fig:graphN*-e_aclose2C}
\end{figure}

Concerning $N^*$, we observe that that there exists a unique value $\varepsilon _{a,0} \in (1.99, 1.991)$ such that:
\begin{itemize}
\item $N^*$ is positive on $(2^{-1/4}, 1)$ as long as $\varepsilon _a <\varepsilon _{a,0}$,
\item $N^*$ has exactly two zeros in $(2^{-1/4}, 1)$ for all $\varepsilon _a \in (\varepsilon _{a,0}, 2)$, it is negative between these two zeros and positive elsewhere,
\end{itemize}
see see Figure \ref{fig:graphN*-e_aclose2C} . Therefore, following the numerical results, the convexity of $\Phi$ is described in Remark \ref{lem-nbre-eq-1*conv-simpl}.
Concerning the number of equilibrium points, it remains to study the case $\varepsilon _a \in (\varepsilon _{a,0}, 2)$, where the function $\Phi$ (from numerical tests) is respectively convex, concave and convex over the interval $[0,+\infty)$. Let us assume that such behaviour is satisfied. Then, since $\Phi$ is strictly increasing, equation $\Phi (T_s) = q \beta _s (T_s)$ has at most one solution in $[0,T_{s,-}]$ and at most one solution in $[T_{s,+}, +\infty)$. If we have four solutions of the equation
$$\Phi(T_s)=q\beta_s(T_s)$$
in the interval $[T_{s,-}, T_{s,+}]$ then, by Rolle's Theorem,
$$\Phi'(T_s)=q\frac{\beta_{s,+}-\beta_{s,-}}{T_{s,+}-T_{s,-}}$$
admits three solutions. Applying again Rolle's Theorem we deduce that $\Phi ''$ vanishes at least twice. And since $\Phi ''$ has at most two zeros, we conclude that there cannot be five solutions of
$$\Phi(T_s)=q\beta_s(T_s)$$
in the interval $[T_{s,-}, T_{s,+}]$.

Therefore, the equation of the equilibrium points can have at most six solutions over all $[0,+\infty)$. In this case there would be one equilibrium point in $[0,T_{s,-}]$, one in $[T_{s,+}, +\infty)$, and four in $(T_{s,-}, T_{s,+})$. Let us denote them by $T_{s,i}$, with $i=1,\dots, 6$. As explained above, there exist three values $\tilde T_{s,j}\in[T_{s,-},T_{s,+}]$, $j=1,2,3$, such that 
$$ \Phi '( \tilde T_{s,j}) = q\beta ' (\tilde T_{s,j})$$
and $ T_{s,2} < \tilde T_{s,1} <  T_{s,3} <\tilde T_{s,2} < T_{s,4} < \tilde T_{s,3} < T_{s,5} $. However, $\Phi$ is strictly convex on $[\tilde T_{s,3}, +\infty)$
and $\Phi '(  T_{s,5}) > q\beta ' ( T_{s,5})$ thus the existence of $T_{s,6}$. Therefore, we conclude that there exist at most five equilibrium points. \qed

\begin{figure}[h!]
\centering
\epsfig{figure=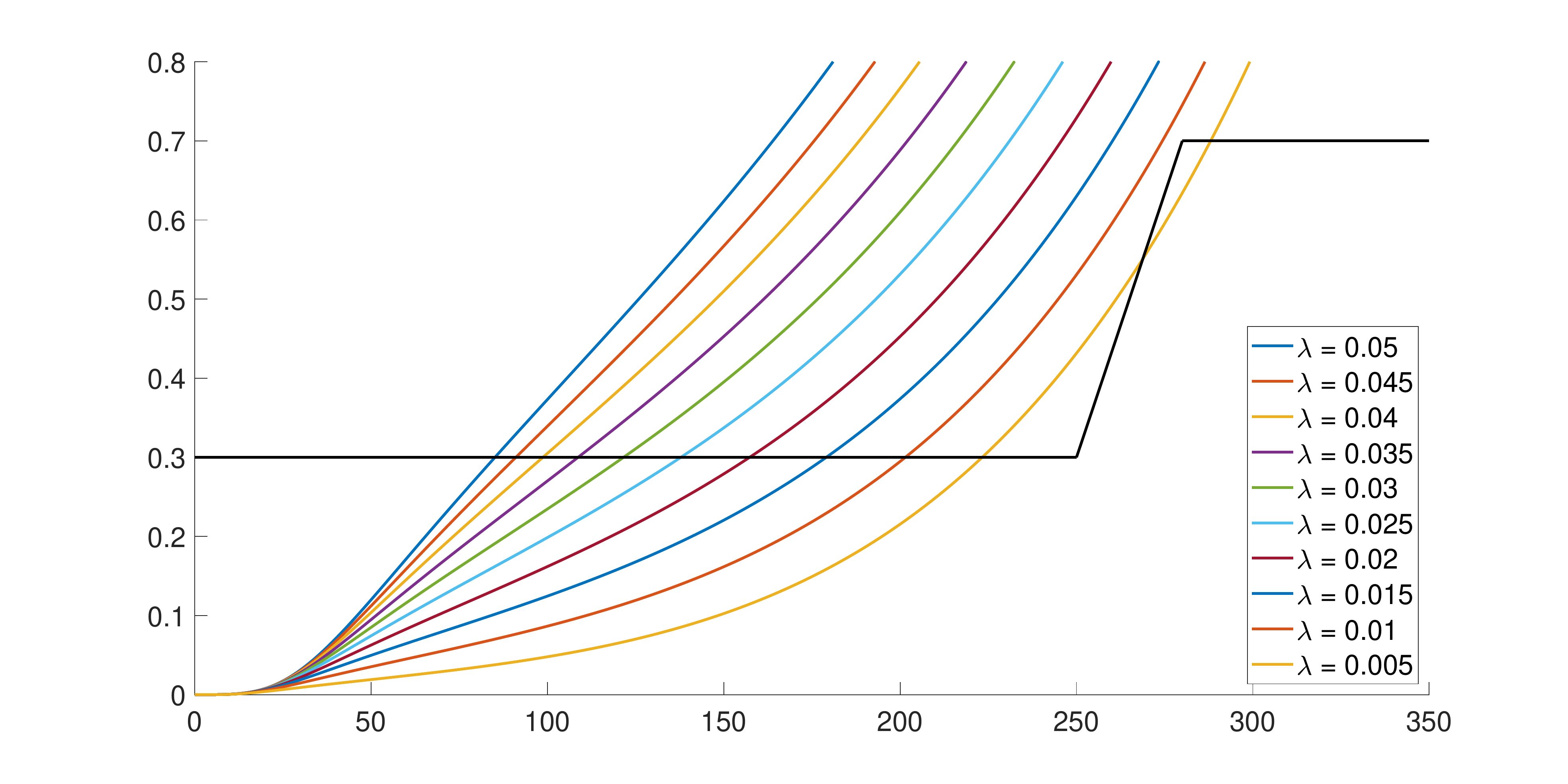,width=8cm}
\caption{Intersection of the graph of $\Phi$ with $\eps_a=1.997$ and $\beta _s$, with ten values of $\lambda$ in $\{0,005,...,0.05\}$.}
\label{PhiBeta}
\end{figure}

In the picture that follows we have plot the graph of $\Phi$, for $\eps_a=1.997$ and ten values of $\lambda$ from $0.005$ to $0.05$, and the graph of $\beta _s$. Even for such a big value of $\eps_a>\eps_{a,0}$ the drawing suggests that for physical relevant parameters of $\beta_s$ there are at most three equilibrium points of our system (see Figure \ref{PhiBeta}).

\qed



\end{document}